\documentclass{article}

\usepackage{amssymb}
\usepackage{amsmath}
\usepackage{arxiv}
\usepackage[utf8]{inputenc}
\usepackage[T1]{fontenc}
\usepackage{url}
\usepackage{booktabs}
\usepackage{amsfonts}
\usepackage{nicefrac}
\usepackage{microtype}
\usepackage{color}
\usepackage{lipsum}
\usepackage{amsthm}
\usepackage{indentfirst}
\usepackage{graphicx}
\usepackage{epstopdf}
\usepackage{fourier}
\usepackage{bm}
\usepackage{mathtools}
\usepackage{latexsym,enumerate}
\usepackage{multicol}
\usepackage[skip=2pt]{caption}
\usepackage[font=small,skip=0pt]{subcaption}
\usepackage{array}
\usepackage{ascii}
\usepackage{algorithm,algorithmic}
\usepackage{mathtools}
\usepackage{arydshln}
\usepackage[pdfstartview=FitH, CJKbookmarks=true,
bookmarksnumbered=true, bookmarksopen=true,
colorlinks,
linkcolor=green,
anchorcolor=blue, citecolor=blue
]{hyperref}

\usepackage{ulem}
\usepackage{multirow}
\newcommand{\comment}[1]{}

\newcommand{\BEA}{\begin{eqnarray}}
	\newcommand{\EEA}{\end{eqnarray}}

\newtheorem{thm}{Theorem}[section]

\newtheorem{example}[thm]{Example}

\newtheorem{rem}[thm]{Remark}

\newcommand{\PreserveBackslash}[1]{\let\temp=\\#1\let\\=\temp}
\newcolumntype{C}[1]{>{\PreserveBackslash\centering}p{#1}}
\newcolumntype{R}[1]{>{\PreserveBackslash\raggedleft}p{#1}}
\newcolumntype{L}[1]{>{\PreserveBackslash\raggedright}p{#1}}
\newcommand{\stkout}[1]{\ifmmode\text{\sout{\ensuremath{#1}}}\else\sout{#1}\fi}

\begin{document}
	
	\title{RBF-Generated Finite Difference Method Coupled with Quadratic
		Programming for Solving PDEs on Surfaces with Derivative Boundary Conditions}
	\author{ Peng Chen \\
		School of Information Science and Technology, ShanghaiTech University,
		Shanghai 201210, China \\
		\texttt{chenpeng2024@shanghaitech.edu.cn}
		\And Shixiao Willing Jiang \footnotemark[1] \\
		Institute of Mathematical Sciences, ShanghaiTech University, Shanghai
		201210, China \\
		\texttt{jiangshx@shanghaitech.edu.cn}
		\And Rongji Li \thanks{Corresponding authors.} \\
		School of Information Science and Technology, ShanghaiTech University,
		Shanghai 201210, China\\
		\texttt{lirj2022@shanghaitech.edu.cn}
		\And Qile Yan \\
		School of Mathematics, University of Minnesota, 206 Church St SE,
		Minneapolis, MN 55455, USA \\
		\texttt{yan00082@umn.edu}
	}
	\date{\today }
	\maketitle
	
	\begin{abstract}
		
		Derivative boundary conditions introduce challenges for mesh-free discretizations of PDEs on surfaces, especially when the domain is represented by randomly sampled point clouds. The recently developed two-step tangent-space RBF-generated finite difference (RBF-FD) method provides high accuracy on closed surfaces. However, it may lose stability when applied directly to surface PDEs with derivative boundary conditions. To enhance numerical stability, we develop a mesh-free method that couples the two-step tangent-space RBF-FD discretization with a quadratic programming (QP) procedure to stabilize the operator approximation for interior points near boundaries. For boundary points, we construct restricted nearest-neighbor stencils biased in the co-normal direction and employ a constrained quadratic program to approximate outward co-normal derivatives. The resulting method avoids using ghost points and does not require quasi-uniform node distributions. We validate the approach on elliptic problems, eigenvalue problems, time-dependent diffusion equations, and elliptic interface problems on surfaces with boundary. Numerical experiments demonstrate stable performance and high-order accuracy across a variety of surfaces.
	\end{abstract}
	
	
	\lhead{} \newpage
	
	\section{Introduction}
	
	Partial differential equations on surfaces or manifolds arise across physics \cite{rauter2018finite}, biology \cite{elliott2010modeling}, and image processing \cite{tian2009segmentation}. In many applications, the surface or manifold has a boundary, and Neumann or Robin conditions naturally encode fluxes, contact-line physics, and exchanges with a surrounding bulk medium. For instance, reaction–diffusion systems on evolving biological surfaces \cite{barreira2011surface}  model pattern formation in developmental biology and cell-scale morphogenesis, with boundary conditions representing no-flux constraints or regulated transport across edges. Closely related bulk–surface reaction–diffusion models \cite{alphonse2018coupled} used in cell signaling impose Robin boundary conditions that couple bulk diffusion to surface reactions. More broadly, when the surface boundary has direct physical meaning, derivative boundary conditions can also encode geometric constraints and force balance. Examples include surface diffusion flow with contact-angle boundary conditions \cite{asai2014self} and surfactant transport near moving contact lines \cite{lai2010numerical},  where an accurate boundary treatment is essential for capturing the correct dynamics.
	
	Many numerical methods have been developed for approximating solutions of PDEs on surfaces without or with boundaries. One major class consists of mesh-based methods. Given a triangular mesh approximating the domain, the surface finite element method (FEM) \cite{dziuk2013finite} provides a robust and efficient framework for discretizing PDEs directly on the mesh. Another major class consists of mesh-free methods, which are particularly effective when the domain is represented only by point cloud data, as in the setting considered in this paper. Examples of such mesh-free methods include radial basis function (RBF) methods \cite{flyer2009radial,Fuselier2009Stability,piret2012orthogonal}, RBF-generated finite differences (RBF-FD) \cite{shankar2015radial,flyer2016role,lehto2017radial,petras2018rbf}, generalized moving least squares (GMLS) \cite{liang2013solving,gross2020meshfree,li2024generalized,hangelbroek2024generalized}, generalized finite difference methods (GFDM) \cite{suchde2021meshfree,halada2025overview,jiang2024generalized}, and graph-based approaches \cite{li2017point,gh2019}, among others.
	
	Nonetheless, surfaces or manifolds with boundaries present additional challenges for kernel-based mesh-free discretizations. Because neighborhoods near the boundary are necessarily one-sided, constructing accurate and stable local discretizations becomes more difficult. In strong-form localized mesh-free collocation methods, this difficulty is particularly pronounced for Neumann boundary conditions, which are known to cause loss of robustness and ill-conditioning \cite{liu2006stabilized}. To address these challenges, a number of strategies have been proposed. One approach is the use of ghost points, where additional points are introduced outside the domain to impose derivative boundary conditions more accurately. Ghost points were introduced for global RBF collocation methods in Euclidean domains in \cite{lin2020radial,chen2020novel}, and related ideas were later extended to surfaces in \cite{jiang2023ghost,yan2023kernel}. However, ghost-point constructions on surfaces are substantially more involved, since they require both a reliable extension of the underlying geometry and accurate estimation of boundary normals or co-normals from point cloud data. Another approach is to improve stability through weak formulations. For example, \cite{TominecLarssonHeryudono2021} proposed an RBF-FD method in a least-squares setting and showed that it is significantly more robust than RBF-FD collocation method without ghost points. For unfitted node sets, additional constraint-based treatments of boundary conditions have also been proposed, including Lagrange multiplier formulations \cite{nielsen2025high}. We note that quasi-uniform node distributions are typically assumed in these methods.

	The main contributions of this work are listed as follows.
		{ First, we extend our two-step RBF-FD method \cite{li2025two} to solve boundary value problems (BVPs) on surfaces represented by randomly sampled point clouds, thereby relaxing the need for quasi-uniform node distributions.}
		{ Second, we introduce a quadratic-programming-based approach to compute stable operator weights for interior points near the boundary, addressing the ill-conditioning caused by one-sided stencils. This approach ensures stability while preserving the required consistency constraints.}
		Third, a stable discretization of the derivative boundary condition is essential. Here we construct stencils for boundary points using restricted $K$-nearest neighbors, selected based on a weighted quadratic norm biased in the co-normal direction.
		Together with an auto-tuning strategy for $K$, we then apply a constraint quadratic optimization formulation to approximate the outward co-normal derivative at boundary points.
		Overall, these strategies yield a stable and robust estimator and are shown to be effective for several classes of surface PDE problems, including eigenvalue problems, time-dependent diffusion equations, and elliptic interface problems. Our resulting mesh-free approach avoids using ghost points, as ghost point extension is nontrivial on curved surfaces, especially for unknown surfaces.

	The paper is organized as follows. In Sec.~\ref{sec:pre}, we provide a short review of the GMLS approach and the two-step tangent-space RBF-FD approach for approximating the Laplace-Beltrami operator. We then present a new RBF-FD framework coupled with quadratic-programming-based stabilization for interior points near the boundary of the surface.
		In Sec.~\ref{sec:bound}, we present the discretization of outward co-normal derivatives at boundary points, and then describe how the interior and boundary discretizations are assembled for solving boundary value problems on surfaces.
		In Sec.~\ref{sec:numerical}, we provide supporting numerical examples for solving boundary value surface PDEs, particularly eigenvalue problems, linear diffusion equations, and elliptic interface problems. In Sec. \ref{sec:conclusion}, we conclude the paper with a summary and some open problems.

	\section{Approximation of the Laplace--Beltrami operator on surfaces with boundaries}
	
	\label{sec:pre} Let $M$ be a two-dimensional surface embedded in $\mathbb{R}%
	^{3}$ with boundary $\partial M$. Let $T_{\mathbf{x}}M$ be the tangent space
	at an arbitrary point $\mathbf{x}\in M$ and let $\left\{ \boldsymbol{t}%
	_{1}\left( \mathbf{x}\right) ,\boldsymbol{t}_{2}\left( \mathbf{x}\right)
	\right\} $ be an orthonormal basis of the tangent space. Let $u:M\rightarrow
	\mathbb{R}$ be a smooth function. Given a set of $N$ (distinct) randomly
	sampled data points $\mathbf{X}_{M}=\left\{ \mathbf{x}_{i}\right\}
	_{i=1}^{N} $ on $M$ and function values $\mathbf{u}:=(u(\mathbf{x}%
	_{1}),\ldots ,u(\mathbf{x}_{N}))^{\top }$ on $\mathbf{X}_{M}$. Among these $%
	N $\ points, let the first $N_{I}$ points be interior points $\mathbf{X}%
	_{I}=\left\{ \mathbf{x}_{i}\right\} _{i=1}^{N_{I}}\subset \mathrm{Int}(M)$\
	and the last $N_{B}=N-N_{I}$ points be boundary points $\mathbf{X}%
	_{B}=\left\{ \mathbf{x}_{j}\right\} _{j=N_{I}+1}^{N}\subset \partial M$. Hence we have
	$\mathbf{X}_{M} = \mathbf{X}_{I} \cup \mathbf{X}_{B}$. In
	Sec. \ref{sec:2_1} and Sec. \ref{sec:rbffd}, we first briefly review the
	GMLS and the two-step RBF-FD methods for approximating functions on surfaces
	from the data points $\mathbf{X}_{M}$ and the function values $\mathbf{u}$.
	In Sec. \ref{sec:rbfaway}, we review RBF-FD methods for approximating the
	Laplace--Beltrami operator at most interior points. In Sec. \ref{sec:rbfnear}%
	, we introduce a new quadratic programming approach to enhance the stability
	of the Laplace--Beltrami operator at the remaining interior points, which
	typically lie near the boundary of a surface.
	
	\subsection{Review of GMLS approaches on surfaces}
	
	\label{sec:2_1} The GMLS approach, which utilizes a set of local polynomials
	for function approximation, has been successfully used to solve PDEs on
	surfaces and manifolds \cite%
	{liang2013solving,gross2020meshfree,mohammadi2021divergence,jones2023generalized,jiang2024generalized}%
	. For an interior base point $\mathbf{x}_{0}\in \mathbf{X}_{I}$, we denote
	its $K$-nearest neighbors in $\mathbf{X}_{M}$ by $S_{\mathbf{x}_{0}}=\left\{
	\mathbf{x}_{0,k}\in \mathbf{X}_{M}\right\} _{k=1}^{K}$, which is generally
	referred to as a "stencil". By definition, we have $\mathbf{x}_{0,1}=\mathbf{%
		x}_{0}$ for the base point. Let $\mathbb{P}_{\mathbf{x}_{0}}^{l}$ denote the
	local polynomial space with degree up to $l$ near the base point ${\mathbf{x}%
		_{0}}$, i.e.,
	\begin{equation}
		\mathbb{P}_{\mathbf{x}_{0}}^{l}=\mathrm{span}\Big\{p_{\boldsymbol{\alpha }%
		}\left( \boldsymbol{\theta }\left( \mathbf{x}\right) \right) :=\theta
		_{1}^{\alpha _{1}}\left( \mathbf{x}\right) \theta _{2}^{\alpha _{2}}\left(
		\mathbf{x}\right) \Big|\ 0\leq \left\vert \boldsymbol{\alpha }\right\vert
		=\alpha _{1}+\alpha _{2}\leq l\ \Big\},  \label{eq:thet12}
	\end{equation}%
	where $\boldsymbol{\theta }\left( \mathbf{x}\right) =(\theta _{1}(\mathbf{x}%
	),\theta _{2}(\mathbf{x}))$ is a projected Monge coordinate system for the
	local tangent space $T_{\mathbf{x}_{0}}M$ \cite%
	{monge1809application,pressley2010elementary,liang2013solving,jones2023generalized,gross2020meshfree}
	with $\theta _{i}(\mathbf{x})=\boldsymbol{t}_{i}({\mathbf{x}_{0}})\cdot (%
	\mathbf{x}-{\mathbf{x}_{0}})$ for\ $i=1,2$. Here, $\boldsymbol{\alpha }%
	=(\alpha _{1},\alpha _{2})\in \mathbb{N}^{2}$ is the nonnegative integer
	exponent for powers of variables $(\theta _{1}(\mathbf{x}),\theta _{2}(%
	\mathbf{x}))$. By definition, the dimension of the space $\mathbb{P}_{%
		\mathbf{x}_{0}}^{l}$ is $m=(l+2)(l+1)/2$.
	
	Given the $K(>m)$-nearest neighbors and the $K$ function values $\mathbf{u}_{%
		{\mathbf{x}_{0}}}:=(u({\mathbf{x}_{0,1}}),...,u(\mathbf{x}_{0,K}))^{\top }$,
	we can find $\mathcal{I}_{p}\mathbf{u}_{{\mathbf{x}_{0}}}\in \mathbb{P}_{%
		\mathbf{x}_{0}}^{l}$ to be the optimal solution of the moving least-squares
	problem:
	\begin{equation}
		\underset{\mathcal{I}_{p}\mathbf{u}_{{\mathbf{x}_{0}}}\in \mathbb{P}_{%
				\mathbf{x}_{0}}^{l}}{\min }\frac{1}{2}\sum_{k=1}^{K}\lambda _{k}\left( u({%
			\mathbf{x}_{0,k}})-(\mathcal{I}_{p}\mathbf{u}_{{\mathbf{x}_{0}}})({\mathbf{x}%
			_{0,k}})\right) ^{2}.  \label{eqn:int_LS}
	\end{equation}%
	The solution to the above least-squares problem (\ref{eqn:int_LS}) is
	represented as
	\begin{equation}
		(\mathcal{I}_{p}\mathbf{u}_{{\mathbf{x}_{0}}})(\mathbf{x}):=%
		\sum_{j=1}^{m}b_{j}p_{\boldsymbol{\alpha }(j)}\left( \boldsymbol{\theta }%
		\left( \mathbf{x}\right) \right) ,\text{ \ \ }\mathbf{x}\in M,
		\label{eqn:gmls}
	\end{equation}%
	where each coefficient $b_{j}$ corresponds to the $j$-th polynomial basis
	function $p_{\boldsymbol{\alpha }(j)}\left( \boldsymbol{\theta }\left(
	\mathbf{x}\right) \right) $. The concatenated coefficients are denoted by $%
	\mathbf{b}=(b_{1},...,b_{m})^{\top }$ and then the moving least-squares
	problem (\ref{eqn:int_LS})\ can be formulated in matrix-vector form as
	\begin{equation}
		\underset{{\mathbf{b\in }}\mathbb{R}^{m}}{\min }\frac{1}{2}\left( \mathbf{u}%
		_{\mathbf{x}_{0}}-\boldsymbol{P\mathbf{b}}\right) ^{\top }\boldsymbol{%
			\Lambda }\left( \mathbf{u}_{\mathbf{x}_{0}}-\boldsymbol{P\mathbf{b}}\right) ,
		\label{eqn:LS_optm}
	\end{equation}%
	where $\boldsymbol{\Lambda }=\mathrm{diag}(\lambda _{1},\cdots ,\lambda
	_{K})\in \mathbb{R}^{K\times K}$ is the $K\times K$\ symmetric
	positive-definite diagonal matrix and $\boldsymbol{P}$ is the $K\times m$\
	Vandermonde-type matrix with components $\boldsymbol{P}_{kj}=p_{\boldsymbol{%
			\alpha }(j)}\left( \boldsymbol{\theta }\left( \mathbf{x}{_{0,k}}\right)
	\right) $ for $1\leq k\leq K,1\leq j\leq m$. Thus, the coefficient $\mathbf{b%
	}$\ can be obtained by the normal equation:
	\begin{equation}
		\mathbf{b}=(\boldsymbol{P}^{\top }\boldsymbol{\Lambda }\boldsymbol{P})^{-1}%
		\boldsymbol{P}^{\top }\boldsymbol{\Lambda }\mathbf{u}_{{\mathbf{x}_{0}}}.
		\label{eqn:Pij}
	\end{equation}%
	The above normal equation admits a unique solution if $\boldsymbol{P}$ has
	full column rank ($\mathrm{rank}(\boldsymbol{P})=m$).
	
	\begin{rem}
		GMLS is a regression technique for approximating functions by solving the
		least-squares problem with a weighted $\ell ^{2}$-norm. In the
		literature, the diagonal element of the weight $\boldsymbol{\Lambda }$ is
		taken to be a positive weight function of the distance $\Vert {\mathbf{x}%
			_{0,k}-\mathbf{x}_{0}}\Vert $ (see e.g., \cite%
		{Wendland2005Scat,lipman2009stable,liang2013solving,gross2020meshfree,jones2023generalized,li2024generalized}%
		). Here, we used a so-called $1/K$ weight function as reported in previous
		works {\cite%
			{liang2012geometric,liang2013solving,jiang2024generalized,li2024generalized}}%
		, since the diagonal weight,
		\begin{equation}
			\boldsymbol{\Lambda }=\mathrm{diag}(\lambda _{1},\cdots ,\lambda
			_{K}),\qquad \text{with diagonal entries\ }\lambda _{k}=\left\{
			\begin{array}{ll}
				1, & \text{if }k=1 \\
				1/K, & \text{if }k=2,\ldots ,K%
			\end{array}%
			\right. ,  \label{eqn:k_weight}
		\end{equation}%
		has been shown to enhance the numerical stability of the Laplace--Beltrami operator across
		various data sets, in particular, for randomly sampled data points or
		non-quasi-uniform data points.
	\end{rem}
	
	
	
	
	\subsection{Review of two-step tangent-space RBF-FD approaches on surfaces}
	
	\label{sec:rbffd}
	
	To approximate the Laplace--Beltrami operator on surfaces, a PHS+Poly
	interpolant of the function $u$\ was first considered in \cite%
	{shaw2019radial,wright2023mgm,jones2023generalized} by using the so-called
	tangent-space RBF-FD method:%
	\begin{equation}
		(\mathcal{I}_{\phi p}\mathbf{u}_{{\mathbf{x}_{0}}})\left( \mathbf{x}\right)
		:=\sum_{k=1}^{K}a_{k}\phi \left( \Vert \boldsymbol{\theta }\left( \mathbf{x}%
		\right) -\boldsymbol{\theta }\left( \mathbf{x}_{0,k}\right) \Vert \right)
		+\sum_{j=1}^{m}b_{j}p_{\boldsymbol{\alpha }(j)}\left( \boldsymbol{\theta }%
		\left( \mathbf{x}\right) \right) ,\ \ \ \mathbf{x}\in M,  \label{eqn:If2}
	\end{equation}%
	where $\phi $ is the Polyharmonic spline (PHS) function and $\Vert \cdot
	\Vert $ is the standard Euclidean norm in $\mathbb{R}^{2}$. By definition,
	we have $\boldsymbol{\theta }\left( \mathbf{x}_{0,1}\right) =\boldsymbol{%
		\theta }\left( \mathbf{x}_{0}\right) =\mathbf{0}$ as the center of the
	projected Monge coordinate system. The method is referred to as the
	tangent-space method since both the PHS radial function $\phi (\cdot )$ and
	the polynomial $p_{\boldsymbol{\alpha }(j)}(\cdot )$ are defined over the
	projected tangent-space Monge coordinate $\boldsymbol{\theta }\left( \mathbf{%
		x}\right) $. In contrast, for RBF-FD methods in most previous studies, the
	RBF interpolant was defined over the ambient Euclidean space (see e.g., \cite%
	{shankar2015radial,lehto2017radial,petras2018rbf, alvarez2021local}). We
	point out that there are two advantages for the tangent-space method. First,
	as noted in \cite{jones2023generalized}, it simplifies the calculation of
	the derivatives of interpolation functions as also discussed in Remark 3.4
	of \cite{li2025two}. Second, the interpolant $\sum_{k=1}^{K}a_{k}\phi \left(
	\cdot \right) +\sum_{j=1}^{m}b_{j}p_{\boldsymbol{\alpha }(j)}\left( \cdot
	\right) $ is approximating the local representation of the target function,
	that is, $u\circ \boldsymbol{\theta }^{-1}$. This is natural in differential
	geometry in order to take derivatives of a function on surfaces.
	
	Here, we only consider the radial function $\phi $\ to be,
	\begin{equation}
		\text{Polyharmonic spline (PHS)}:\qquad \phi (r)=r^{2\kappa +1},
		\label{eqn:phs}
	\end{equation}%
	where the variable $r$ is the radial distance and the PHS parameter $\kappa $
	$\in \mathbb{N}$ controls the smoothness of $\phi $. For RBF-FD methods, the
	PHS function has been used together with polynomials in many studies (see
	e.g., \cite{flyer2016role,bayona2017role,bayona2019role,jones2023generalized}%
	). PHS is attractive due to its absence of the shape parameter. Some other
	common choices of the radial $\phi $ include Gaussian function \cite%
	{fasshauer2012stable}, inverse quadratic function \cite{harlim2023radial},
	inverse multiquadric function \cite{fuselier2013high} and Mat\'{e}rn class
	function \cite{fuselier2013high}. These radial basis functions contain a
	shape parameter that controls the flatness of the RBF. In general, a small
	value of the shape parameter achieves better accuracy but also results in an
	ill-conditioned interpolation matrix. When the shape parameter is properly
	chosen, the local RBF-FD shows convergence \cite%
	{shankar2015radial,lehto2017radial,petras2018rbf, alvarez2021local}.
	
	For the classical RBF-FD approach, the expansion coefficients, $\mathbf{a}%
	=\left( a_{1},a_{2},\cdots ,a_{K}\right) ^{\top }$\ and $\mathbf{b}=\left(
	b_{1},b_{2},\cdots ,b_{m}\right) ^{\top }$, are obtained by enforcing $K$
	interpolation conditions together with $m$ additional moment conditions
	ensuring uniqueness, which can be written in matrix form as%
	\begin{equation}
		\begin{bmatrix}
			\boldsymbol{\Phi } & \boldsymbol{P} \\
			\boldsymbol{P}^{\top } & \mathbf{0}%
		\end{bmatrix}%
		\begin{bmatrix}
			\mathbf{a} \\
			\mathbf{b}%
		\end{bmatrix}%
		=%
		\begin{bmatrix}
			\mathbf{u}_{\mathbf{x}_{0}} \\
			\mathbf{0}%
		\end{bmatrix}%
		,  \label{eqn:mat_PRFD}
	\end{equation}%
	where $\boldsymbol{P}\in \mathbb{R}^{K\times m}$ is defined in (\ref%
	{eqn:LS_optm}) and $\boldsymbol{\Phi }=[\phi \left( \Vert \boldsymbol{\theta
	}\left( \mathbf{x}_{0,i}\right) -\boldsymbol{\theta }\left( \mathbf{x}%
	_{0,j}\right) \Vert \right) ]_{i,j=1}^{K}\in \mathbb{R}^{K\times K}$ is the
	symmetric PHS matrix. For the PHS $\phi (r)=r^{2\kappa +1}$, the chosen
	parameter $\kappa $ should be less than or equal to the polynomial degree$\
	l $ in (\ref{eq:thet12}) \cite%
	{iske2003approximation,Wendland2005Scat,jones2023generalized} so that $%
	\boldsymbol{\Phi }$ is conditionally positive definite on the subspace
	satisfying the $m$ moment conditions in (\ref{eqn:mat_PRFD}). If the stencil
	points satisfy $\mathrm{rank}(\boldsymbol{P})=m$, then the system (\ref%
	{eqn:mat_PRFD}) is non-singular and the PHS+Poly interpolant is well-posed.
	
	We recently consider a two-step tangent-space RBF-FD method \cite{li2025two}
	which takes exactly the same approximation space as that in equation (\ref%
	{eqn:If2}), but approximates the interpolant coefficients $\mathbf{a}%
	=(a_{1},a_{2},\cdots ,a_{K})^{\top }$ and $\mathbf{b}=(b_{1},...,b_{m})^{%
		\top }$\ in a different way. For the first step, we take the GMLS regression
	in (\ref{eqn:LS_optm}) to approximate the function $u$ and then the residual
	of GMLS becomes%
	\begin{equation*}
		\mathbf{s}_{{\mathbf{x}_{0}}}:=\mathbf{u}_{{\mathbf{x}_{0}}}-(\mathcal{I}_{p}%
		\mathbf{u}_{{\mathbf{x}_{0}}})|_{S_{\mathbf{x}_{0}}}=\mathbf{u}_{{\mathbf{x}%
				_{0}}}-\boldsymbol{P}\mathbf{b},
	\end{equation*}%
	where the operator $(\mathcal{I}_{p}\mathbf{u}_{{\mathbf{x}_{0}}})(\mathbf{x}%
	)=\sum_{j=1}^{m}b_{j}p_{\boldsymbol{\alpha }(j)}\left( \boldsymbol{\theta }%
	\left( \mathbf{x}\right) \right) $ is defined in (\ref{eqn:gmls}) and the
	coefficient $\mathbf{b}$\ is given in (\ref{eqn:Pij}). For the second step,
	we fit the above residual $\mathbf{s}_{{\mathbf{x}_{0}}}$ by taking the PHS
	interpolant,
	\begin{equation}
		\left( \mathcal{I}_{\phi }\mathbf{s}_{\mathbf{x}_{0}}\right) \left( \mathbf{x%
		}\right) =\sum_{k=1}^{K}a_{k}\phi \left( \Vert \boldsymbol{\theta }\left(
		\mathbf{x}\right) -\boldsymbol{\theta }\left( \mathbf{x}_{0,k}\right) \Vert
		\right) ,  \label{eqn:Ipph}
	\end{equation}%
	and then the coefficient $\mathbf{a}=\left( a_{1},a_{2},\cdots ,a_{K}\right)
	^{\top }$ are determined by the interpolation condition,
	\begin{equation}
		\boldsymbol{\Phi }\mathbf{a=s}_{{\mathbf{x}_{0}}}=\mathbf{u}_{{\mathbf{x}_{0}%
		}}-\boldsymbol{P}\mathbf{b=u}_{{\mathbf{x}_{0}}}-\boldsymbol{P}\left(
		\boldsymbol{P}^{\top }\boldsymbol{\Lambda P}\right) ^{-1}\boldsymbol{P}%
		^{\top }\boldsymbol{\Lambda }\mathbf{u}_{\mathbf{x}_{0}},  \label{eqn:PHisa}
	\end{equation}%
	where $\boldsymbol{\Phi }$ is given in (\ref{eqn:mat_PRFD}). Numerically, we
	use the \textit{weighted ridge regression} for the inverse of the PHS matrix
	$\boldsymbol{\Phi }$ since $\boldsymbol{\Phi }$ might be singular (see Sec.
	4.3 of \cite{li2025two}). Then the coefficient $\mathbf{a}$ can be computed
	by
	\begin{equation}
		\mathbf{a}=\boldsymbol{\Phi }^{\dag }\left( \mathbf{u}_{{\mathbf{x}_{0}}}-%
		\boldsymbol{P}\mathbf{b}\right) =\boldsymbol{\Phi }^{\dag }\left( \mathbf{I}-%
		\boldsymbol{P}\left( \boldsymbol{P}^{\top }\boldsymbol{\Lambda P}\right)
		^{-1}\boldsymbol{P}^{\top }\boldsymbol{\Lambda }\right) \mathbf{u}_{\mathbf{x%
			}_{0}},  \label{eqn:a}
	\end{equation}%
	where $\boldsymbol{\Phi }^{\dag }$ is the pseudoinverse of $\boldsymbol{\Phi
	}$ given by
	\begin{equation}
		\boldsymbol{\Phi }^{\dag }=\left( \boldsymbol{\Phi }^{\top }\boldsymbol{%
			\Lambda }\boldsymbol{\Phi }+\delta ^{2}\mathbf{I}\right) ^{-1}\boldsymbol{%
			\Phi }^{\top }\boldsymbol{\Lambda },  \label{eqn:PhiIn}
	\end{equation}%
	with $\delta =10^{-5}$ being the regularization parameter for a stable and
	unique inverse of $\boldsymbol{\Phi }$. We also apply the
	normalization technique for stable computation (see Sec. 4.2 of \cite%
	{li2025two}).
	
	
	Basically, the two-step RBF-FD is to regress the function $\mathbf{u}_{%
		\mathbf{x}_{0}}$\ using GMLS followed by fitting the residual $\mathbf{s}_{{%
			\mathbf{x}_{0}}}$ using the PHS interpolant. GMLS captures the smooth,
	leading component of the Taylor expansion of $u$, while PHS compensates for
	the expansion's remainder contained within the residual. As a result, the two-step RBF-FD method takes exactly the same
	PHS+Poly interpolant as in (\ref{eqn:If2}),
	\begin{equation}
		(\mathcal{I}_{\phi p}\mathbf{u}_{{\mathbf{x}_{0}}})\left( \mathbf{x}\right)
		=(\mathcal{I}_{\phi }\mathbf{s}_{{\mathbf{x}_{0}}})\left( \mathbf{x}\right)
		+\left( \mathcal{I}_{p}\mathbf{u}_{\mathbf{x}_{0}}\right) \left( \mathbf{x}%
		\right) :=\sum_{k=1}^{K}a_{k}\phi \left( \Vert \boldsymbol{\theta }\left(
		\mathbf{x}\right) -\boldsymbol{\theta }\left( \mathbf{x}_{0,k}\right) \Vert
		\right) +\sum_{j=1}^{m}b_{j}p_{\boldsymbol{\alpha }(j)}\left( \boldsymbol{%
			\theta }\left( \mathbf{x}\right) \right) ,\text{ \ }\mathbf{x}\in M,
		\label{eqn:GPRFD_inerpolant}
	\end{equation}%
	where the coefficients $\mathbf{a}$ and $\mathbf{b}$ are determined by (\ref%
	{eqn:a}) and (\ref{eqn:Pij}), respectively.
	
	
	
	\subsection{RBF-FD approach for approximating the Laplace--Beltrami operator}
	
	\label{sec:rbfaway}
	
	We first review the results for approximating the Laplace--Beltrami using the
	RBF-FD method \cite%
	{shankar2015radial,flyer2016role,lehto2017radial,jones2023generalized,li2025two}%
	. The Laplace--Beltrami operator at the interior base point $\mathbf{x}_{0}$\
	can be approximated using the RBF-FD,
	\begin{equation}
		\Delta _{M}u\left( \mathbf{x}_{0}\right) \approx \Delta _{\boldsymbol{\theta
		}}(\mathcal{I}_{\phi p}\mathbf{u}_{{\mathbf{x}_{0}}})\left( \mathbf{x}%
		_{0}\right) =\left. \left( \sum_{k=1}^{K}a_{k}\Delta _{\boldsymbol{\theta }%
		}\phi \left( \Vert \boldsymbol{\theta }\left( \mathbf{x}\right) -\boldsymbol{%
			\theta }\left( \mathbf{x}_{0,k}\right) \Vert \right)
		+\sum_{s=1}^{m}b_{s}\Delta _{\boldsymbol{\theta }}p_{\boldsymbol{\alpha }%
			(s)}\left( \boldsymbol{\theta }\left( \mathbf{x}\right) \right) \right)
		\right\vert _{\mathbf{x}=\mathbf{x}_{0}},  \label{eqn:LB_intr}
	\end{equation}%
	where $\Delta _{\boldsymbol{\theta }}=\partial _{\theta _{1}\theta
		_{1}}^{2}+\partial _{\theta _{2}\theta _{2}}^{2}$. Here, the first term in (%
	\ref{eqn:LB_intr}) can be calculated as
	\begin{equation}
		\Delta _{\boldsymbol{\theta }}\phi \left( \Vert \boldsymbol{\theta }\left(
		\mathbf{x}\right) -\boldsymbol{\theta }\left( \mathbf{x}_{0,k}\right) \Vert
		\right) \big|_{\mathbf{x}=\mathbf{x}_{0}}=\left( 4\kappa ^{2}+4\kappa
		+1\right) \Vert \boldsymbol{\theta }\left( \mathbf{x}_{0}\right) -%
		\boldsymbol{\theta }\left( \mathbf{x}_{0,k}\right) \Vert ^{2\kappa -1},
		\label{eqn:Lphi}
	\end{equation}%
	where the PHS $\phi \left( \Vert \boldsymbol{\theta }\left( \mathbf{x}%
	\right) -\boldsymbol{\theta }\left( \mathbf{x}_{0,k}\right) \Vert \right)
	=\Vert \boldsymbol{\theta }\left( \mathbf{x}\right) -\boldsymbol{\theta }%
	\left( \mathbf{x}_{0,k}\right) \Vert ^{2\kappa +1}$ has been used. The
	second term in (\ref{eqn:LB_intr}) can be calculated as%
	\begin{equation}
		\Delta _{\boldsymbol{\theta }}p_{\boldsymbol{\alpha }(s)}\left( \boldsymbol{%
			\theta }\left( \mathbf{x}\right) \right) \big|_{\mathbf{x}=\mathbf{x}%
			_{0}}=\left\{
		\begin{array}{ll}
			2, & \boldsymbol{\alpha }(s)\in \mathbf{E,} \\
			0, & \boldsymbol{\alpha }(s)\notin \mathbf{E,}%
		\end{array}%
		\right.  \label{eqn:2}
	\end{equation}%
	where the index set $\mathbf{E}=\left\{ \left( \alpha _{1},\alpha
	_{2}\right) \in \mathbb{N}^{2}|\left( \alpha _{1},\alpha _{2}\right) =\left(
	2,0\right) \text{ or }\left( \alpha _{1},\alpha _{2}\right) =\left(
	0,2\right) \right\} $. Substituting (\ref{eqn:Lphi}) and (\ref{eqn:2}) into (%
	\ref{eqn:LB_intr}), we arrive at
	\begin{equation}
		\Delta _{M}u\left( \mathbf{x}_{0}\right) \approx \sum_{k=1}^{K}\underbrace{%
			a_{k}\left( 4\kappa ^{2}+4\kappa +1\right) \Vert \boldsymbol{\theta }\left(
			\mathbf{x}_{0}\right) -\boldsymbol{\theta }\left( \mathbf{x}_{0,k}\right)
			\Vert ^{2\kappa -1}}_{(\Delta _{\boldsymbol{\theta }}\boldsymbol{\phi }_{%
				\mathbf{x}_{0}})\mathbf{a}}+\sum_{\left\{ s|\boldsymbol{\alpha }(s)\in
			\mathbf{E}\right\} }\underbrace{2b_{s}}_{(\Delta _{\boldsymbol{\theta }}%
			\boldsymbol{p}_{\mathbf{x}_{0}})\mathbf{b}}  \label{eqn:Lp2}
	\end{equation}%
	where $\Delta _{\boldsymbol{\theta }}\boldsymbol{\phi }_{\mathbf{x}_{0}}\in
	\mathbb{R}^{1\times K}$ with the $k$th entry given by (\ref{eqn:Lphi}) and $%
	\Delta _{\boldsymbol{\theta }}\boldsymbol{p}_{\mathbf{x}_{0}}\in \mathbb{R}%
	^{1\times m}$ with only two nonzero values of 2. Substituting $\mathbf{a}$
	in (\ref{eqn:a}) and $\mathbf{b}$\ in (\ref{eqn:Pij}) into above equation,
	we obtain the Laplacian coefficients $\left\{ w_{k}\right\} _{k=1}^{K}$ such
	that
	\begin{equation}
		\Delta _{M}u\left( \mathbf{x}_{0}\right) \approx \sum_{k=1}^{K}w_{k}u\left(
		\mathbf{x}_{0,k}\right) =\left( \left( \Delta _{\boldsymbol{\theta }}%
		\boldsymbol{\phi }_{\mathbf{x}_{0}}\right) \boldsymbol{\Phi }^{\dag }\big(%
		\mathbf{I}-\boldsymbol{P}\left( \boldsymbol{P}^{\top }\boldsymbol{\Lambda P}%
		\right) ^{-1}\boldsymbol{P}^{\top }\boldsymbol{\Lambda }\big)+\left( \Delta
		_{\boldsymbol{\theta }}\boldsymbol{p}_{\mathbf{x}_{0}}\right) \left(
		\boldsymbol{P}^{\top }\boldsymbol{\Lambda P}\right) ^{-1}\boldsymbol{P}%
		^{\top }\boldsymbol{\Lambda }\right) \mathbf{u}_{\mathbf{x}_{0}}.
		\label{eqn:LB_weight}
	\end{equation}
	
	\begin{example}
		\label{ex1} In Sec. \ref{sec:pre} and the following Sec. \ref{sec:bound}, we
		consider solving the Poisson equation $\Delta _{M}u=f$ with a Robin boundary
		condition $\frac{\partial u}{\partial \boldsymbol{n}}+u=h$ on a 2D
		semi-torus $M$ embedded in $\mathbb{R}^{3}$ via the following map:
		\begin{equation}
			\mathbf{x}:=\left( x^{1},x^{2},x^{3}\right) =((2+\cos ( \psi
			^{1} ))\cos (\psi ^{2}),(2+\cos (\psi ^{1}))\sin (\psi ^{2}),\sin
			(\psi ^{1}))\in M\subset \mathbb{R}^{3},  \label{eqn:torpar}
		\end{equation}%
		for $0\leq \psi ^{1}<2\pi$ and $0\leq \psi ^{2}\leq \pi$, where $(\psi^1, \psi^2)$ are the two intrinsic coordinates. The two boundaries are $\psi^{2}=0$ and $\psi^{2}=\pi$. The Riemannian metric is given by
		\begin{equation}
			g( \psi^{1},\psi^{2} ) =\left(
			\begin{array}{cc}
				1 & 0 \\
				0 & (2+\cos (\psi ^{1}))^{2}%
			\end{array}%
			\right) .  \label{Eqn:tormet}
		\end{equation}%
		The data points are randomly sampled in the $(\psi^{1}, \psi^{2})$ parameter space following a uniform distribution and then mapped onto the surface $M$ via \eqref{eqn:torpar}. In Sec. \ref{sec:pre}, we illustrate the consistency of the discrete operator obtained from Secs. \ref{sec:rbfaway} and \ref{sec:rbfnear} for approximating the Laplace--Beltrami operator $\Delta _{M}$. In the next section (Sec. \ref{sec:bound}), we will illustrate the consistency of derivative boundary conditions and the convergence of the numerical solutions. We set the analytic solution to be $u=\sin (\psi ^{1})\cos (\psi ^{2} + \frac{\pi}{4})$ and then calculate the source term as $f := \Delta _{M}u$, where the Laplace--Beltrami operator is defined as
		\begin{equation}
			\Delta _{M}u = \frac{1}{\sqrt{|g|}}\frac{\partial }{\partial \psi ^{i}}%
			\left( \sqrt{|g|}g^{ij}\frac{\partial u}{\partial \psi^{j}}\right).
		\end{equation}
		We impose the Robin boundary condition and calculate $h:=\frac{\partial u}{\partial \boldsymbol{n}}+u$ at the two boundaries $\psi^{2}=0$ and $\psi^{2}=\pi$. Finally, we solve for the approximate solution of the PDE subjected to the manufactured $f$ and $h$.
	\end{example}
	
	\begin{figure*}[tbp]
		\centering
		\begin{subfigure}{0.48\textwidth}
			\centering
			\caption{RBF-FD without QP}
			\label{fig:bad_points_no_qp}
			\includegraphics[width=\linewidth]{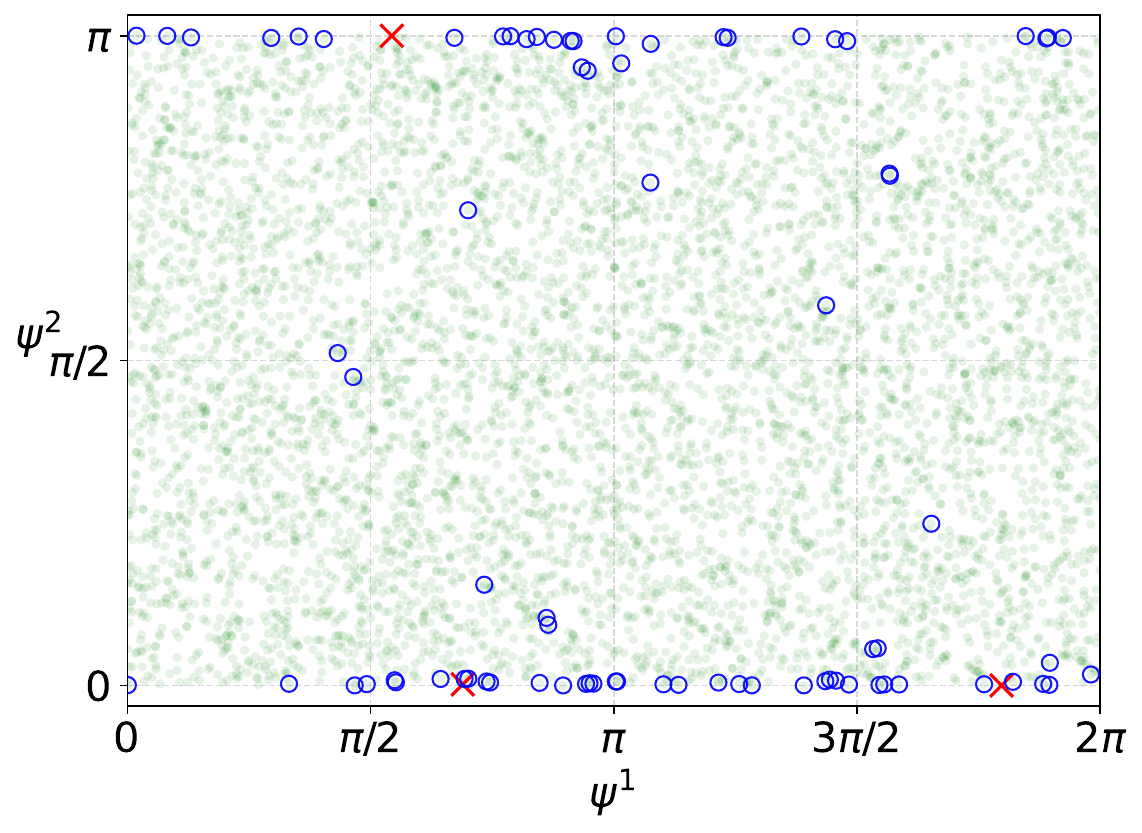}
		\end{subfigure}\hfill
		\begin{subfigure}{0.48\textwidth}
			\centering
			\caption{RBF-FD with QP}
			\label{fig:bad_points_qp}
			\includegraphics[width=\linewidth]{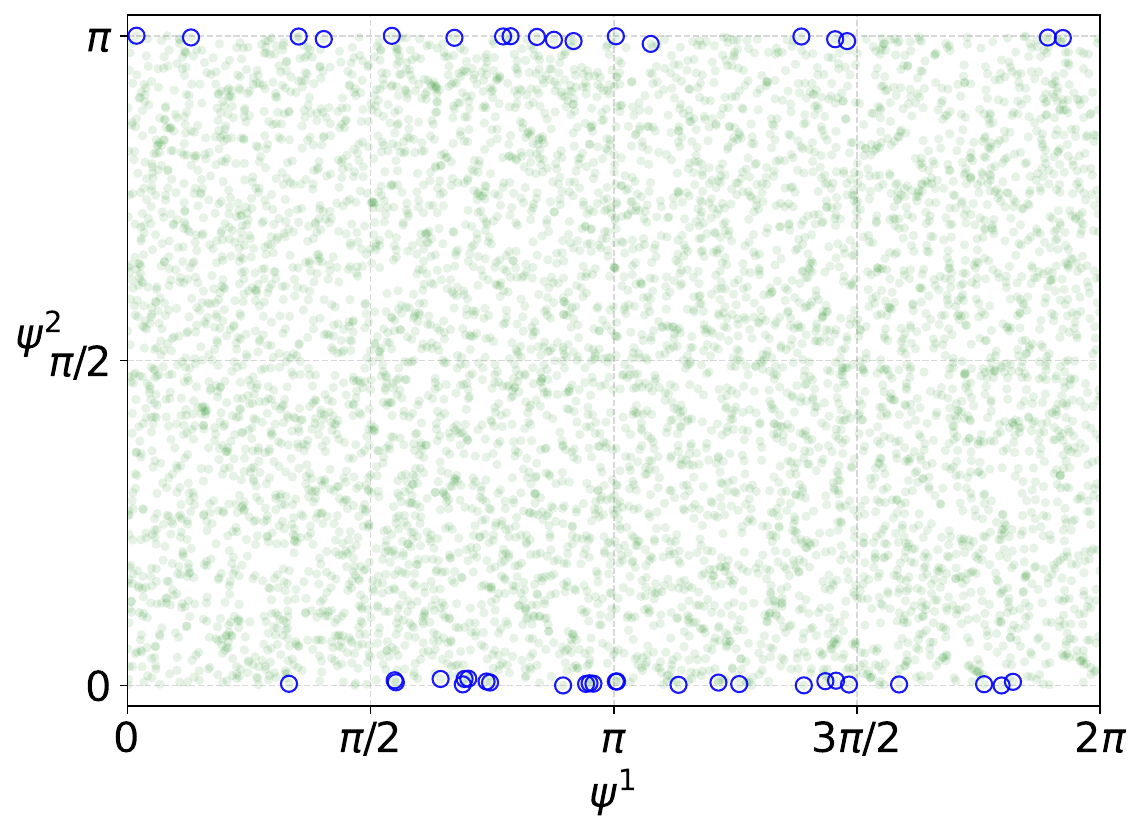}
		\end{subfigure}
		
		\vspace{2ex}
		\includegraphics[width=0.6\textwidth]{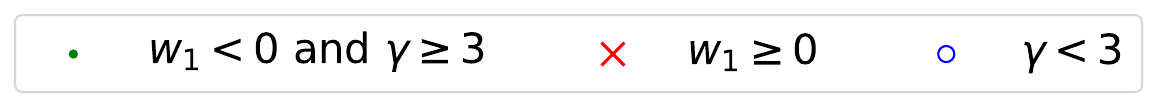}
		\caption{\textbf{Poisson problems on 2D semi-torus in $\mathbb{R}^{3}$}.
			Panels (a) and (b) show three types of interior points in intrinsic
			coordinates $(\psi^1,\psi^2)$ without and with quadratic
			programming, respectively: green dots satisfy both nearly diagonal dominance
			conditions, red crosses violate the first condition $w_1<0$, and blue
			circles violate the second condition $\gamma\geq 3$. The quadratic
			optimization approach completely eliminates points with $w_1 > 0$ while only
			partially eliminating those with $\gamma < 3$. Here no
			boundary points are shown and all points shown are interior. The degree is $l=4$, the number of interior
			points is $N_I = 6240$ and the number of boundary points is $N_B=160$.}
		\label{fig2:torus}
	\end{figure*}
	
	\begin{figure*}[tbp]
		\centering
		\begin{minipage}[c]{0.32\textwidth}
			\centering
			\begin{subfigure}{\linewidth}
				\centering
				\caption{point away from boundary}
				\label{fig:torus_good}
				\includegraphics[width=\linewidth]{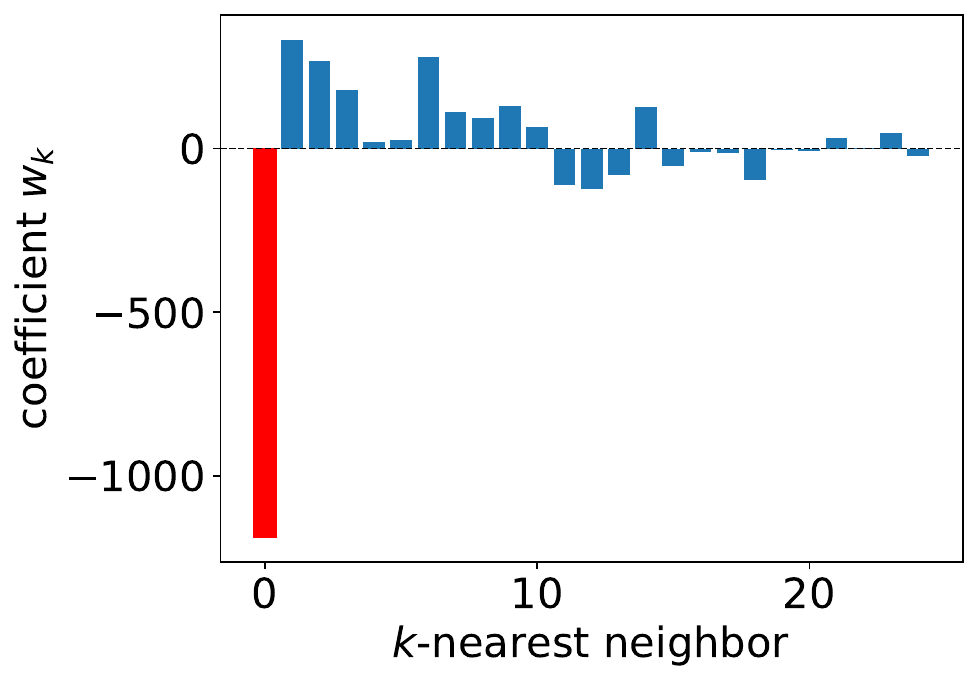}
			\end{subfigure}
		\end{minipage}
		\hfill
		\begin{minipage}[c]{0.66\textwidth}
			\centering
			\begin{subfigure}{0.485\linewidth}
				\centering
				\caption{point near boundary, with $w_1 > 0$}
				\label{fig:torus_pos}
				\includegraphics[width=\linewidth]{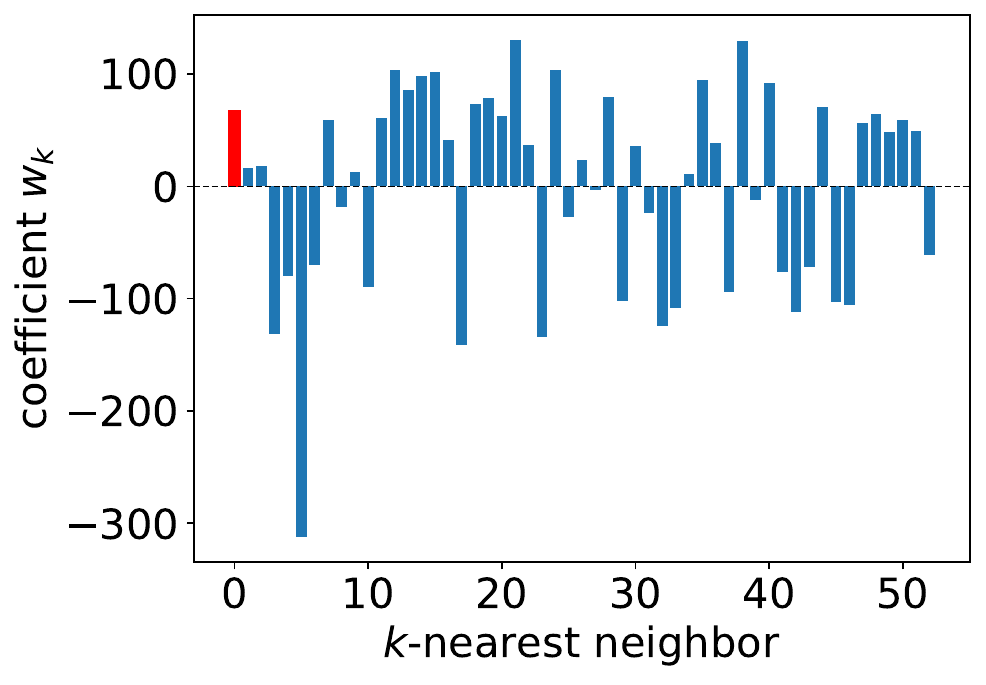}
			\end{subfigure}
			\hfill
			\begin{subfigure}{0.485\linewidth}
				\centering
				\caption{the same point as (\subref{fig:torus_pos}), using QP}
				\label{fig:torus_fixed} 
				\includegraphics[width=\linewidth]{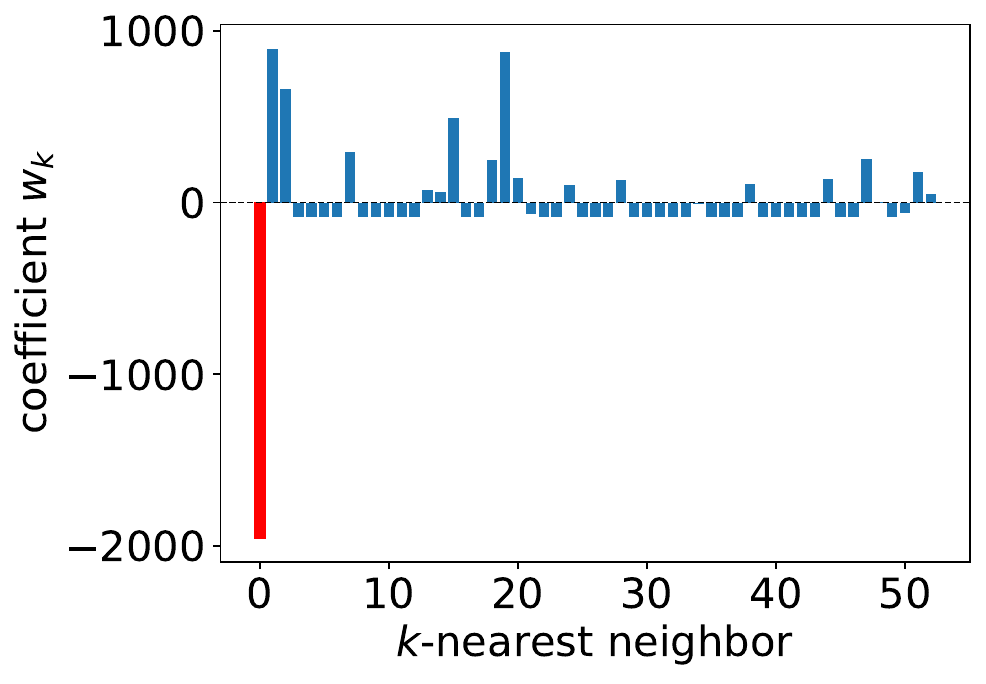}
			\end{subfigure}
			
			\vspace{2ex}
			
			\begin{subfigure}{0.485\linewidth}
				\centering
				\caption{point near boundary, ratio $\gamma<3$}
				\label{fig:torus_ill}
				\includegraphics[width=\linewidth]{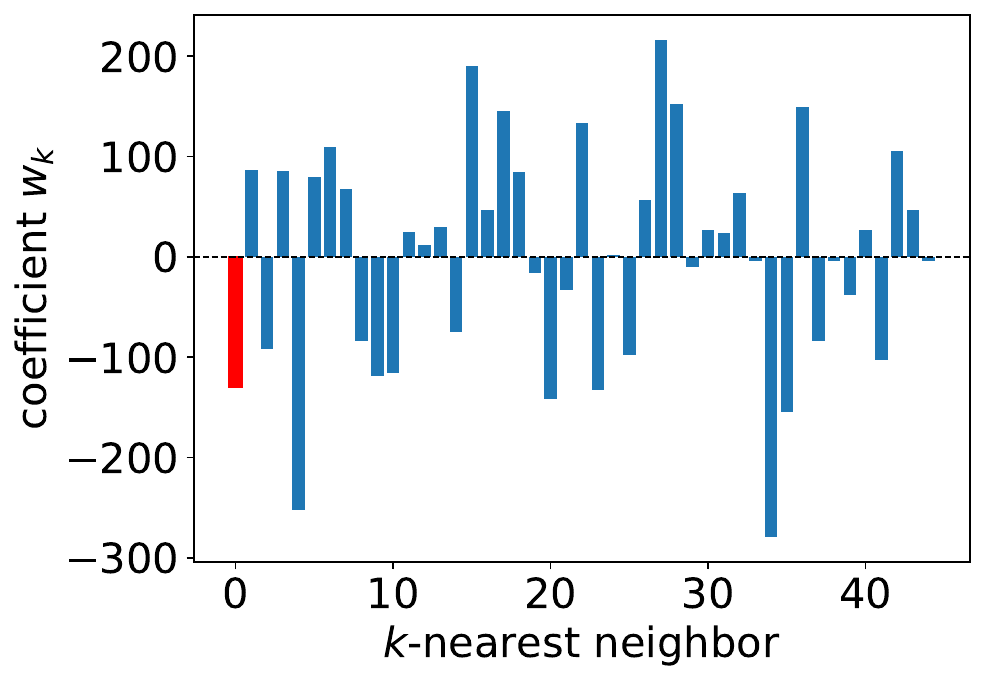}
			\end{subfigure}
			\hfill
			\begin{subfigure}{0.485\linewidth}
				\centering
				\caption{the same point as (\subref{fig:torus_ill}), using QP}
				\label{fig:torus_fixed_ill}
				\includegraphics[width=\linewidth]{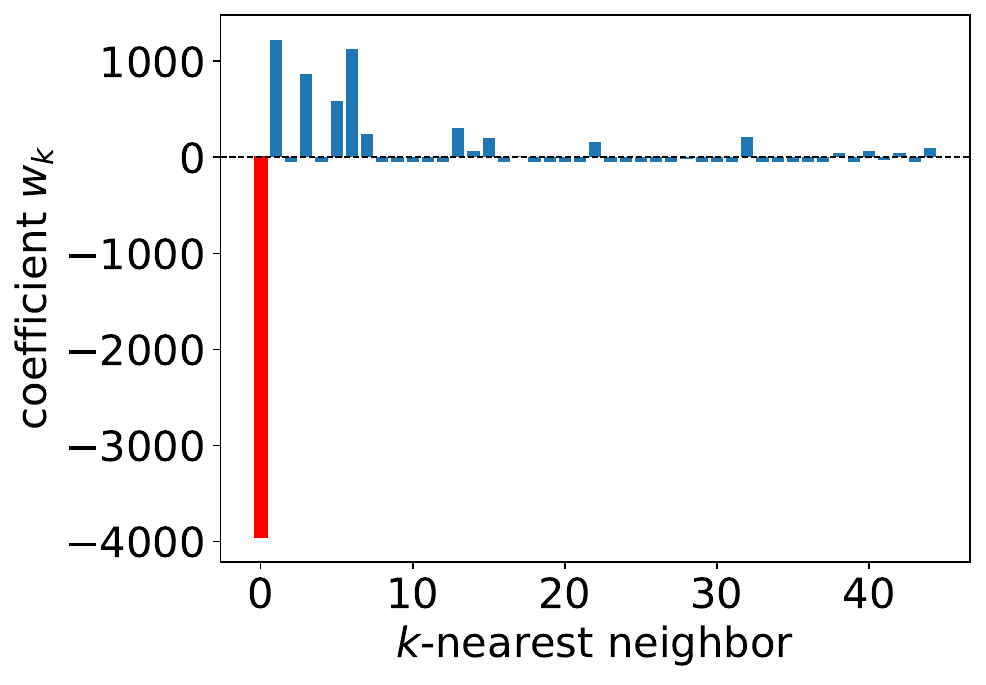}
			\end{subfigure}
		\end{minipage}		
		\caption{\textbf{Poisson problems on 2D semi-torus in $\mathbb{R}^{3}$}.
			Comparison of Laplacian coefficients $w_{1},\ldots,w_k$ for points away from
			and close to the boundary. (a) For most interior points (green dots in Fig.~%
			\ref{fig2:torus}), both nearly diagonal dominance conditions ($w_1<0$
			and $\gamma\geq 3$) can be satisfied. However, for some points close
			to the boundary (red crosses in Fig.~\ref{fig2:torus}(a)), the first
			condition $w_1 < 0$ may not be satisfied (panel (b)) for a certain range of $K$.
			After quadratic programming, $w_1 < 0$ can be satisfied for all these
			points (panel (c)). For other points close to the boundary (blue circles in
			Figs.~\ref{fig2:torus}(a)(b)), the condition $\gamma \geq 3$
			may not be met (panel (d)). After quadratic programming, some points can be
			improved (panel (e)) whereas some points remain to be $\gamma < 3$.
			The degree is $l=4$, the number of interior points is $N_I = 6240$ and the
			number of boundary points is $N_B=160$.}
		\label{fig1:torus}
	\end{figure*}
	
	For the choice of $K$ at each interior base point, we apply the auto-tuned $%
	K $-nearest neighbors \cite{li2025two} such that the following two
	conditions are satisfied:
	\begin{enumerate}
		\item (Sign constraint). The Laplacian coefficient $w_{1}$ at the base point $\mathbf{x}%
		_{0}$ satisfies $w_{1}<0$.
		\item (Ratio constraint). The ratio $\gamma :=\displaystyle{\frac{\left\vert w_{1}\right\vert }{%
				\max_{2\leq k\leq K}\left\vert w_{k}\right\vert }}$ satisfies $\gamma \geq 3$
		for each base point $\mathbf{x}_{0}$.
	\end{enumerate}
	We refer to these two conditions as \textit{nearly diagonal dominance
		conditions}. The details of the auto-tuned $K$-nearest neighbors will be
	shown in the following Algorithm \ref{algo:intrin-LB_ad}. We expect that the
	Laplacian matrix generated by the RBF-FD method is nearly diagonally
	dominant so that the algorithm provides a stable and convergent estimator.
	In our numerical implementation, we see from Fig.~\ref{fig2:torus}(a) that
	most points satisfy both nearly diagonal dominance conditions (green dots in
	Fig.~\ref{fig2:torus}(a)). However, we observe that some points, especially those close to the
	boundary, may fail to satisfy one of the two conditions for all values of $K$
	in a certain range (e.g. $20\sim 60$) (red crosses and blue circles in Fig.~\ref%
	{fig2:torus}(a)). Panels (a)(b)(d) in Fig. \ref{fig1:torus} show the Laplacian
	coefficients given by RBF-FD for the three types of points corresponding to
	green dots, red crosses, and blue circles in Fig.~\ref{fig2:torus}(a),
	respectively. We can see from Fig.~\ref{fig1:torus}(b) and Fig.~\ref%
	{fig1:torus}(d) that the Laplacian coefficients fail to satisfy the first
	condition ($w_{1}<0$) and the second condition ($\gamma \geq 3$),
	respectively. If these RBF-FD Laplacian coefficients are directly employed
		to solve the Poisson problem in Example~\ref{ex1}, the resulting estimator
		will not be stable and the solution will not converge (as will be shown later in Sec. \ref{sec:study_QPstab}). To 	make the algorithm stable and convergent, we consider the following
	quadratic programming approach for points that violate either of the two
	nearly diagonal dominance conditions.
	
	
	
	\subsection{Quadratic programming (QP) approach for approximating the
		Laplace--Beltrami operator}\label{sec:QP}
	
	\label{sec:rbfnear} In this section, we first review the equivalence between
	the GMLS approach and the quadratic programming with equality constraints
	\cite{levin1998approximation,liang2013solving}. We then introduce a new
	quadratic programming formulation that incorporates both equality and
	inequality constraints to enhance the stability. Let us now consider the
	following quadratic programming problem: find the Laplacian coefficients $%
	w_{1},\ldots ,w_{K}$ for approximating $\Delta _{M}u({\mathbf{x}_{0}})$
	using $\sum_{k=1}^{K}w_{k}u({\mathbf{x}_{0,k}})$ by minimizing a quadratic
	form%
	\begin{equation*}
		\underset{w_{1},\ldots ,w_{K}}{\min }\sum_{k=1}^{K}\frac{1}{2}\frac{w_{k}^{2}%
		}{\lambda _{k}},
	\end{equation*}%
	subject to the equality constraint,%
	\begin{equation*}
		\sum_{k=1}^{K}w_{k}p_{\alpha (j)}(\boldsymbol{\theta }\left( \mathbf{x}%
		_{0,k}\right) )=\Delta _{M}p_{\alpha (j)}(\boldsymbol{\theta }\left( \mathbf{%
			x}_{0}\right) ),\text{ \ \ }j=1,\ldots ,m.
	\end{equation*}%
	The above quadratic optimization problem can be written in the compact
	matrix form,
	\begin{equation}
		\begin{array}{ll}
			\underset{\boldsymbol{w}}{\min } & \frac{1}{2}\boldsymbol{w}^{\top }%
			\boldsymbol{\Lambda }^{-1}\boldsymbol{w} \\
			\mathrm{subject\text{ }to} & \boldsymbol{P}^{\top }\boldsymbol{w}=(\Delta _{%
				\boldsymbol{\theta }}\boldsymbol{p}_{\mathbf{x}_{0}})^{\top },%
		\end{array}
		\label{eqn:prim}
	\end{equation}%
	where $\Delta _{\boldsymbol{\theta }}\boldsymbol{p}_{\mathbf{x}_{0}}=\left(
	\Delta _{M}p_{\alpha (1)}(\boldsymbol{\theta }\left( \mathbf{x}_{0}\right)
	),\ldots ,\Delta _{M}p_{\alpha (m)}(\boldsymbol{\theta }\left( \mathbf{x}%
	_{0}\right) )\right) \in \mathbb{R}^{1\times m}$ is defined in (\ref{eqn:Lp2}%
	) and $\boldsymbol{w}$ $=(w_{1},\ldots ,w_{K})^{\top }\in \mathbb{R}%
	^{K\times 1}$ is the Laplacian coefficient. The Lagrange dual function
	associated with the primal problem (\ref{eqn:prim}) is defined as,%
	\begin{equation}
		\underset{\boldsymbol{w}}{\min } \quad \frac{1}{2}\boldsymbol{w}^{\top }%
		\boldsymbol{\Lambda }^{-1}\boldsymbol{w+z}^{\top }\left( \boldsymbol{P}%
		^{\top }\boldsymbol{w}-(\Delta _{\boldsymbol{\theta }}\boldsymbol{p}_{%
			\mathbf{x}_{0}})^{\top }\right)  \label{eqn:w}
	\end{equation}%
	where $\boldsymbol{z}\in \mathbb{R}^{m\times 1}$ is the Lagrange multiplier.
	Taking the derivative with $\boldsymbol{w}$, we can find the minimizer is $%
	\boldsymbol{w}=-\boldsymbol{\Lambda }\boldsymbol{Pz}$. By substituting $%
	\boldsymbol{w}$\ into (\ref{eqn:w}), then the Lagrange dual problem is
	defined as%
	\begin{equation*}
		\max_{\boldsymbol{z}} \quad -\frac{1}{2}\boldsymbol{z}^{\top }\boldsymbol{P}%
		^{\top }\mathbf{\Lambda }\boldsymbol{Pz}-(\Delta _{\boldsymbol{\theta }}%
		\boldsymbol{p}_{\mathbf{x}_{0}})\boldsymbol{z},
	\end{equation*}%
	and the solution to the dual problem is $\boldsymbol{z}=-\left( \boldsymbol{P%
	}^{\top }\boldsymbol{\Lambda }\boldsymbol{P}\right) ^{-1}(\Delta _{%
		\boldsymbol{\theta }}\boldsymbol{p}_{\mathbf{x}_{0}})^{\top }$. Then the
	Laplacian coefficient $\boldsymbol{w}$\ can be calculated by%
	\begin{equation*}
		\boldsymbol{w^{\top }}=(\Delta _{\boldsymbol{\theta }}\boldsymbol{p}_{%
			\mathbf{x}_{0}})\left( \boldsymbol{P}^{\top }\boldsymbol{\Lambda }%
		\boldsymbol{P}\right) ^{-1}\boldsymbol{P}^{\top }\boldsymbol{\Lambda }\in
		\mathbb{R}^{1\times K}.
	\end{equation*}%
	This is equivalent to the GMLS formula as can be seen from the second term
	of Eq. (\ref{eqn:LB_weight}).
	
	
	We now impose additional constraints to stabilize the Laplacian
	approximation for the points which cannot satisfy the two nearly diagonal
	dominance conditions. For those points, we consider the following quadratic
	optimization problem with additional inequality constraints:%
	\begin{equation}
		\begin{array}{ll}
			\underset{w_{1},\ldots ,w_{K},C}{\min } & \sum_{k=1}^{K}\displaystyle{\frac{1%
				}{2}\frac{w_{k}^{2}+C^{2}}{\lambda _{k}}} \\
			\mathrm{subject\text{ }to} & \sum_{k=1}^{K}w_{k}p_{\alpha (j)}(\boldsymbol{%
				\theta }\left( \mathbf{x}_{0,k}\right) )=\Delta _{M}p_{\alpha (j)}(%
			\boldsymbol{\theta }\left( \mathbf{x}_{0}\right) ),\text{ \ \ }j=1,\ldots ,m,
			\\
			& w_{1}\leq -C, \\
			& w_{k}\geq -C,\text{ \ \ }k=2,\ldots ,K, \\
			& C\geq 0.%
		\end{array}
		\label{eqn:inequad}
	\end{equation}%
	For these inequality constraints, we relax the diagonal dominance
	restriction to allow that the Laplacian matrix to be close enough to a
	diagonally dominant matrix. This is expected to enhance the stability of
	approximation and improve the convergence of numerical solution. As a result, we
	construct the sparse Laplacian matrix$\ \boldsymbol{L}_{I,M}\in \mathbb{R}%
	^{N_{I}\times N}$ using an RBF-FD approach coupled with quadratic
	programming: for each interior node, we first attempt an auto-tuned $K$ search to enforce nearly diagonal dominance. If this attempt fails, we then solve the QP formulation in (\ref{eqn:inequad}) based on an appropriate stencil size $K$ yielding nearly diagonal dominance or, otherwise, the largest ratio $\gamma$. We refer to this approach  as \textbf{RBF-FD-QP} method and summarize the algorithm in Algorithm~%
	\ref{algo:intrin-LB_ad}.
	
	\begin{algorithm}[ht]
		\caption{RBF-FD-QP using auto-tuned $K$ for the Laplace--Beltrami operator}
		\begin{algorithmic}[1]
			\STATE {\bf Input:} A point cloud $\{\mathbf{x}_{i}\}_{i=1}^{N}\subset M$, bases of analytic tangent vectors at each node $\left\{ \boldsymbol{t}_{1}(\mathbf{x}_{i}),\boldsymbol{t}_{2}(%
			\mathbf{x}_{i})\right\} _{i=1}^{N}$, the degree $l$ of  polynomials, the PHS smoothness parameter $\kappa$, and an initial $K_0>m$-nearest neighbors where $%
			m=(l+2)(l+1)/2 $ is the number of monomial basis functions $\left\{ \boldsymbol{%
				\theta }^{\boldsymbol{\alpha }}|0\leq \left\vert \boldsymbol{\alpha }%
			\right\vert \leq l\right\} $, the maximum number of  nearest neighbors, $K_{\max}$.
			\STATE Set $\boldsymbol{L}_{I,M}$ to be a sparse $N_I\times N$ matrix with initial $N_I K_0$ nonzeros.
			\FOR{$i\in \{1,2,...,N_I\}$}
			\STATE \%\%\% RBF-FD for computing Laplacian weight
			\STATE Set $\gamma = 0$,  $w_1 = 0$, $K=K_0$.
			\WHILE{$\gamma < 3$ or $w_1>0$ }
			\STATE Find the $K$-nearest neighbors of the interior point $\mathbf{x}_i$, denoted by the stencil $S_{\mathbf{x}_i}=\{\mathbf{x}_{i,k}\}_{k=1}^K$.
			\STATE  Construct the $K$ by $m$ Vandermonde-type matrix $\boldsymbol{P}$ with entries
			$\boldsymbol{P}_{kj} = p_{\boldsymbol{\alpha }(j)}\left( \boldsymbol{%
				\theta}\left( \mathbf{x}{_{i,k}}\right) \right) $ defined in (\ref{eq:thet12}).
			\STATE Construct the $K$ by $K$ PHS matrix $\boldsymbol{\Phi}$ with entries
			$\boldsymbol{\Phi}_{ks} = \Vert \boldsymbol{\theta}\left(\mathbf{x}_{i,k}\right) - \boldsymbol{\theta}\left(\mathbf{x}_{i,s}\right)\Vert^{2\kappa+1}$.
			\STATE Calculate the Laplacian coefficients $\{w_k\}_{k=1}^{K}$ by the formula in \eqref{eqn:LB_weight},%
			\begin{equation*}
				(w_1,\ldots,w_K) =\left( \Delta _{\boldsymbol{\theta }}%
				\boldsymbol{\phi }_{\mathbf{x}_{i,0}}\right) \boldsymbol{\Phi }^{\dag }\big(%
				\mathbf{I}-\boldsymbol{P}\left( \boldsymbol{P}^{\top }\boldsymbol{\Lambda P}%
				\right) ^{-1}\boldsymbol{P}^{\top }\boldsymbol{\Lambda }\big)+\left( \Delta
				_{\boldsymbol{\theta }}\boldsymbol{p}_{\mathbf{x}_{i,0}}\right) \left(
				\boldsymbol{P}^{\top }\boldsymbol{\Lambda P}\right) ^{-1}\boldsymbol{P}%
				^{\top }\boldsymbol{\Lambda }.
			\end{equation*}
			\STATE Calculate the ratio $\gamma = \frac{|w_1|}{ \max_{2\leq k\leq K}|w_k |}$. 
			\STATE Increase the $K$ by a small positive integer (e.g., 2).
			\IF{$K \geq K_{\max}$}
			\STATE Break.
			\ENDIF
			\ENDWHILE
			
			\STATE \%\%\% QP for computing Laplacian weight for point that violate the nearly diagonal dominance
			conditions.
			\STATE Set $K=K_0$.
			\WHILE{$\gamma < 3$ or $w_1>0$ }
			\STATE Solve the quadratic optimization problem (\ref{eqn:inequad}).
			\IF{The problem (\ref{eqn:inequad}) has a solution and $\gamma\geq3$  }
			\STATE Break.
			\ELSIF{The problem (\ref{eqn:inequad}) has a solution and $K<K_{\max}$}
			\STATE Record the current values of $(w_1,\ldots,w_K)$, $\gamma$ and $K$.
			\STATE Increase the $K$ by a small positive integer (e.g., 2).
			\ELSE
			\STATE Let $(w_1,\ldots,w_K)$ be the Laplacian weights  for the recorded $K$ that yields the largest $\gamma$.
			\STATE Break.
			\ENDIF
			\ENDWHILE
			
			\STATE Arrange the Laplacian coefficients $\{w_k\}_{k=1}^{K}$ into the corresponding rows and columns of $\boldsymbol{L}_{I,M}$.
			\ENDFOR
			\STATE {\bf Output:} The approximate operator matrix $\boldsymbol{L}_{I,M}$.
		\end{algorithmic}
		\label{algo:intrin-LB_ad}
	\end{algorithm}

	\begin{figure*}[tbp]
		\centering
		\begin{subfigure}{0.5\textwidth}
			\centering
			\caption{pointwise absolute error with $N=6400$, $l=4$}
			\label{fig:fe_interior_pointwise}
			\includegraphics[height=6cm, keepaspectratio]{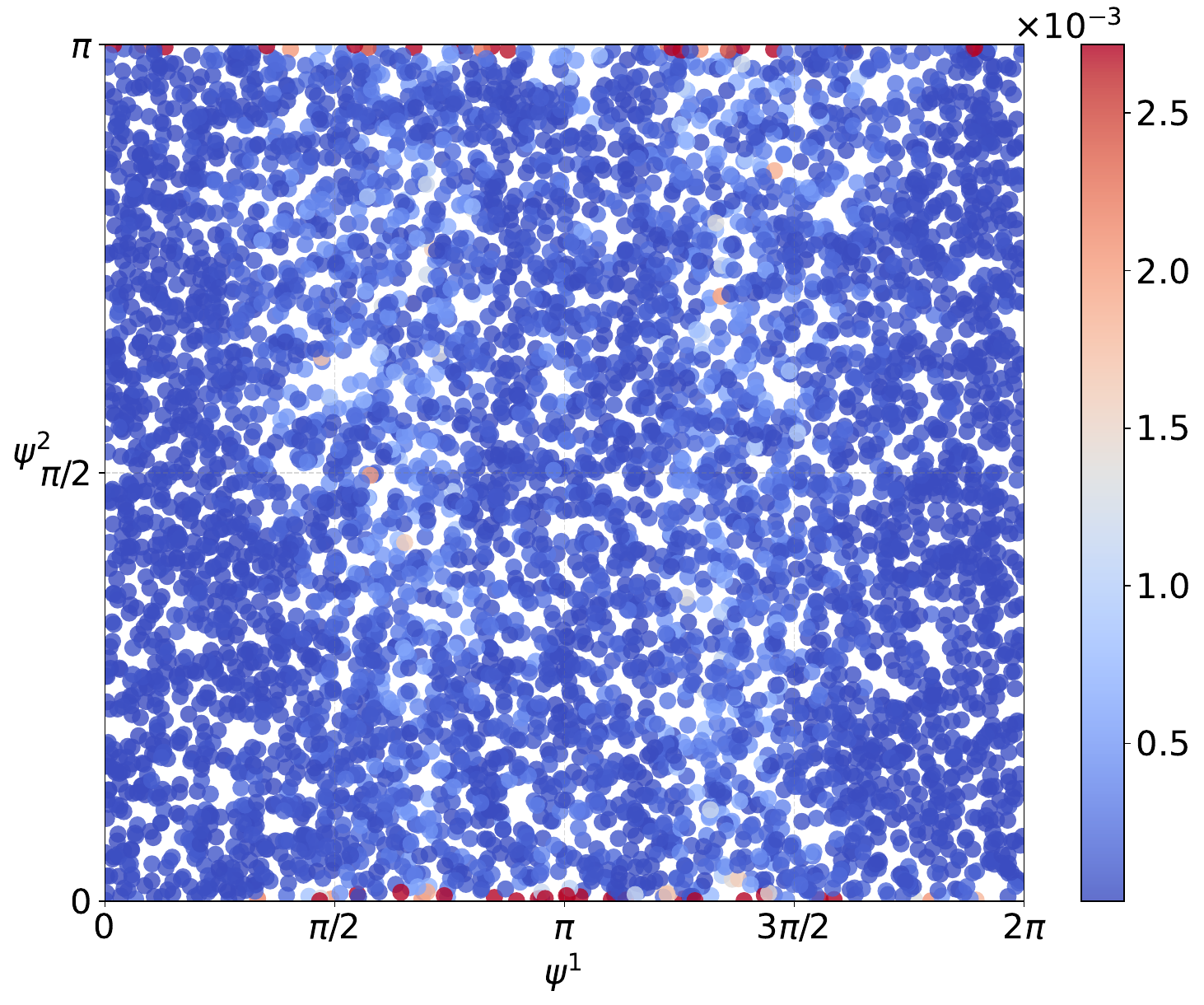}
		\end{subfigure}\hfill
		\begin{subfigure}{0.45\textwidth}
			\centering
			\caption{consistent \textbf{FE}s with different degrees}
			\label{fig:fe_interior_convergence}
			\includegraphics[height=6cm, keepaspectratio]{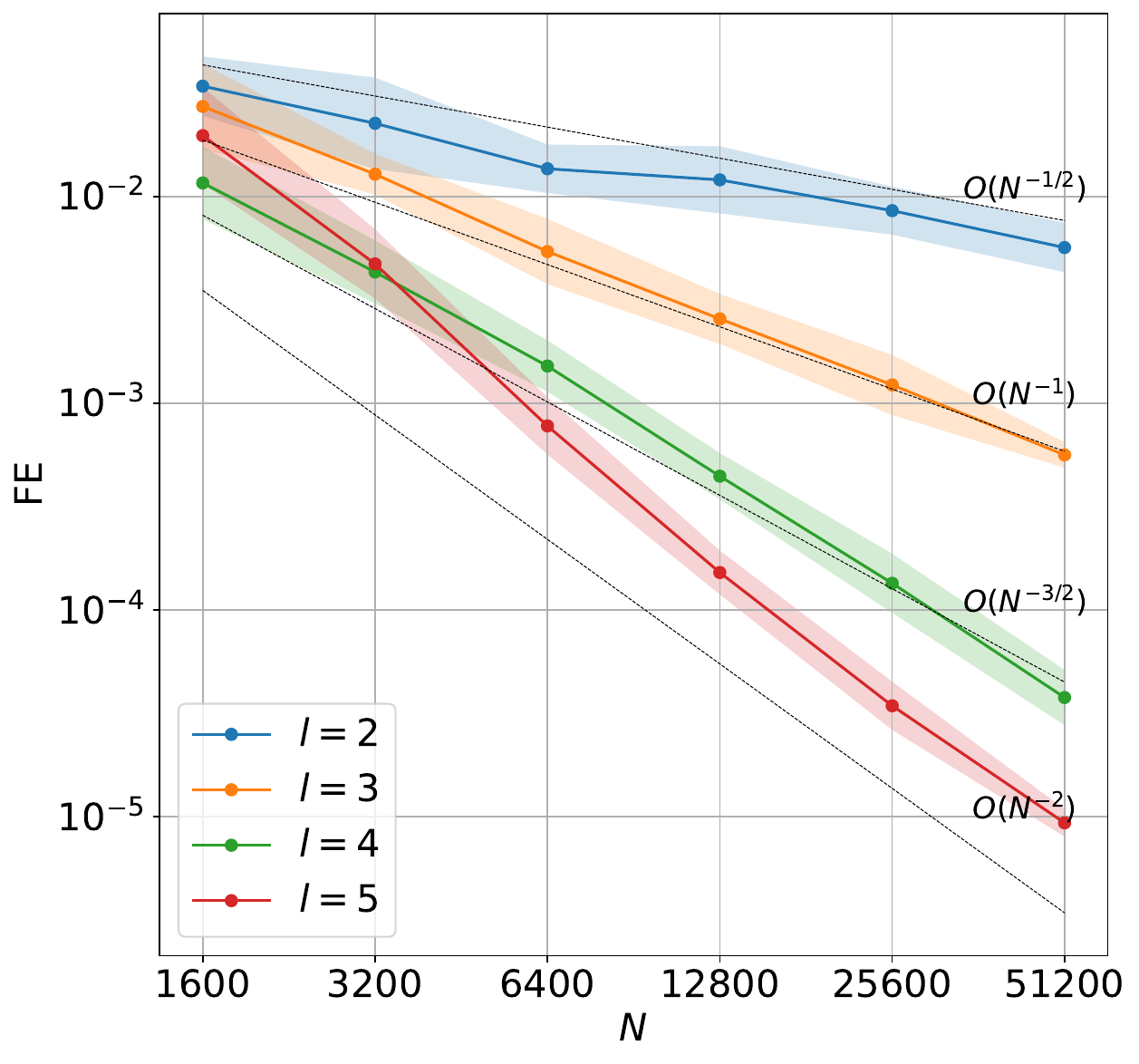}
		\end{subfigure}
		\caption{\textbf{Poisson problems on 2D semi-torus in $\mathbb{R}^{3}$}.
			Panel (a) shows the pointwise absolute error of $\Delta_M$ using our
			combined RBF-FD and QP approach with $N=6400$ and $l=4$. Panel (b) shows the
			consistency of \textbf{FE}s for different polynomial degrees. All simulations are run with 12 independent trials, each with a set of randomly sampled data points. }
		\label{fig3:torus}
	\end{figure*}
	
	After applying the QP approach, Fig. \ref{fig2:torus}(b) shows that all
		points with $w_1 > 0$ are completely eliminated, as indicated by the absence
		of red crosses. Indeed, this complete elimination of points with $w_1 > 0$
		was consistently observed across all simulation trials.
		However, points with $\gamma < 3$ are only partially removed, with a small
		number still remaining (see blue circles in Fig.~\ref{fig2:torus}(b) and see also the resulting Laplacian coefficients obtained from	the QP approach in Fig.~\ref{fig1:torus}(c)). While some of the ratios are still $\gamma < 3$, we numerically observe that
		all of them have $\gamma \geq 1$  (see the coefficients in Fig.~\ref{fig1:torus}(c)). Fortunately,
		despite the presence of some points with $1\leq \gamma < 3$, the estimator becomes
		sufficiently stable to yield convergent solutions, as will be demonstrated
		below in Sec.~\ref{sec:bound}.
	
	
	To numerically verify the validity of operator approximation, we define the
	forward error (\textbf{FE}) for {the} operator approximation as
	\begin{equation}
		\mathbf{FE}=\left( \frac{1}{N_I}\sum\limits_{i=1}^{N_I}\left( \Delta _{M}u\left(
		\mathbf{x}_{i}\right) -\left( \boldsymbol{L}_{I,M}\mathbf{u}\right)
		_{i}\right) ^{2}\right) ^{1/2},  \label{eqn:LapAp}
	\end{equation}%
	where $\mathbf{u}:=(u(\mathbf{x}_{1}),\ldots ,u(\mathbf{x}_{N}))^{\top }\in
	\mathbb{R}^{N\times 1}$ is the analytic solution and $\left( \boldsymbol{L}%
	_{I,M}\mathbf{u}\right) _{i}$ is the $i$th component of $\boldsymbol{L}_{I,M}%
	\mathbf{u}\in \mathbb{R}^{N_{I}\times 1}$. Fig.~\ref{fig3:torus} shows the
	pointwise absolute error of $\Delta _{M}$ (panel (a)) and the consistency of
	\textbf{FE}s for different polynomial degrees (panel (b)). As can be seen from Fig. %
		\ref{fig3:torus}(a), points near the boundary, whose coefficients are obtained from the QP approach, exhibit  relatively large pointwise errors compared to those away from the boundary.
	As shown in Fig. %
	\ref{fig3:torus}(b), the consistency rate is approximately $O(N^{-(l-1)/2})$%
	, which is in good agreement with the theoretical error bound of $O((\log
	N/N)^{(l-1)/2})$ (see e.g., \cite%
	{levin1998approximation,mirzaei2012generalized,liang2013solving}).
	
	\section{Application to solving boundary problems}
	
	\label{sec:bound}
	
	In this section, we extend the RBF-FD-QP approach, introduced in Sec.~\ref{sec:pre} for interior points, to boundary points for solving BVPs on surfaces with derivative boundary conditions. The task here is to approximate the outward co-normal derivative  at boundary points using one-sided stencils. We first derive the weights generated by our two-step RBF-FD approach on the boundary (Sec.~\ref{sec:rbf_normal}), and then describe the numerical implementation of our RBF-FD-QP approach for boundary operator approximation (Sec.~\ref{sec:algo_boundary}). Last, we assemble the discrete estimators of interior and boundary operators into a sparse linear system for solving the Poisson problems in Example~\ref{ex1}.

	
	\subsection{RBF-FD approach for approximating the outward co-normal derivative on boundaries}
	
	\label{sec:rbf_normal}
	
	To solve boundary value problems with derivative boundary conditions on
	surfaces, we need to approximate the outward co-normal derivative for points
	located on the boundary $\partial M$. 	
	Let $\mathbf{x}_{b}\in \partial M$ be
	a boundary base point. The outward co-normal derivative of a smooth function $u$
	at $\mathbf{x}_{b}$ is defined as $\frac{\partial u}{\partial \boldsymbol{n}}%
	(\mathbf{x}_{b})=\nabla _{M}u(\mathbf{x}_{b})\cdot \boldsymbol{n}$, where $%
	\nabla _{M}u(\mathbf{x}_{b})\in \mathbb{R}^{3}$ is the surface gradient and $%
	\boldsymbol{n}\in \mathbb{R}^{3}$ is the unit outward co-normal vector at $%
	\mathbf{x}_{b}$. Note that $\boldsymbol{n}$ lies in the tangent space $T_{%
		\mathbf{x}_{b}}M$ and is orthogonal to the boundary curve $\partial M$.
	
	We first approximate the surface gradient $\nabla _{M}u(\mathbf{x}_{b})\in
	T_{\mathbf{x}_{b}}M$. Using the local Monge coordinate system $\boldsymbol{%
		\theta }(\mathbf{x})=(\theta _{1}(\mathbf{x}),\theta _{2}(\mathbf{x}))$, the
	surface gradient can be expressed by pushing forward the local intrinsic
	gradient $\nabla _{\boldsymbol{\theta }}u(\mathbf{x}_{b})\in \mathbb{R}%
	^{2\times 1}$ via the tangent basis matrix $\mathbf{T}_{\mathbf{x}_{b}}=[%
	\boldsymbol{t}_{1}(\mathbf{x}_{b}),\boldsymbol{t}_{2}(\mathbf{x}_{b})]\in
	\mathbb{R}^{3\times 2}$. Applying the two-step tangent-space RBF-FD
	interpolant defined in (\ref{eqn:GPRFD_inerpolant}), we have
	\begin{equation}
		\begin{aligned} \nabla_M u(\mathbf{x}_b) &\approx \mathbf{T}_{\mathbf{x}_b}
			\nabla_{\boldsymbol{\theta}} (\mathcal{I}_{\phi p}
			\mathbf{u}_{\mathbf{x}_b})(\mathbf{x}) \big|_{\mathbf{x}=\mathbf{x}_b} \\ &=
			\mathbf{T}_{\mathbf{x}_b} \left( \sum_{k=1}^K a_k
			\nabla_{\boldsymbol{\theta}} \phi \left( \| \boldsymbol{\theta}(\mathbf{x})
			- \boldsymbol{\theta}(\mathbf{x}_{b,k}) \| \right) + \sum_{s=1}^m b_s
			\nabla_{\boldsymbol{\theta}}
			p_{\boldsymbol{\alpha}(s)}(\boldsymbol{\theta}(\mathbf{x})) \right)
			\Bigg|_{\mathbf{x}=\mathbf{x}_b}, \end{aligned}  \label{eqn:grad_approx}
	\end{equation}%
	where $\left\{ \mathbf{x}_{b,k}\in \mathbf{X}_{M}\right\} _{k=1}^{K}$ is a specifically selected stencil of
		$K$ neighbors of $\mathbf{x}_{b}=\mathbf{x}_{b,1}$ and $\mathbf{u}_{%
			{\mathbf{x}_{b}}}:=(u({\mathbf{x}_{b,1}}),...,u(\mathbf{x}_{b,K}))^{\top }$
		is the corresponding $K$ function values. The choice of the $K$ neighbors will be clarified below in the following Sec.~\ref{sec:algo_boundary}. Evaluating the gradients at the boundary base
	point $\mathbf{x}_{b}\ $(or correspondingly, $\boldsymbol{\theta }(\mathbf{x}%
	_{b})=\mathbf{0}$\textbf{)}, the derivative of the PHS function $\phi
	(r)=r^{2\kappa +1}$ yields a $2\times 1$ column vector:
	\begin{equation}
		\nabla _{\boldsymbol{\theta }}\phi \left( \Vert \boldsymbol{\theta }(\mathbf{%
			x})-\boldsymbol{\theta }(\mathbf{x}_{b,k})\Vert \right) \big|_{\mathbf{x}=%
			\mathbf{x}_{b}}=-(2\kappa +1)\Vert \boldsymbol{\theta }(\mathbf{x}%
		_{b,k})\Vert ^{2\kappa -1}\boldsymbol{\theta }(\mathbf{x}_{b,k}) \in \mathbb{%
			R}^2.  \label{eqn:grad_phi}
	\end{equation}%
	For the polynomial terms $p_{\boldsymbol{\alpha }(s)}(\boldsymbol{\theta }(%
	\mathbf{x}))$, all quadratic and higher-order terms vanish at the origin $%
	\boldsymbol{\theta }(\mathbf{x}_{b})=\mathbf{0}$, leaving only the gradients
	of the linear terms:
	\begin{equation}
		\nabla _{\boldsymbol{\theta }}p_{\boldsymbol{\alpha }(s)}(\boldsymbol{\theta
		}(\mathbf{x}))\big|_{\mathbf{x}=\mathbf{x}_{b}}=%
		\begin{cases}
			(1,0)^{\top }, & \text{if }\boldsymbol{\alpha }(s)=(1,0), \\
			(0,1)^{\top }, & \text{if }\boldsymbol{\alpha }(s)=(0,1), \\
			(0,0)^{\top }, & \text{otherwise}.%
		\end{cases}
		\label{eqn:grad_poly}
	\end{equation}%
	Substituting (\ref{eqn:grad_phi}) and (\ref{eqn:grad_poly}) into (\ref%
	{eqn:grad_approx}) gives the matrix-vector form:
	\begin{equation}
		\nabla _{M}u(\mathbf{x}_{b})\approx \mathbf{T}_{\mathbf{x}_{b}}\left(
		(\nabla _{\boldsymbol{\theta }}\boldsymbol{\Phi }_{\mathbf{x}_{b}})\mathbf{a}%
		+(\nabla _{\boldsymbol{\theta }}\boldsymbol{P}_{\mathbf{x}_{b}})\mathbf{b}%
		\right) ,  \label{eqn:grad_matrix_form}
	\end{equation}%
	where $\nabla _{\boldsymbol{\theta }}\boldsymbol{\Phi }_{\mathbf{x}_{b}}\in
	\mathbb{R}^{2\times K}$ contains the PHS derivatives in (\ref{eqn:grad_phi})
	and $\nabla _{\boldsymbol{\theta }}\boldsymbol{P}_{\mathbf{x}_{b}}\in
	\mathbb{R}^{2\times m}$ contains the polynomial derivatives in (\ref%
	{eqn:grad_poly}).
	
	To obtain the co-normal derivative, we take the inner product of the ambient
	surface gradient with the outward co-normal vector $\boldsymbol{n}$:
	\begin{equation}
		\frac{\partial u}{\partial \boldsymbol{n}}(\mathbf{x}_{b})=\boldsymbol{n}%
		^{\top }\nabla _{M}u(\mathbf{x}_{b})\approx \boldsymbol{n}^{\top }\mathbf{T}%
		_{\mathbf{x}_{b}}\left( (\nabla _{\boldsymbol{\theta }}\boldsymbol{\Phi }_{%
			\mathbf{x}_{b}})\mathbf{a}+(\nabla _{\boldsymbol{\theta }}\boldsymbol{P}_{%
			\mathbf{x}_{b}})\mathbf{b}\right) .  \label{eqn:normal_deriv_inner}
	\end{equation}%
	Notice that the vector-matrix product $\boldsymbol{n}^{\top }\mathbf{T}_{%
		\mathbf{x}_{b}}$ computes the projection of the co-normal vector onto the local
	tangent space basis. Finally, by substituting the coefficients $\mathbf{a}$
	in (\ref{eqn:a}) and $\mathbf{b}$ in (\ref{eqn:Pij}) into (\ref%
	{eqn:normal_deriv_inner}), we obtain the weights $%
	\{v_{k}\}_{k=1}^{K}$ such that
	\begin{equation}
		\begin{aligned} \frac{\partial u}{\partial \boldsymbol{n}}(\mathbf{x}_b)
			&\approx \sum_{k=1}^K v_k u(\mathbf{x}_{b,k}) \\ &= \boldsymbol{n}^\top
			\mathbf{T}_{\mathbf{x}_b} \left[ (\nabla_{\boldsymbol{\theta}}
			\boldsymbol{\Phi}_{\mathbf{x}_b}) \boldsymbol{\Phi}^\dagger \big(\mathbf{I}
			- \boldsymbol{P}(\boldsymbol{P}^\top \boldsymbol{\Lambda}
			\boldsymbol{P})^{-1} \boldsymbol{P}^\top \boldsymbol{\Lambda}\big) +
			(\nabla_{\boldsymbol{\theta}} \boldsymbol{P}_{\mathbf{x}_b})
			(\boldsymbol{P}^\top \boldsymbol{\Lambda} \boldsymbol{P})^{-1}
			\boldsymbol{P}^\top \boldsymbol{\Lambda} \right] \mathbf{u}_{\mathbf{x}_b}.
		\end{aligned}  \label{eqn:normal_weight}
	\end{equation}%
	This formulation yields the weights required to enforce Neumann or Robin
	boundary conditions by applying the two-step RBF-FD method.
	
	\subsection{Algorithmic implementation for boundary operators}
	
	\label{sec:algo_boundary}
	
	With the mathematical formulation established in Sec. \ref{sec:rbf_normal},
	we now describe the practical implementation for approximating the outward
	co-normal derivative at boundary points. The procedure shares similarities with
	the evaluation of the interior Laplace-Beltrami operator but requires
	specific considerations for stencil selection, weight function $\boldsymbol{%
		\Lambda}$, and the auto-tuning of the neighborhood size $K$. Unlike the interior Laplace-Beltrami operator, the boundary operator should be constructed from a one-sided local stencil that accounts for the directional derivative along the outward co-normal direction.
	

	When dealing with a boundary base point $\mathbf{x}_{b}\in \partial M$, the
	local geometry is heavily one-sided due to the truncation of the domain. To
	collect the stencil $S_{\mathbf{x}_{b}}=\{\mathbf{x}_{b,k}\}_{k=1}^{K}$, we
	consider a restricted $K$-nearest neighbor (KNN) search strategy that includes only
	the boundary base point $\mathbf{x}_{b}$ and its neighboring interior points, as illustrated in Fig. \ref%
	{fig:stencil_and_fe}(a). That is to say, the restricted KNN method excludes all other boundary points except the base $\mathbf{x}_b$ itself. This restricted KNN search strategy ensures that the value of the boundary base point can be written as a linear combination of the values of neighboring interior points after discretization (as will be shown below in (\ref{eqn:uB})).

\begin{figure*}[htbp]
		\centering
		\begin{minipage}[t]{0.34\textwidth}
			\centering
			\begin{subfigure}{\linewidth}
				\centering
				\caption{restricted KNN}
				\label{fig:stencil_a}
				\includegraphics[height=3.8cm]{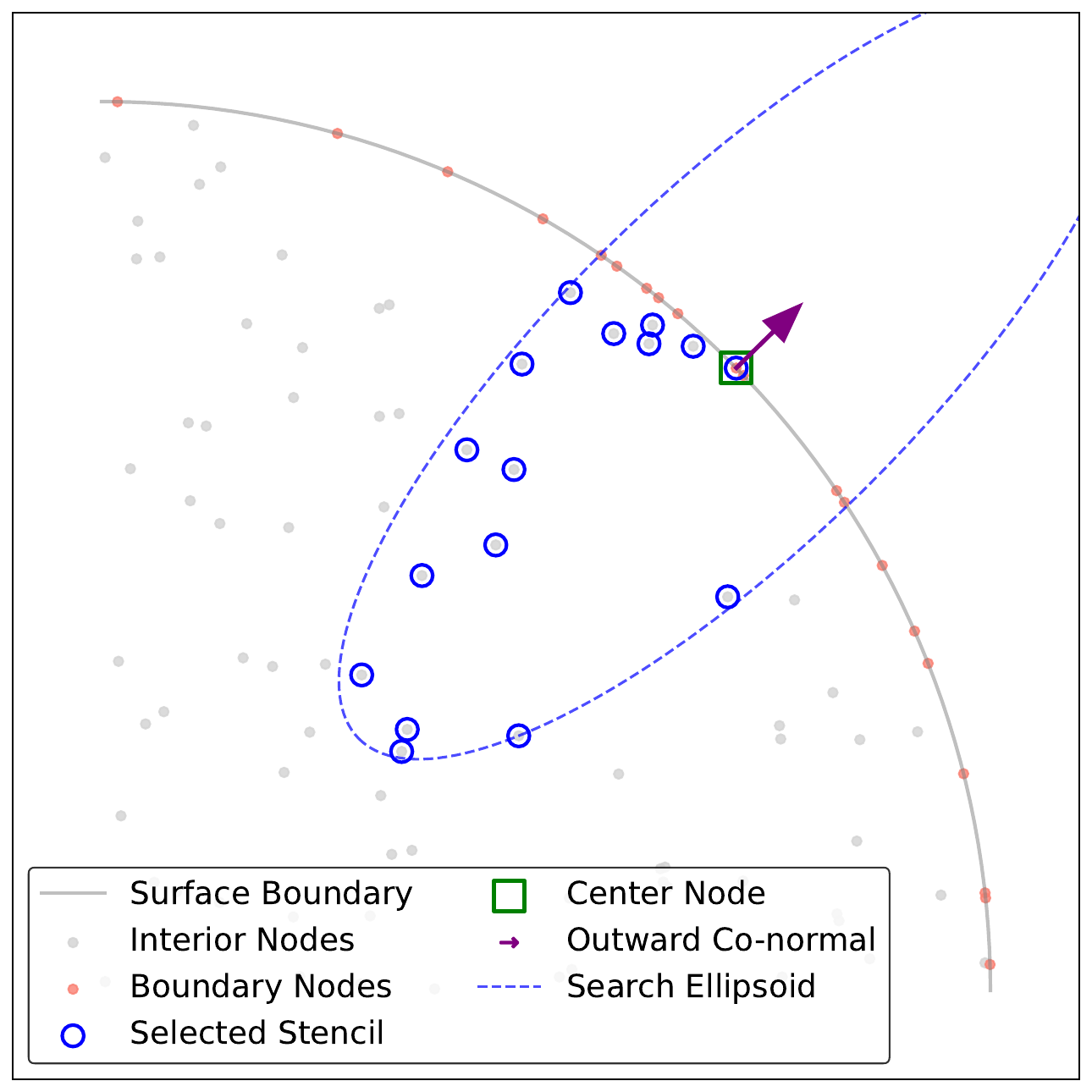}
			\end{subfigure}
			
			\vspace{3ex}
			\begin{subfigure}{\linewidth}
				\centering
				\caption{typical distribution}
				\label{fig:stencil_b}
				\includegraphics[height=3.8cm]{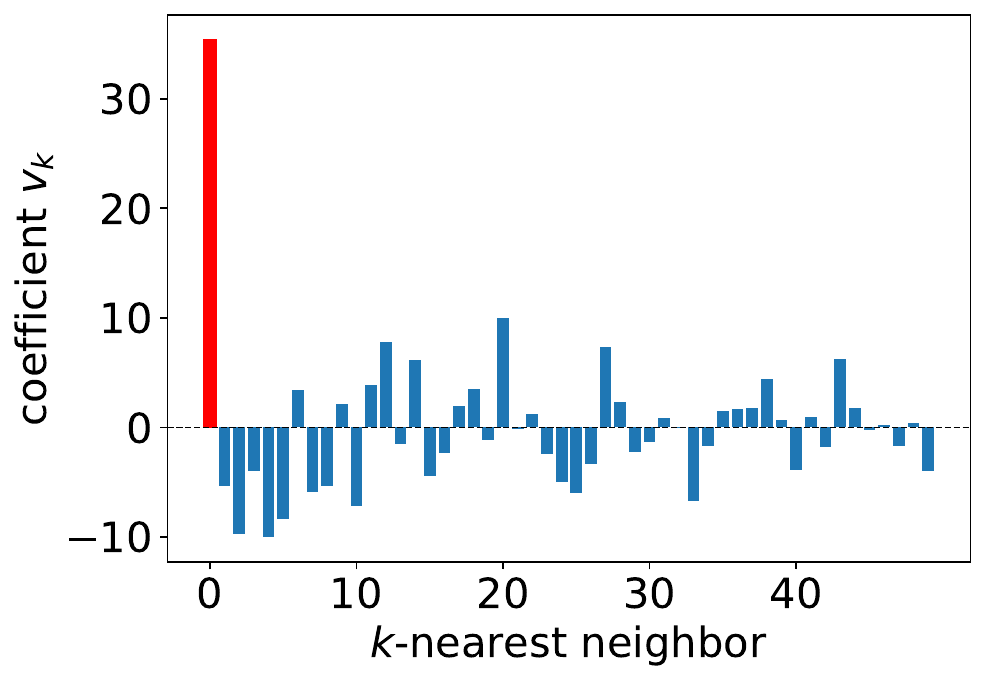}
			\end{subfigure}
		\end{minipage}\hfill
		\begin{minipage}[t]{0.34\textwidth}
			\centering
			\begin{subfigure}{\linewidth}
				\centering
				\caption{ RBF-FD weight, $\gamma_{\mathrm{bd}}<3$ }
				\label{fig:stencil_c}
				\includegraphics[height=3.8cm]{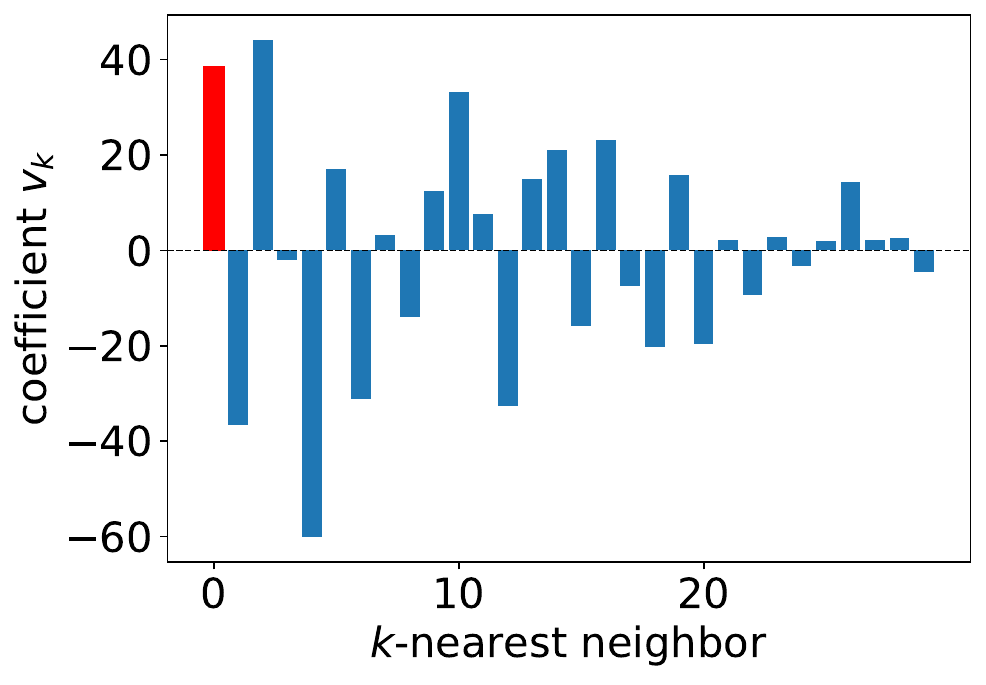}
			\end{subfigure}
			
			\vspace{3ex}
			\begin{subfigure}{\linewidth}
				\centering
				\caption{the same point as (c), using QP}
				\label{fig:stencil_d}
				\includegraphics[height=3.8cm]{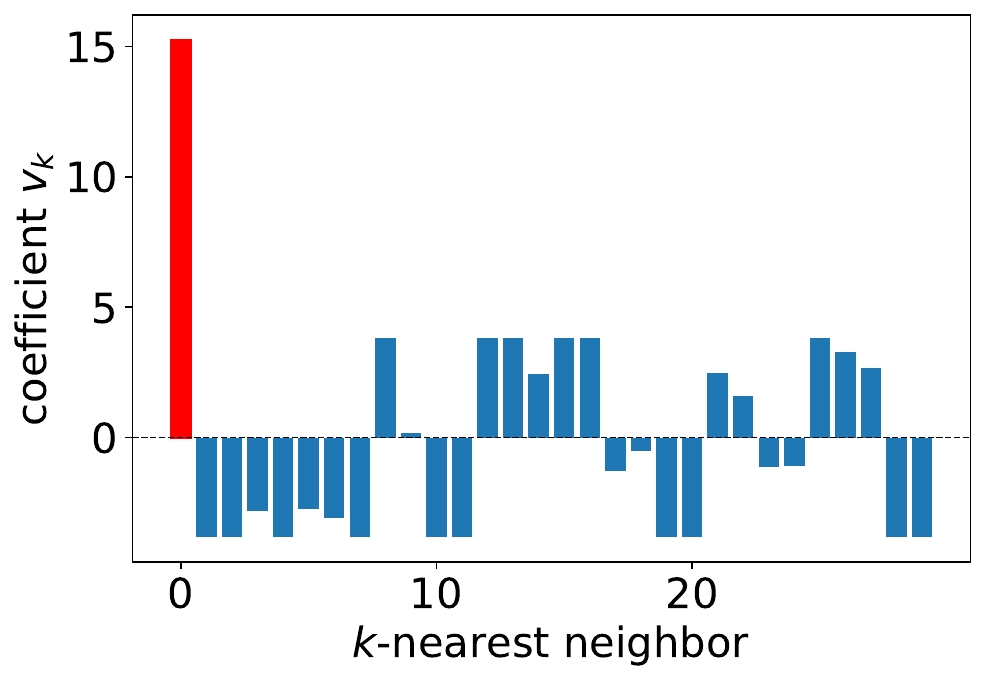}
			\end{subfigure}
		\end{minipage}\hfill
		\begin{minipage}[t]{0.28\textwidth}
			\centering
			\begin{subfigure}{\linewidth}
				\centering
				\caption{pointwise forward error}
				\label{fig:stencil_e}
				\includegraphics[height=3.8cm]{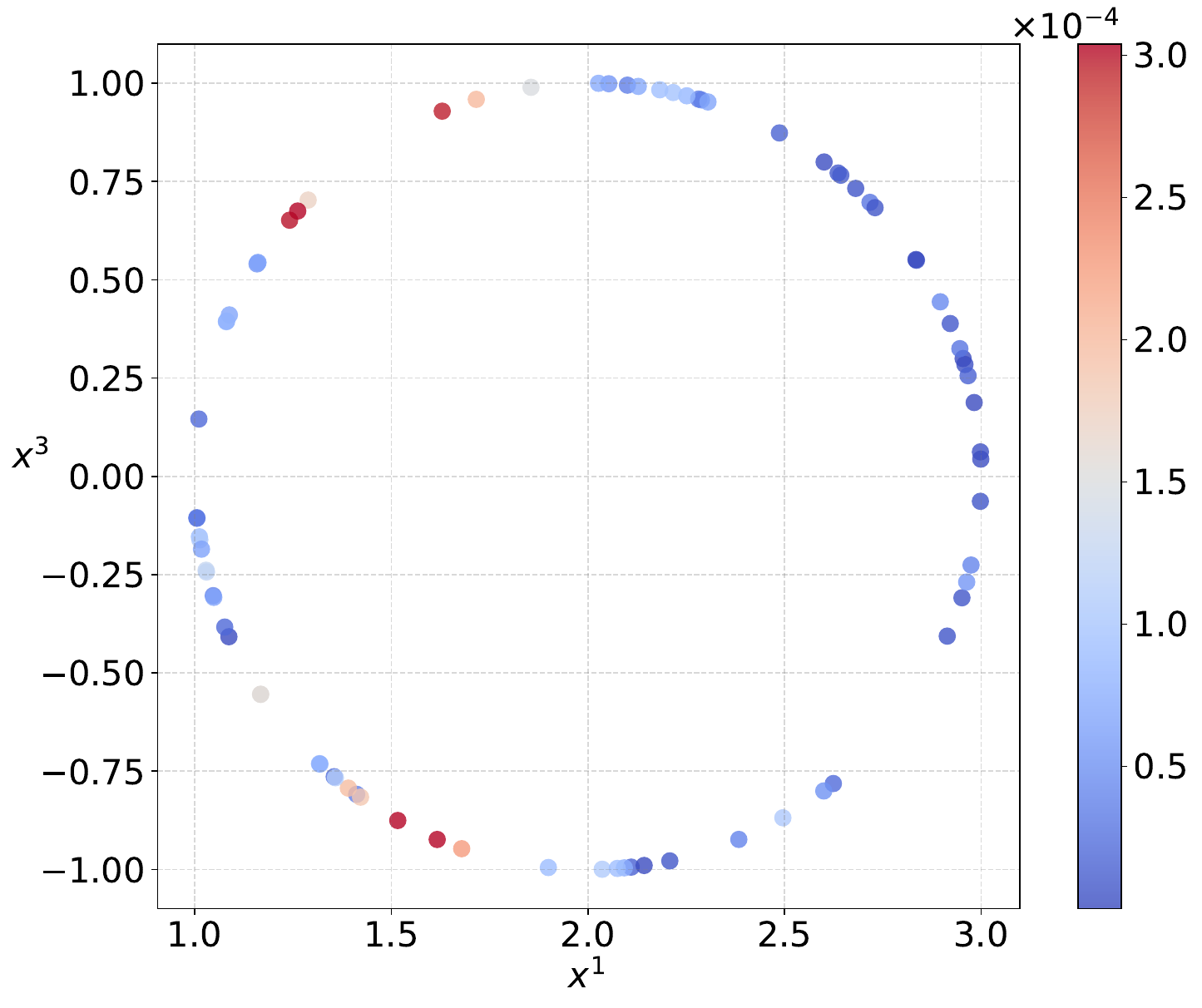}
			\end{subfigure}
			
			\vspace{3ex}
			\begin{subfigure}{\linewidth}
				\centering
				\caption{consistency}
				\label{fig:stencil_f}
				\includegraphics[height=3.8cm]{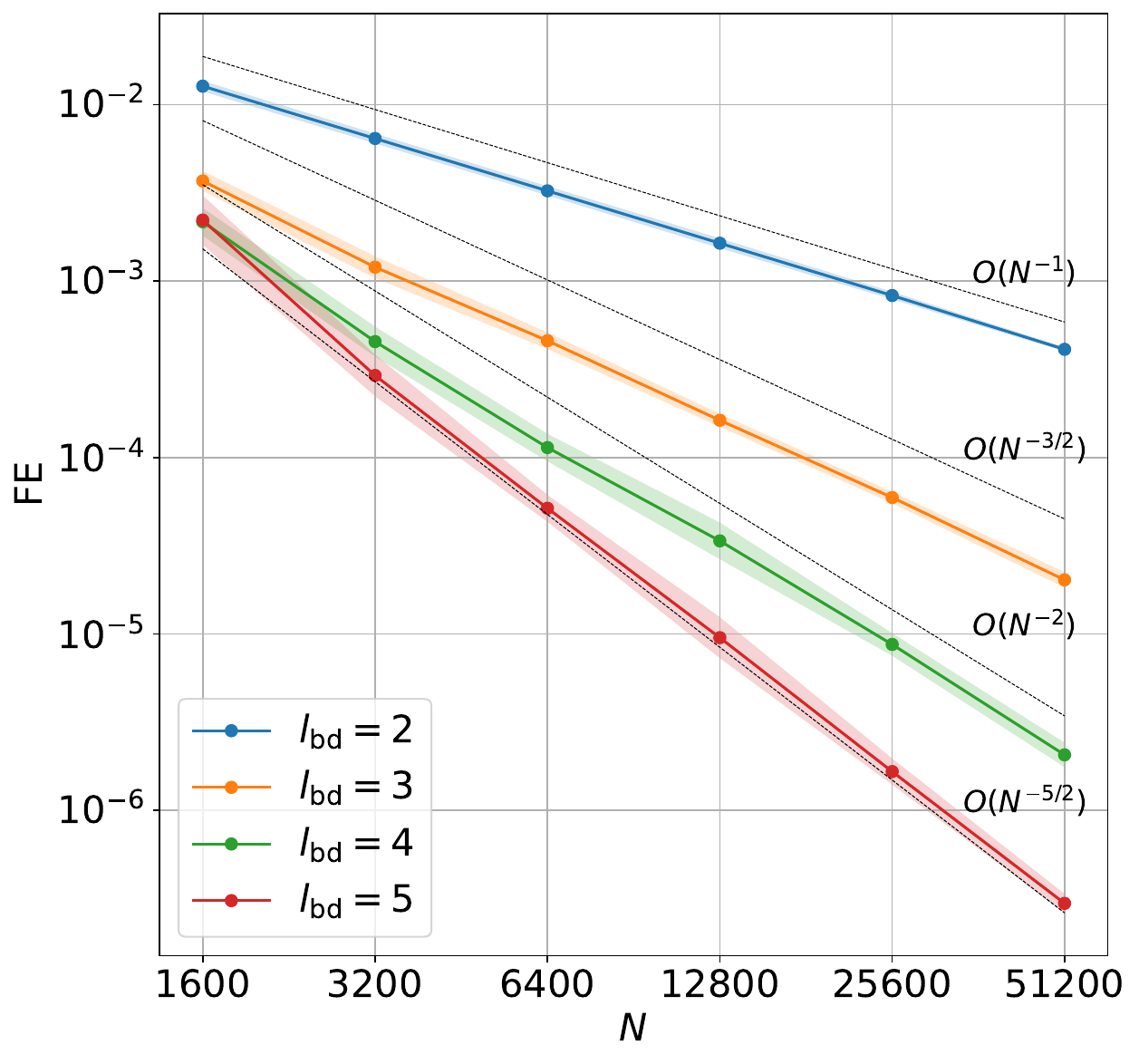}
			\end{subfigure}
		\end{minipage}
		
		\vspace{2ex}
		\caption{\textbf{Stencil illustration and numerical consistency for boundary operators.}
			Panel (a) shows the restricted KNN stencil. Panel (b) shows a typical distribution of weights obtained by the RBF-FD approach for most boundary points.
			Panels (c) and (d) compare the weights for the outward co-normal derivative without and with quadratic  optimization,
			demonstrating the elimination of unstable  weights.
			Panel (e) illustrates the pointwise forward error (\textbf{FE}) for a circular boundary of the semi-torus,
			while panel (f) shows the consistency across different polynomial degrees $l_{\mathrm{bd}}$,
			matching the expected theoretical convergence rate $O(N^{-l/2})$.}
		\label{fig:stencil_and_fe}
	\end{figure*}

Moreover, the restricted KNN search is based on a weighted metric biased in the outward co-normal direction. For a boundary base point $\mathbf{x}_b$, let $\boldsymbol{n}$ be the unit outward co-normal vector. The weighted distance metric is defined as:
\begin{equation}
	\Vert \mathbf{x}_{b,k} - \mathbf{x}_{b} \Vert_\mathbf{W} :=  \Vert \omega \boldsymbol{n}\boldsymbol{n}^\top(\mathbf{x}_{b,k} - \mathbf{x}_{b}) + (\mathbf{I} - \boldsymbol{n}\boldsymbol{n}^\top)(\mathbf{x}_{b,k} - \mathbf{x}_{b}) \Vert,
\end{equation}
where $\Vert \cdot \Vert$ is the standard Euclidean norm in $\mathbb{R}^3$ and $\mathbf{I}$ is the $3 \times 3$ identity matrix. The parameter $\omega \in (0, 1)$ is a scaling factor that intentionally shrinks the distance measurement along the co-normal direction. By taking small $\omega$, the effective search neighborhood forms an ellipsoid elongated along $\boldsymbol{n}$ (see Fig.~\ref{fig:stencil_and_fe}(a)). In our numerical implementation, we set $\omega = 1/3$. The stencil $S_{\mathbf{x}_{b}}=\{\mathbf{x}_{b,k} \in \mathbf{X}_{I} \cup \{\mathbf{x}_b\} \}_{k=1}^{K}$ is then efficiently constructed by querying the $K$-nearest neighbors based on this weighted metric.
	We also utilize the $1/K$ diagonal weight matrix $\boldsymbol{\Lambda }$ as defined
	in Eq. (\ref{eqn:k_weight}) for a stable approximation of the co-normal derivative (\ref{eqn:normal_weight}).
	
	
	Similar to the interior Laplacian approximation, we employ an auto-tuning algorithm
	to choose an appropriate $K$ at each boundary point. For the approximation of
	the outward co-normal derivative $\frac{\partial u}{\partial \boldsymbol{n}}(%
	\mathbf{x}_{b})\approx \sum_{k=1}^{K}v_{k}u(\mathbf{x}_{b,k})$, we consider
	the following two nearly diagonal dominance conditions:
	\begin{enumerate}
		\item (Sign constraint). The weight at the center node is strictly positive, i.e., $v_1 > 0$.		
		\item (Ratio constraint). The ratio $%
		\gamma_\mathrm{bd} := \displaystyle{\frac{|v_1|}{\max_{2 \leq k \leq K} |v_k|}}$ satisfies $%
		\gamma_\mathrm{bd} \geq 3$.
	\end{enumerate}
	This sign constraint is motivated by the fact that the derivative is taken in the outward co-normal direction. For a one-sided boundary stencil, the resulting discrete co-normal derivative typically has a positive coefficient at the center node.
	
	
	Algorithmically, the boundary routine reuses Algorithm~\ref{algo:intrin-LB_ad} from Sec.~\ref{sec:QP} with three modifications: the interior Laplacian weights are replaced by the co-normal derivative weights in \eqref{eqn:normal_weight}, the conditions \(v_1>0\) and \(\gamma_\mathrm{bd}\geq 3\) are used during the same auto-tuned search for \(K\), and the quadratic optimization formulation is given for the boundary base point by the following:
		\begin{equation}
			\begin{array}{ll}
				\underset{v_{1},\ldots ,v_{K},C}{\min } & \sum_{k=1}^{K}\displaystyle{\frac{1%
					}{2}\frac{v_{k}^{2}+C^{2}}{\lambda _{k}}} \\
				\mathrm{subject\text{ }to} & \sum_{k=1}^{K}v_{k}p_{\alpha (j)}(\boldsymbol{%
					\theta }\left( \mathbf{x}_{b,k}\right) )=\frac{\partial}{\partial \boldsymbol{n}} p_{\alpha (j)}(%
				\boldsymbol{\theta }\left( \mathbf{x}_{b}\right) ),\text{ \ \ }j=1,\ldots ,m,
				\\
				& v_{1}\geq 0, \\
				& -C \leq v_{k}\leq C,\text{ \ \ }k=2,\ldots ,K, \\
				& C\geq 0.%
			\end{array}
			\label{eqn:inequadbd}
		\end{equation}%
		Then the  weights $(v_1,\ldots,v_K)$ are inserted row-wise into the corresponding boundary operator matrix.
	
	Fig.~\ref{fig:stencil_and_fe}(b) shows a typical distribution of weights obtained by the RBF-FD approach for most boundary points satisfying the two nearly diagonal dominance conditions.  Figs.~\ref{fig:stencil_and_fe}(c) and (d) compare the weights for the outward co-normal derivative before and after quadratic  optimization, demonstrating the elimination of unstable  weights. To verify the validity of the boundary operator approximation, we define the
		forward error (\textbf{FE}) for the co-normal derivative similarly to (\ref{eqn:LapAp}). Fig. ~\ref{fig:stencil_and_fe}(f) demonstrates the consistency of \textbf{FE}s for
		different polynomial degrees $l_\mathrm{bd}$. The numerical consistency rate is approximately $%
		O(N^{-l_\mathrm{bd}/2})$, which is in good agreement with the theoretical prediction.
	
	
	\subsection{Matrix assembly and PDE solution}
	
	\label{sec:global_assembly}
	
	With the discretized coefficients for the interior Laplace--Beltrami operator
	(Algorithm~\ref{algo:intrin-LB_ad}) and the boundary co-normal derivative (Sec. \ref{sec:algo_boundary}) computed, the final step is to
	assemble the sparse linear system and solve the boundary value problem.
	
	To solve the Poisson problem $\Delta _{M}u=f$ on the surface $M$ subjected
	to the Robin boundary condition $u+\frac{\partial u}{\partial \boldsymbol{n}}%
	=h$ on the boundary $\partial M$ (Example {\ref{ex1}}), we discretize the
	differential operators to construct a coupled algebraic system. Let $\mathbf{u}_I$ and $\mathbf{u}_B$ denote the vectors of unknown values at the interior and boundary nodes, respectively. Then the discrete system consists of  equations for interior points obtained from  the Laplace--Beltrami matrix and equations for boundary points obtained from  the discretization of the derivative boundary conditions. The
	discretized Laplace--Beltrami matrix $\boldsymbol{L}_{I,M}\in \mathbb{R}%
	^{N_{I}\times N}$ is constructed for all interior points via the RBF-FD-QP
	method. We partition the Laplacian matrix $\boldsymbol{L}_{I,M}$ and the
	solution vector $\mathbf{u}\in \mathbb{R}^{N}$ as follows:
	\begin{equation*}
		\boldsymbol{L}_{I,M}=%
		\begin{bmatrix}
			\boldsymbol{L}_{I,I} & \boldsymbol{L}_{I,B}%
		\end{bmatrix}%
		,\quad \mathbf{u}=%
		\begin{bmatrix}
			\mathbf{u}_{I} \\
			\mathbf{u}_{B}%
		\end{bmatrix}%
		,
	\end{equation*}%
	where $\boldsymbol{L}_{I,I}\in \mathbb{R}^{N_{I}\times N_{I}}$, $\boldsymbol{%
		L}_{I,B}\in \mathbb{R}^{N_{I}\times N_{B}}$, $\mathbf{u}_{I}\in \mathbb{R}%
	^{N_{I}\times 1}$, and $\mathbf{u}_{B}\in \mathbb{R}^{N_{B}\times 1}$. This
	yields a set of discretized equations for interior points:
	\begin{equation}
		\boldsymbol{L}_{I,I}\mathbf{u}_{I}+\boldsymbol{L}_{I,B}\mathbf{u}_{B}=%
		\mathbf{f}_{I},  \label{eqn:interior_discrete}
	\end{equation}%
	where $\mathbf{f}_{I}\in \mathbb{R}^{N_{I}}$ contains the forcing term $f$
	at the interior points.
	
	For the Robin boundary conditions, we define a discretized boundary operator
	matrix $\boldsymbol{B}\in \mathbb{R}^{N_{B}\times N}$ representing the
	discretization of the operator $(\mathcal{I}+\frac{\partial }{\partial
		\boldsymbol{n}})$, where $\mathcal{I}$ is the identity map. The boundary
	matrix $\boldsymbol{B}$\ is partitioned to be:
	\begin{equation*}
		\boldsymbol{B}=%
		\begin{bmatrix}
			\boldsymbol{B}_{B,I} & \boldsymbol{B}_{B,B}%
		\end{bmatrix}%
		,
	\end{equation*}%
	where $\boldsymbol{B}_{B,I}\in \mathbb{R}^{N_{B}\times N_{I}}$ and $%
	\boldsymbol{B}_{B,B}\in \mathbb{R}^{N_{B}\times N_{B}}$. Note
	that the matrix $\boldsymbol{B}_{B,B}$
	is diagonal with positive diagonal entries because we employ the restricted KNN strategy for the boundary stencil
	selection (as discussed in Sec. \ref{sec:algo_boundary}) as well as the sign constraint in the nearly dominance condition for the boundary base point.
	
	
	We then couple the discretized Poisson equation (\ref{eqn:interior_discrete}%
	) with the discretized boundary conditions to obtain the following block
	system:
	\begin{equation}
		\begin{bmatrix}
			\boldsymbol{L}_{I,I} & \boldsymbol{L}_{I,B} \\
			\boldsymbol{B}_{B,I} & \boldsymbol{B}_{B,B}%
		\end{bmatrix}%
		\begin{bmatrix}
			\mathbf{u}_{I} \\
			\mathbf{u}_{B}%
		\end{bmatrix}%
		=%
		\begin{bmatrix}
			\mathbf{f}_{I} \\
			\mathbf{h}_{B}%
		\end{bmatrix}%
		,  \label{eqn:LB}
	\end{equation}%
	where $\mathbf{h}_{B}\in \mathbb{R}^{N_{B}\times 1}$ is the nonhomogeneous
	boundary data vector. To solve this coupled system (\ref{eqn:LB})
	efficiently, we eliminate the boundary unknowns $\mathbf{u}_{B}$. Solving
	the second block row for $\mathbf{u}_{B}$, we obtain:
	\begin{equation}
		\mathbf{u}_{B}=\boldsymbol{B}_{B,B}^{-1}(\mathbf{h}_{B}-\boldsymbol{B}_{B,I}%
		\mathbf{u}_{I}).  \label{eqn:uB}
	\end{equation}%
	Substituting this expression into the first block row yields a reduced
	system for the interior unknowns $\mathbf{u}_{I}$ only:
	\begin{equation}
		\underbrace{(\boldsymbol{L}_{I,I}-\boldsymbol{L}_{I,B}\boldsymbol{B}%
			_{B,B}^{-1}\boldsymbol{B}_{B,I})}_{\boldsymbol{A}'}\mathbf{u}_{I}=\underbrace{%
			\mathbf{f}_{I}-\boldsymbol{L}_{I,B}\boldsymbol{B}_{B,B}^{-1}\mathbf{h}_{B}}_{%
			\mathbf{b}'}.  \label{eqn:schur}
	\end{equation}%
	By solving the reduced $N_{I}\times N_{I}$ system $\boldsymbol{A}'\mathbf{u}%
	_{I}=\mathbf{b}'$, we obtain the solution $\mathbf{u}_{I}$\ at the interior
	points. The boundary values $\mathbf{u}_{B}$ are then recovered using
	back-substitution in (\ref{eqn:uB}).
	
	To numerically verify the overall accuracy and convergence of our proposed
	approach for solving the boundary value problem, we evaluate the inverse
	error (\textbf{IE}):
	\begin{equation*}
		\mathbf{IE}=\left( \frac{1}{N}\sum_{i=1}^{N}(\mathbf{u}_{i}-u(\mathbf{x}%
		_{i}))^{2}\right) ^{1/2},
	\end{equation*}%
	where $u(\mathbf{x}_{i})$ is the exact analytic solution evaluated at node $%
	\mathbf{x}_{i}$ and $\mathbf{u}_{i}$ is the corresponding numerical solution.
	
\begin{figure*}[tbp]
		\centering
		\begin{minipage}[t]{0.32\textwidth}
			\centering
			\begin{subfigure}{\linewidth}
				\centering
				\caption{pointwise error, $l=3,l_{\mathrm{bd}}=2$}
				\includegraphics[height=4.4cm]{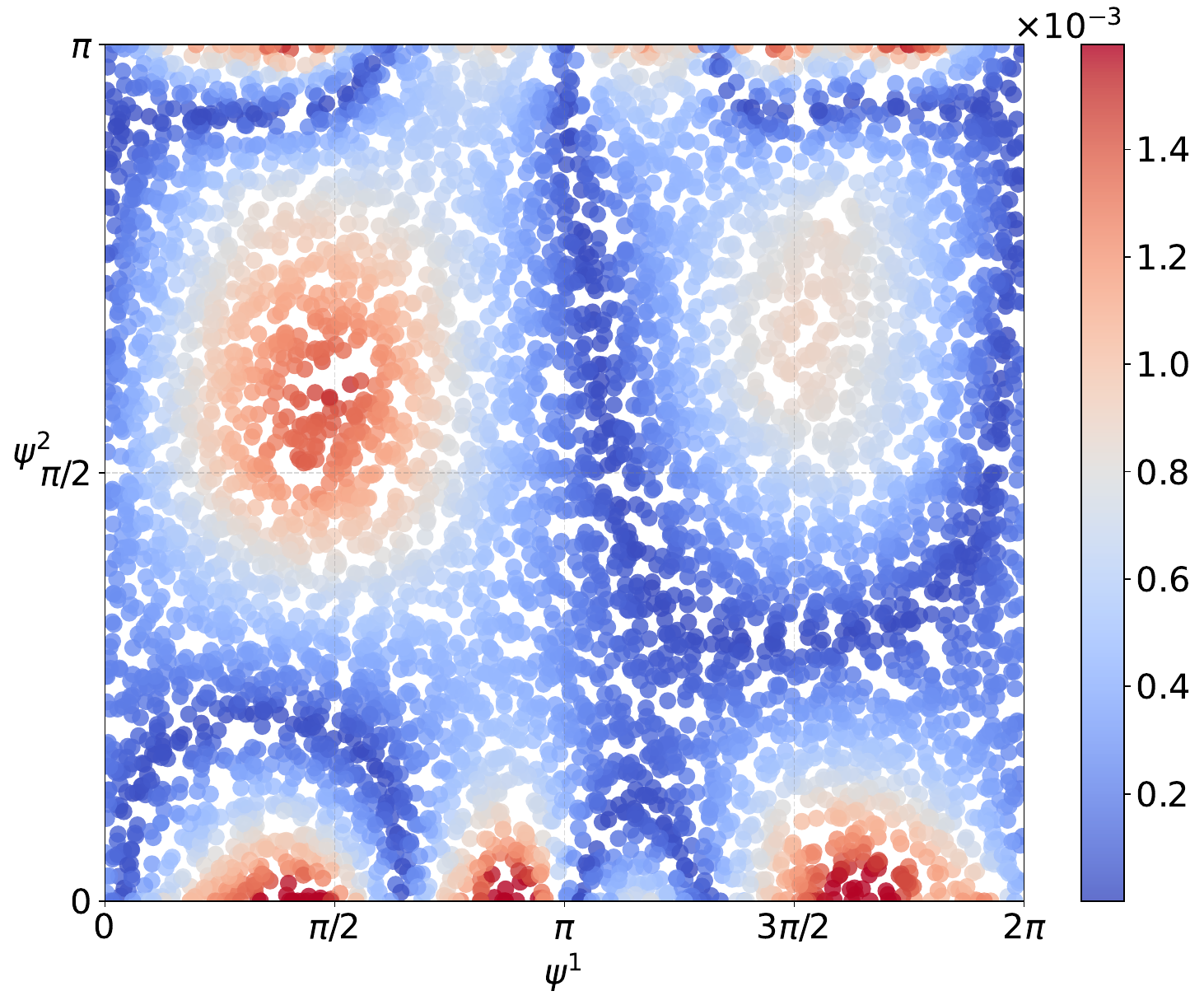}
			\end{subfigure}
			
			\vspace{3ex} 
			\begin{subfigure}{\linewidth}
				\centering
				\caption{pointwise error, $l=4,l_{\mathrm{bd}}=4$}
				\includegraphics[height=4.4cm]{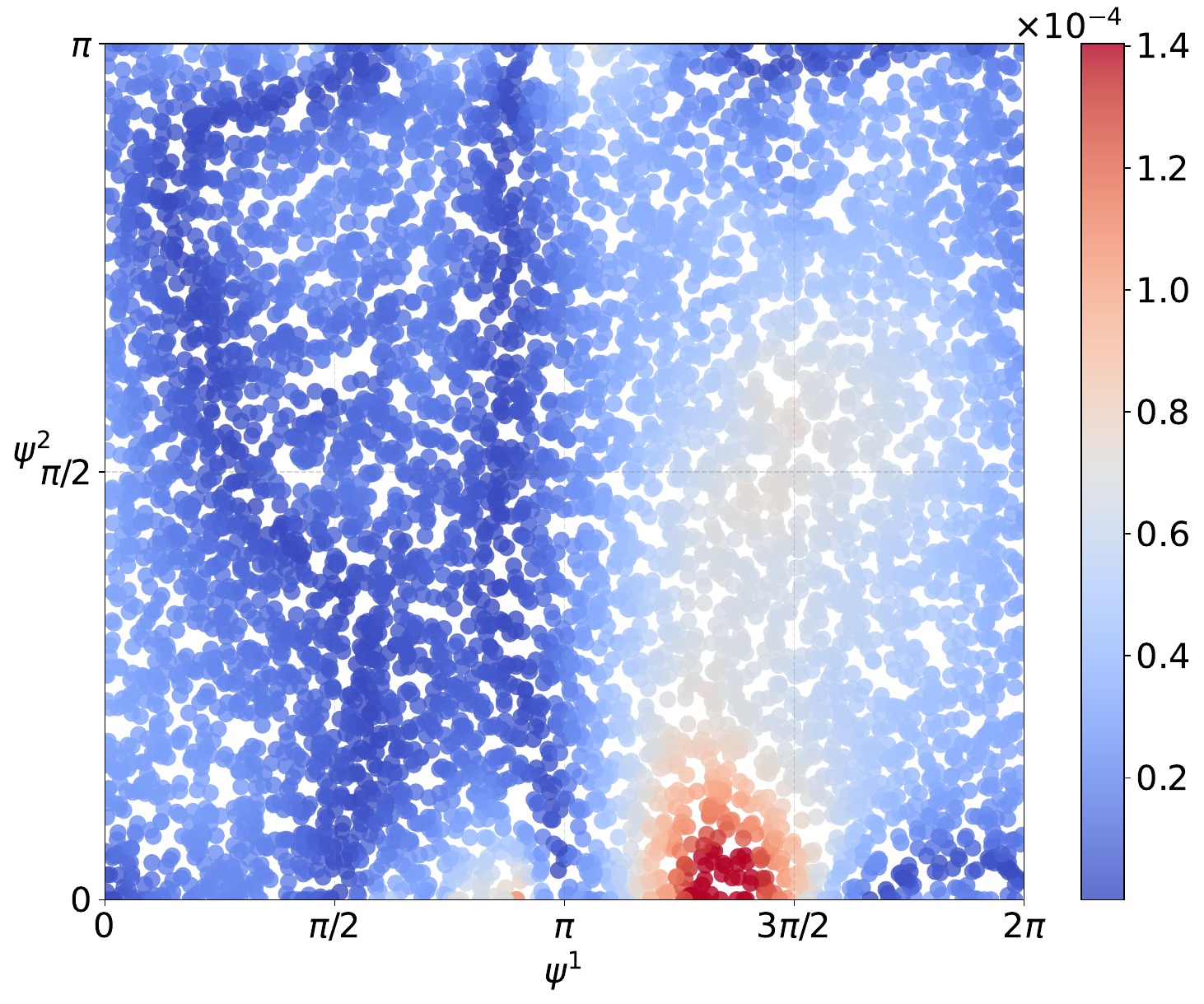}
			\end{subfigure}
		\end{minipage}\hfill
		\begin{minipage}[t]{0.32\textwidth}
			\centering
			\begin{subfigure}{\linewidth}
				\centering
				\caption{\textbf{IE} for $l=2$}
				\includegraphics[height=4.4cm]{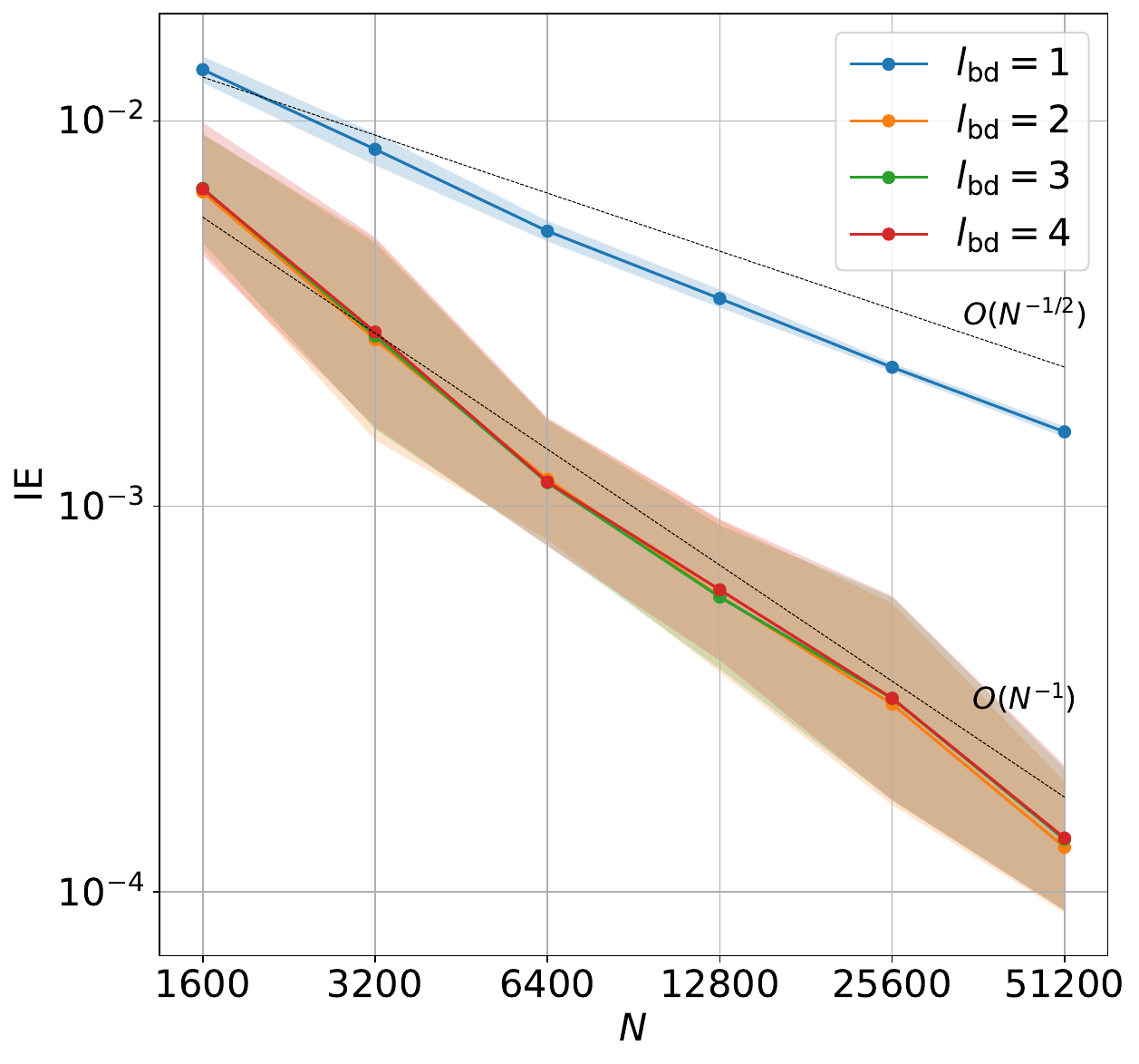}
			\end{subfigure}
			
			\vspace{3ex}
			\begin{subfigure}{\linewidth}
				\centering
				\caption{\textbf{IE} for $l=3$}
				\includegraphics[height=4.4cm]{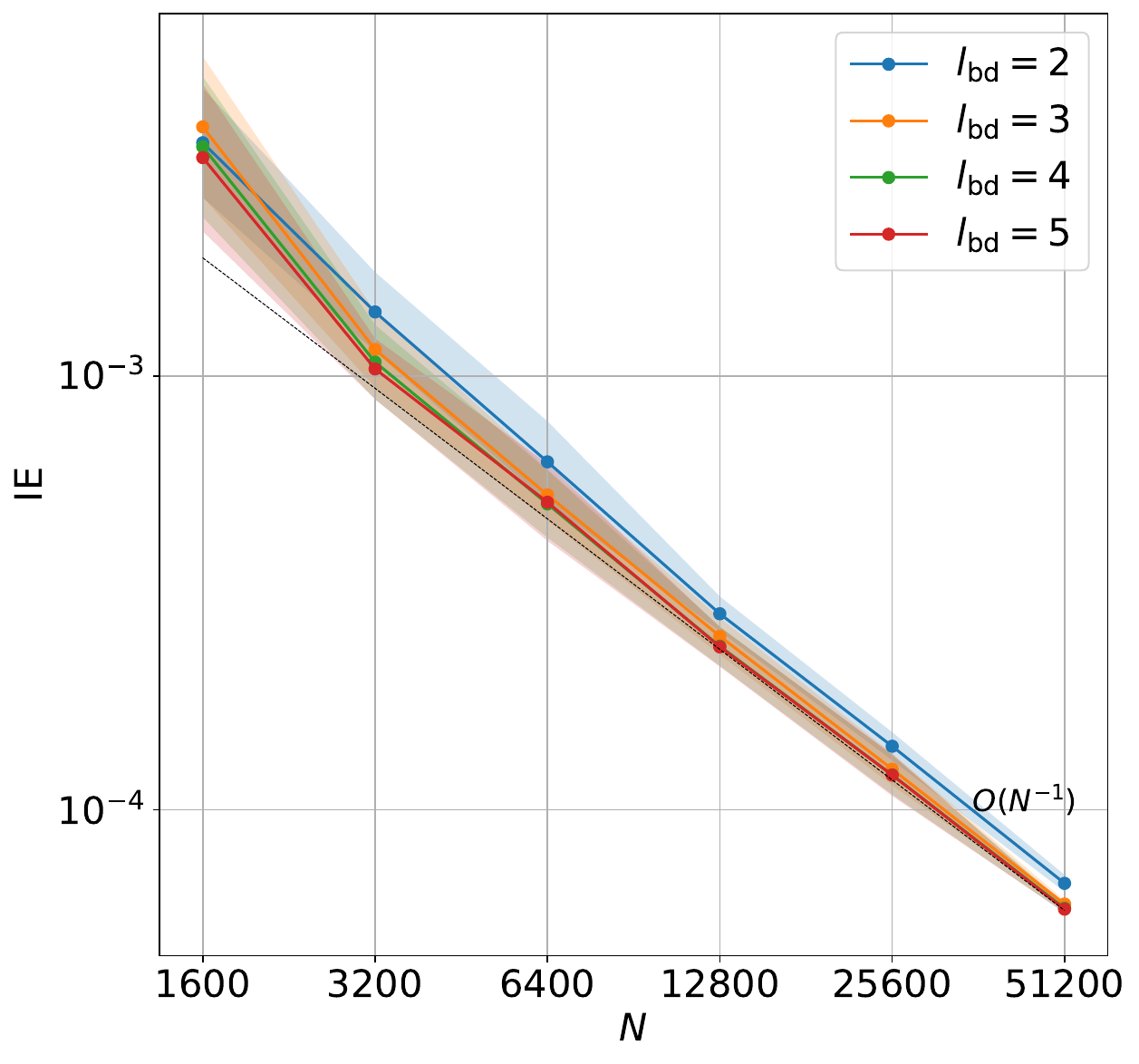}
			\end{subfigure}
		\end{minipage}\hfill
		\begin{minipage}[t]{0.32\textwidth}
			\centering
			\begin{subfigure}{\linewidth}
				\centering
				\caption{\textbf{IE} for $l=4$}
				\includegraphics[height=4.4cm]{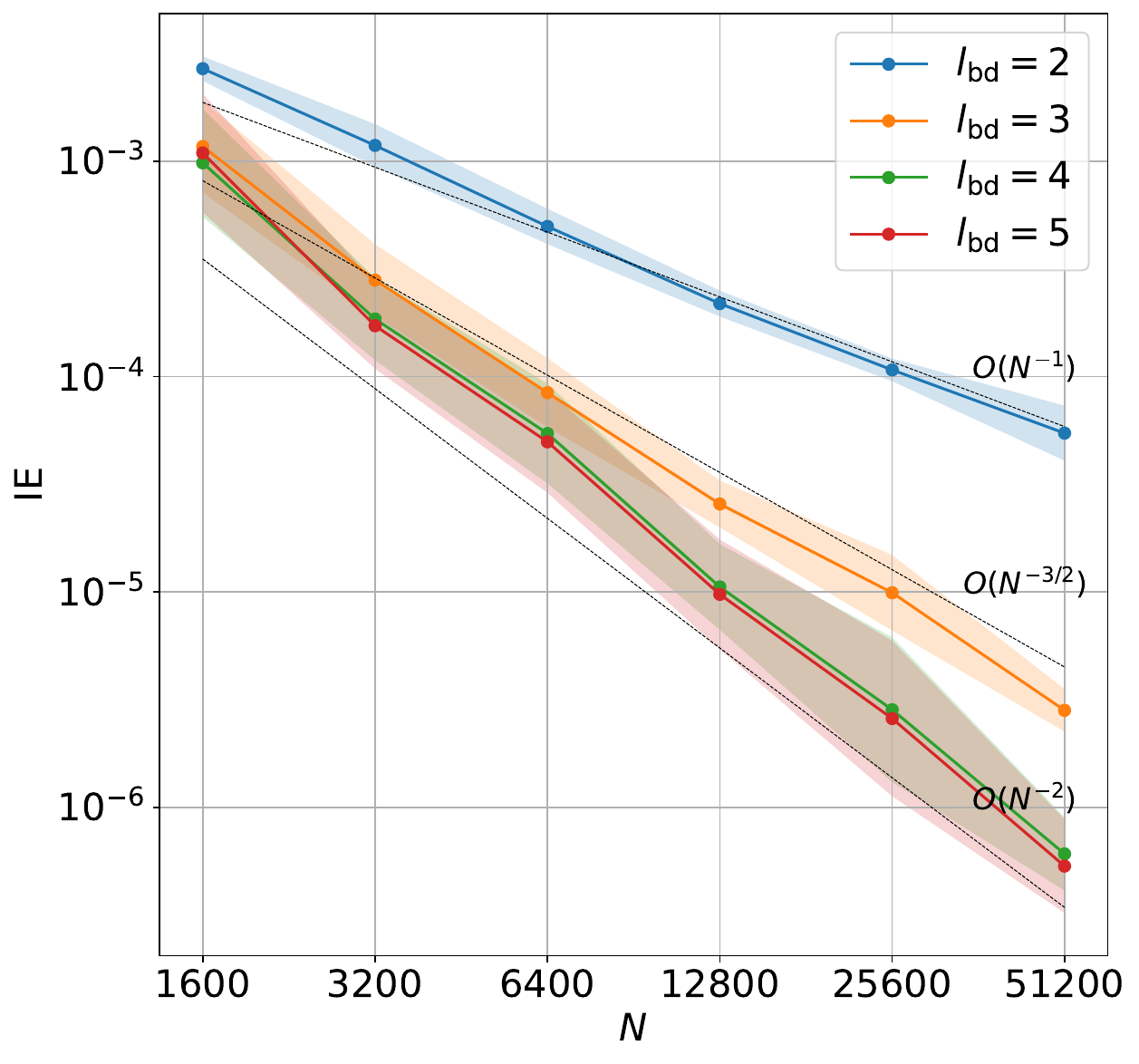}
			\end{subfigure}
			
			\vspace{3ex}
			\begin{subfigure}{\linewidth}
				\centering
				\caption{\textbf{IE} for $l=5$}
				\includegraphics[height=4.4cm]{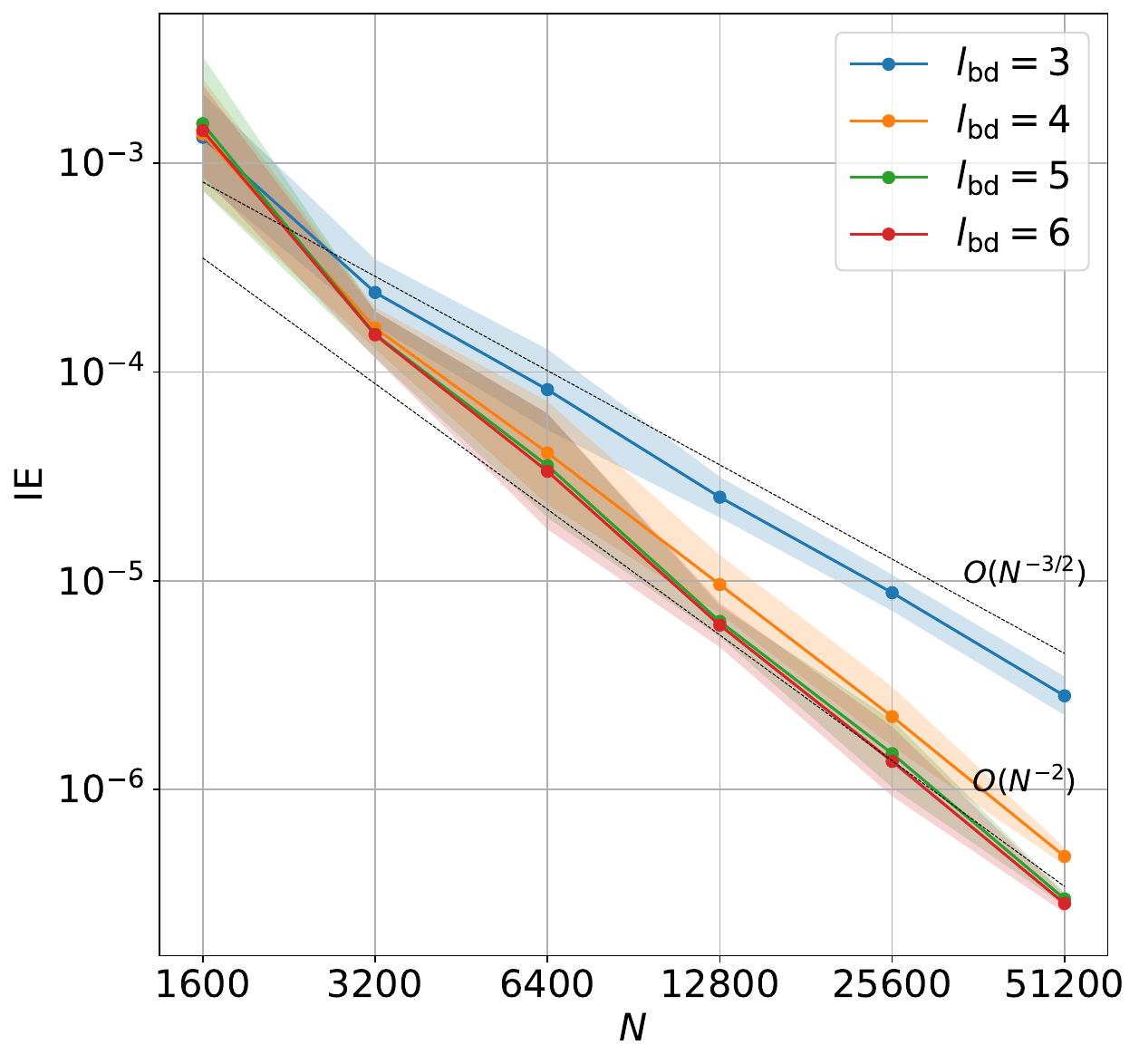}
			\end{subfigure}
		\end{minipage}
		
		\vspace{2ex}
		\caption{\textbf{Solution accuracy for Poisson problems on a 2D semi-torus in $\mathbb{R}^{3}$.}
			Panels (a) and (b) display the spatial distribution of pointwise absolute errors for degrees of $(l,l_{\mathrm{bd}})=(3,2)$ and $(l,l_{\mathrm{bd}})=(4,4)$, respectively.
			Panels (c)--(f) illustrate the convergence of inverse errors (\textbf{IE}s) for $l \in \{2, 3, 4, 5\}$. All simulations are run with 12 independent trials, each with a set of randomly sampled data points.}
		\label{ie_torus}
	\end{figure*}

	{Fig.~\ref{ie_torus} demonstrates the numerical performance of the RBF-FD-QP
		solver for the Poisson equation on the 2D semi-torus in Example %
		\ref{ex1}. Figs.~\ref{ie_torus}(a) and \ref{ie_torus}(b) illustrate the spatial distribution of the pointwise
		absolute error, while Figs.~\ref{ie_torus}(c)-(f) show the  convergence of the
		inverse error (\textbf{IE}) under different polynomial degrees $l$ for interior points and
		the polynomial degrees $l_{\mathrm{bd}}$ for boundary points. Here, the interior degree $l$ and boundary degree $l_{\mathrm{bd}}$ should be chosen to match each other. It can be observed from Fig.~\ref{ie_torus} that when the interior degree $l$ is odd, the optimal boundary degree is $l_{\mathrm{bd}}=l-1$. When the interior degree $l$ is even, the optimal boundary degree is $l_{\mathrm{bd}}=l$. In the following section for numerical results, we only focus on the numerical performance in the setting of the pairs of $(l,l_{\mathrm{bd}})=(3,2)$ and $(l,l_{\mathrm{bd}})=(4,4)$.}

	\begin{figure*}[tbp]
		\centering
		\begin{subfigure}{0.32\textwidth}
			\centering
			\caption{\textbf{FE}}
			\label{fig:compare_fe}
			\includegraphics[height=5cm]{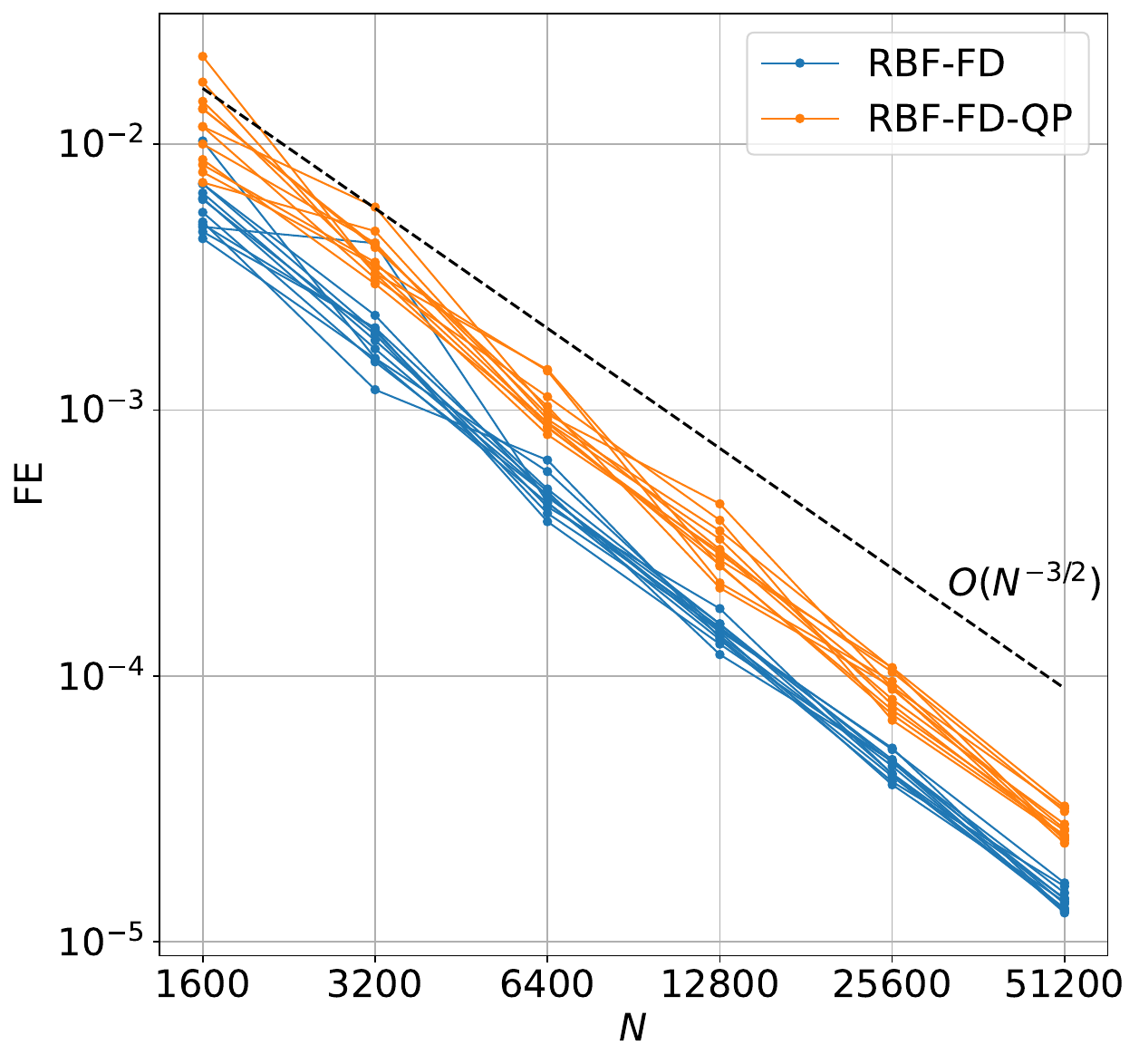}
		\end{subfigure}\hfill
		\begin{subfigure}{0.32\textwidth}
			\centering
			\caption{Stability}
			\label{fig:compare_st}
			\includegraphics[height=5cm]{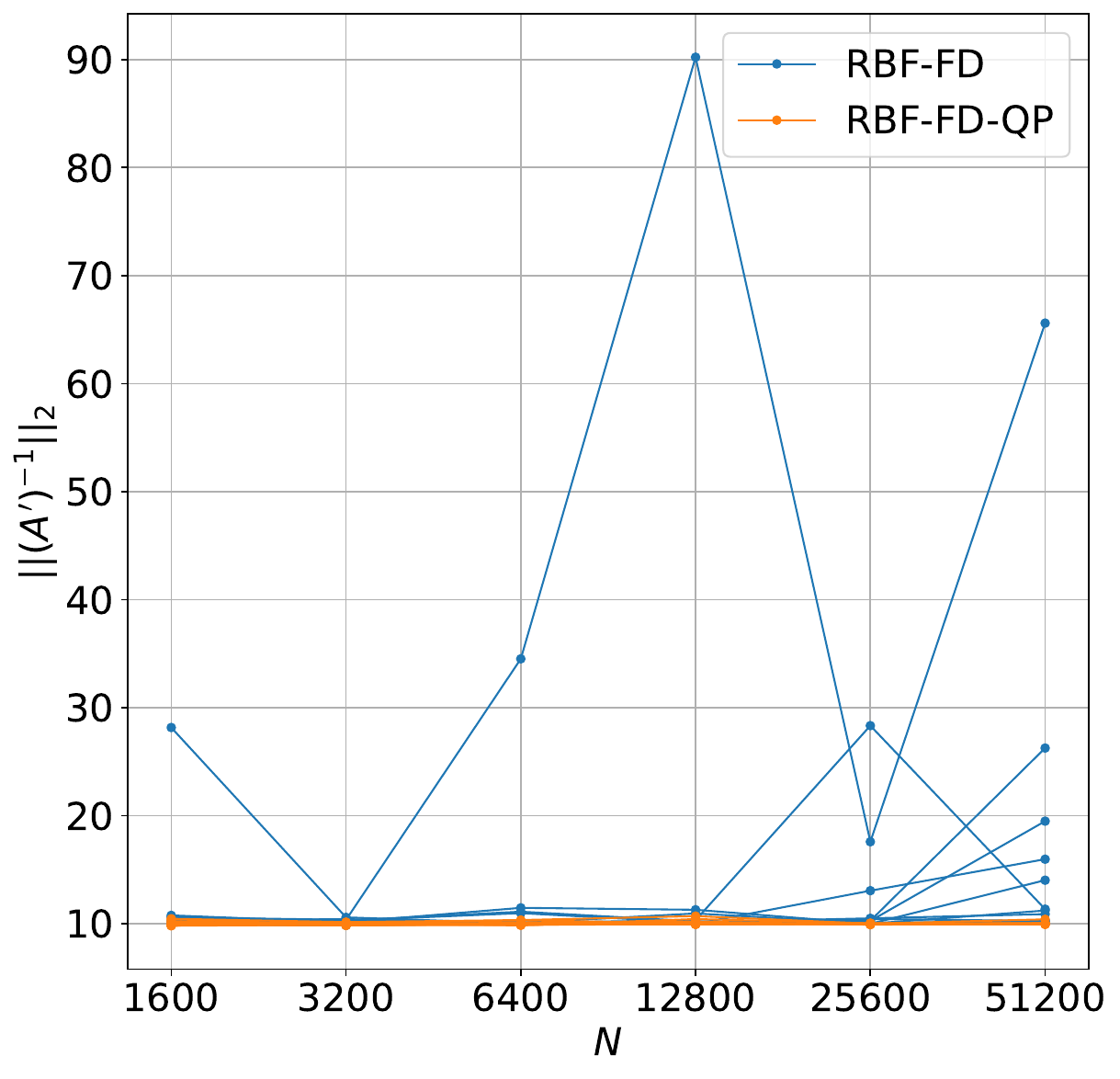}
		\end{subfigure}\hfill
		\begin{subfigure}{0.32\textwidth}
			\centering
			\caption{\textbf{IE}}
			\label{fig:compare_ie}
			\includegraphics[height=5cm]{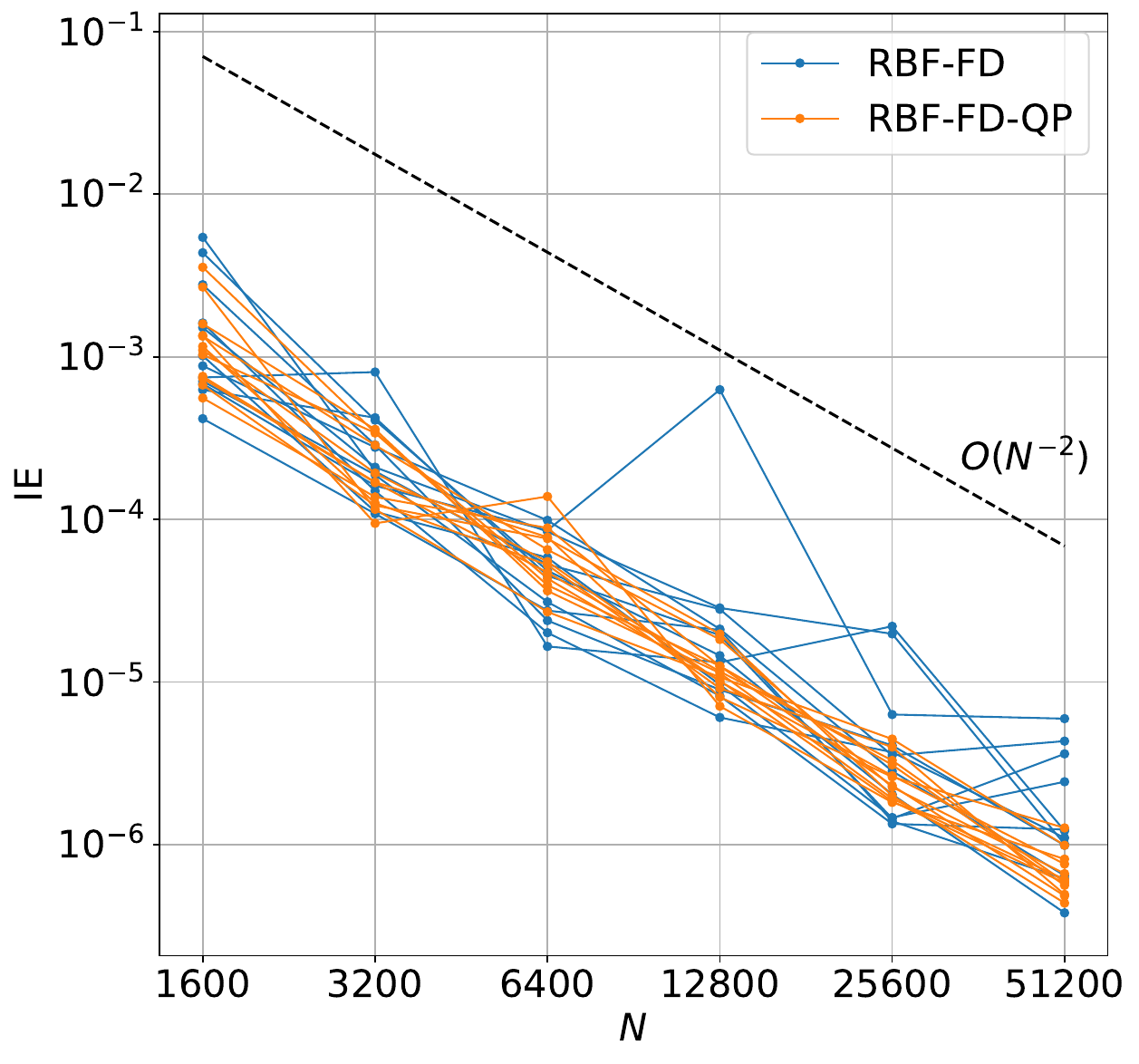}
		\end{subfigure}
		
		\vspace{2ex}
		\caption{\textbf{Ablation study on the necessity of the QP stabilization.} Comparison between the standard RBF-FD and the proposed RBF-FD-QP method for solving the Poisson problem on the 2D semi-torus. Panels (a), (b), and (c) display the forward error (\textbf{FE}), the stability measured by the spectral norm of the inverse reduced matrix $\|(\boldsymbol{A}')^{-1}\|_2$, and the inverse error (\textbf{IE}), respectively. Results from 12 independent random-sampling trials are plotted directly (without averaging or statistics) for visualization.}
		\label{fig:ablation}
	\end{figure*}
	
	
	\subsection{Comparison of RBF-FD and RBF-FD-QP} \label{sec:study_QPstab}
	
		Before proceeding to more complex numerical experiments, we  demonstrate the necessity and effectiveness of the proposed quadratic-programming-based stabilization. We conduct an ablation study comparing the standard RBF-FD method with our proposed RBF-FD-QP approach using the semi-torus Poisson problem (Example \ref{ex1}).
		
		To comprehensively evaluate the effect of the QP module, we compare three quantities: the interior forward error (\textbf{FE}) of the Laplace-Beltrami operator approximation, the stability of the reduced interior system characterized by the spectral norm of the inverse matrix $\|(\boldsymbol{A}')^{-1}\|_2$, and the inverse error (\textbf{IE}) of the PDE solution. For comparison, 12 independent trials with different randomly sampled point clouds are conducted, and the result of each trial is directly plotted for visualization (see Fig.~\ref{fig:ablation}).
		
		
		
		The standard RBF-FD method provides a reasonable and convergent approximation of the Laplacian operator (blue curves in Fig.~\ref{fig:ablation}(a)). However, the spectral norm of the inverse matrix $\|(\boldsymbol{A}')^{-1}\|_2$ obtained with the standard RBF-FD method may not be bounded around a constant for several trials (blue curves in Fig.~\ref{fig:ablation}(b)), implying that the reduced linear system (\ref{eqn:schur}) may lose stability. Meanwhile, for the trials with relatively large norms, the corresponding inverse errors fail to decay as $N$ increases (blue curves in Fig.~\ref{fig:ablation}(c)), indicating a loss of convergence. Upon inspection, we find that such instability and non-convergence are induced by the inappropriate local RBF-FD approximation to differential operators for some points near the boundary (see Fig.~\ref{fig2:torus}). This motivates us to propose the RBF-FD-QP approach to improve the coefficients and stabilize the system. As can be seen from Fig.~\ref{fig:ablation}(a), the proposed RBF-FD-QP method also provides a consistent approximation of the Laplacian operator, although its error is slightly larger than that of the standard RBF-FD. This behavior is expected, since the QP formulation imposes only equality constraints on polynomials, yielding less accurate coefficients than the standard RBF-FD, which enforces exactness for both polynomials and polyharmonic splines. More importantly, compared with the standard RBF-FD method, the proposed RBF-FD-QP method provides uniform bounds for the norm of the inverse matrix and, consequently, convergent solutions for the PDE across all random samplings (see orange curves in Figs.~\ref{fig:ablation}(b) and (c)).

	\section{Numerical experiments}
	\label{sec:numerical}
	
	In this section, we evaluate the proposed RBF-FD-QP approach across three types of problems on surfaces. In Sec. \ref{sec:eigen}, we consider the eigenvalue problems, in which we examine the convergence of the eigenvalues of the Laplace--Beltrami operator on surfaces with boundary. In Sec. \ref{sec:diffusion}, we evaluate the performance for solving linear time-dependent diffusion equations and verify that the proposed spatial discretization remains stable and accurate. In Sec. \ref{sec:interface}, we focus on elliptic interface problems on surfaces to demonstrate the ability of the proposed framework to handle solution discontinuities and jump conditions across the interface using point cloud data. To test the robustness of our approach, all numerical experiments are conducted on various independent trials of randomly sampled point clouds. The results demonstrate the ability of the RBF-FD-QP approach to maintain stability and high-order accuracy without relying on regular meshes or introducing ghost points. For data generation in all examples, the intrinsic data are first randomly generated from a uniform  sampling in the parameter space and then  mapped onto the surfaces in the ambient space via the parametrization. For the RBF-FD-QP approach, we set the PHS parameter $\kappa=3$ and the polynomial degrees $(l,l_{\mathrm{bd}})$ to be $(3,2)$ and $(4,4)$ across all experiments.
	
	
	
	
	\subsection{Eigenvalue problem} \label{sec:eigen}
	We first consider the eigenvalue problem for the Laplace--Beltrami operator on a two-dimensional surface $M$ embedded in $\mathbb{R}^3$, which is governed by:
	\begin{equation}
		\begin{aligned}
			-\Delta_M u_k &= \lambda_k u_k, \quad &\text{on}\ M, \\
			u_k + \frac{\partial u_k}{\partial \boldsymbol{n}} &= 0, \quad &\text{on} \ \partial M,
		\end{aligned}
	\end{equation}
	where $
		\{\lambda_k\}_{k=1}^{\infty}$ denote the eigenvalues and $\{u_k\}_{k=1}^{\infty}$ the corresponding eigenfunctions.
	We test the proposed RBF-FD-QP method on two smooth surfaces with boundary: a semi-torus (as defined in \eqref{eqn:torpar}) and a semi-sphere parameterized by $\mathbf{x}:=(x^1,x^2,x^3)=(\sin \psi^1 \cos \psi^2, \sin \psi^1 \sin \psi^2, \cos \psi^1)$ for $0 \leq \psi^1 \leq \frac{\pi}{2}$ and $0 \leq \psi^2 < 2\pi$.
	
	To validate the accuracy of the RBF-FD-QP approach, we compare our results against reference eigenvalues and eigenfunctions computed by the Finite Element Method (FEM) via the FreeFem++ software \cite{MR3043640}. For the FEM reference, the parametric domains are uniformly discretized into a $300 \times 300$ structured grid, utilizing continuous piecewise quadratic finite elements (P2) to construct the functional spaces.
	
	
	\begin{table*}[htbp]
		\centering
		\caption{Convergence of the mean and standard deviation of absolute errors for the selected eigenvalues on the semi-sphere and semi-torus. Here the polynomial degrees $(l,l_{\mathrm{bd}})=(4,4)$.}
		\label{tab:eigenv_convergence}
		
		\begin{subtable}{\textwidth}
			\centering
			\caption{Semi-sphere}
			\label{tab:eigenv_sphere}
			\resizebox{\textwidth}{!}{
				\begin{tabular}{l c c c c c c c c c c}
					\toprule
					\multirow{2}{*}{$N$} & \multicolumn{2}{c}{$k=1$} & \multicolumn{2}{c}{$k=2$} & \multicolumn{2}{c}{$k=4$} & \multicolumn{2}{c}{$k=8$} & \multicolumn{2}{c}{$k=20$} \\
					\cmidrule(lr){2-3} \cmidrule(lr){4-5} \cmidrule(lr){6-7} \cmidrule(lr){8-9} \cmidrule(lr){10-11}
					& Error mean (std) & Rate & Error mean (std) & Rate & Error mean (std) & Rate & Error mean (std) & Rate & Error mean (std) & Rate \\
					\midrule
					6400   & 2.277e-05 (2.699e-06) & -    & 3.714e-05 (1.397e-05) & -    & 9.647e-05 (1.032e-05) & -    & 6.333e-04 (1.443e-04) & -    & 1.670e-03 (4.465e-04) & -    \\
					12800  & 7.600e-06 (6.807e-07) & 1.58 & 3.072e-06 (1.329e-06) & 3.60 & 3.096e-05 (2.526e-06) & 1.64 & 1.219e-04 (1.704e-05) & 2.38 & 3.935e-04 (1.077e-04) & 2.09 \\
					25600  & 2.227e-06 (1.612e-07) & 1.77 & 3.306e-07 (2.479e-07) & 3.22 & 8.844e-06 (6.389e-07) & 1.81 & 2.425e-05 (2.848e-06) & 2.33 & 8.249e-05 (4.300e-05) & 2.25 \\
					51200  & 6.052e-07 (5.374e-08) & 1.88 & 3.449e-07 (5.693e-08) & -0.06& 2.362e-06 (1.936e-07) & 1.90 & 4.354e-06 (4.435e-07) & 2.48 & 2.294e-05 (7.569e-06) & 1.85 \\
					102400 & 1.685e-07 (1.511e-08) & 1.84 & 1.182e-07 (1.885e-08) & 1.54 & 6.477e-07 (5.592e-08) & 1.87 & 8.708e-07 (7.673e-08) & 2.32 & 9.793e-06 (2.208e-06) & 1.23 \\
					\midrule
					FEM Ref    & 0.7055983          & -    & 3.156904          & -    & 7.153786          & -    & 13.24604          & -    & 32.31781          & -    \\
					\bottomrule
				\end{tabular}
			}
		\end{subtable}
		
		\vspace{4ex}
		
		\begin{subtable}{\textwidth}
			\centering
			\caption{Semi-torus}
			\label{tab:eigenv_torus}
			\resizebox{\textwidth}{!}{
				\begin{tabular}{l c c c c c c c c c c}
					\toprule
					\multirow{2}{*}{$N$} & \multicolumn{2}{c}{$k=1$} & \multicolumn{2}{c}{$k=2$} & \multicolumn{2}{c}{$k=4$} & \multicolumn{2}{c}{$k=8$} & \multicolumn{2}{c}{$k=20$} \\
					\cmidrule(lr){2-3} \cmidrule(lr){4-5} \cmidrule(lr){6-7} \cmidrule(lr){8-9} \cmidrule(lr){10-11}
					& Error mean (std) & Rate & Error mean (std) & Rate & Error mean (std) & Rate & Error mean (std) & Rate & Error mean (std) & Rate \\
					\midrule
					6400   & 2.887e-06 (9.265e-07) & -    & 1.246e-05 (4.932e-06) & -    & 3.674e-05 (1.358e-05) & -    & 4.623e-05 (5.414e-05) & -    & 1.326e-04 (9.891e-05) & -    \\
					12800  & 5.234e-07 (2.707e-07) & 2.46 & 2.822e-06 (8.937e-07) & 2.14 & 1.017e-05 (1.652e-06) & 1.85 & 1.574e-05 (1.163e-05) & 1.55 & 2.868e-05 (1.646e-05) & 2.21 \\
					25600  & 1.580e-07 (7.903e-08) & 1.73 & 6.710e-07 (2.819e-07) & 2.07 & 2.651e-06 (7.814e-07) & 1.94 & 6.419e-06 (1.568e-05) & 1.29 & 2.003e-05 (1.318e-05) & 0.52 \\
					51200  & 2.775e-08 (1.528e-08) & 2.51 & 1.173e-07 (2.468e-08) & 2.52 & 6.993e-07 (2.187e-07) & 1.92 & 4.819e-07 (2.568e-07) & 3.74 & 5.279e-06 (1.143e-06) & 1.92 \\
					102400 & 5.859e-09 (2.531e-09) & 2.24 & 2.636e-08 (8.834e-09) & 2.15 & 1.721e-07 (3.936e-08) & 2.02 & 1.278e-07 (1.094e-07) & 1.91 & 1.541e-06 (2.378e-07) & 1.78 \\
					\midrule
					FEM Ref    & 0.1468539          & -    & 0.5501428          & -    & 1.137836          & -    & 2.268830          & -    & 5.571232          & -    \\
					\bottomrule
				\end{tabular}
			}
		\end{subtable}
	\end{table*}
	
	\begin{figure*}[tbp]
		\centering
		\begin{subfigure}{0.24\textwidth}
			\centering
			\caption{Mode 1 (Sphere)\\FEM}
			\includegraphics[width=\linewidth]{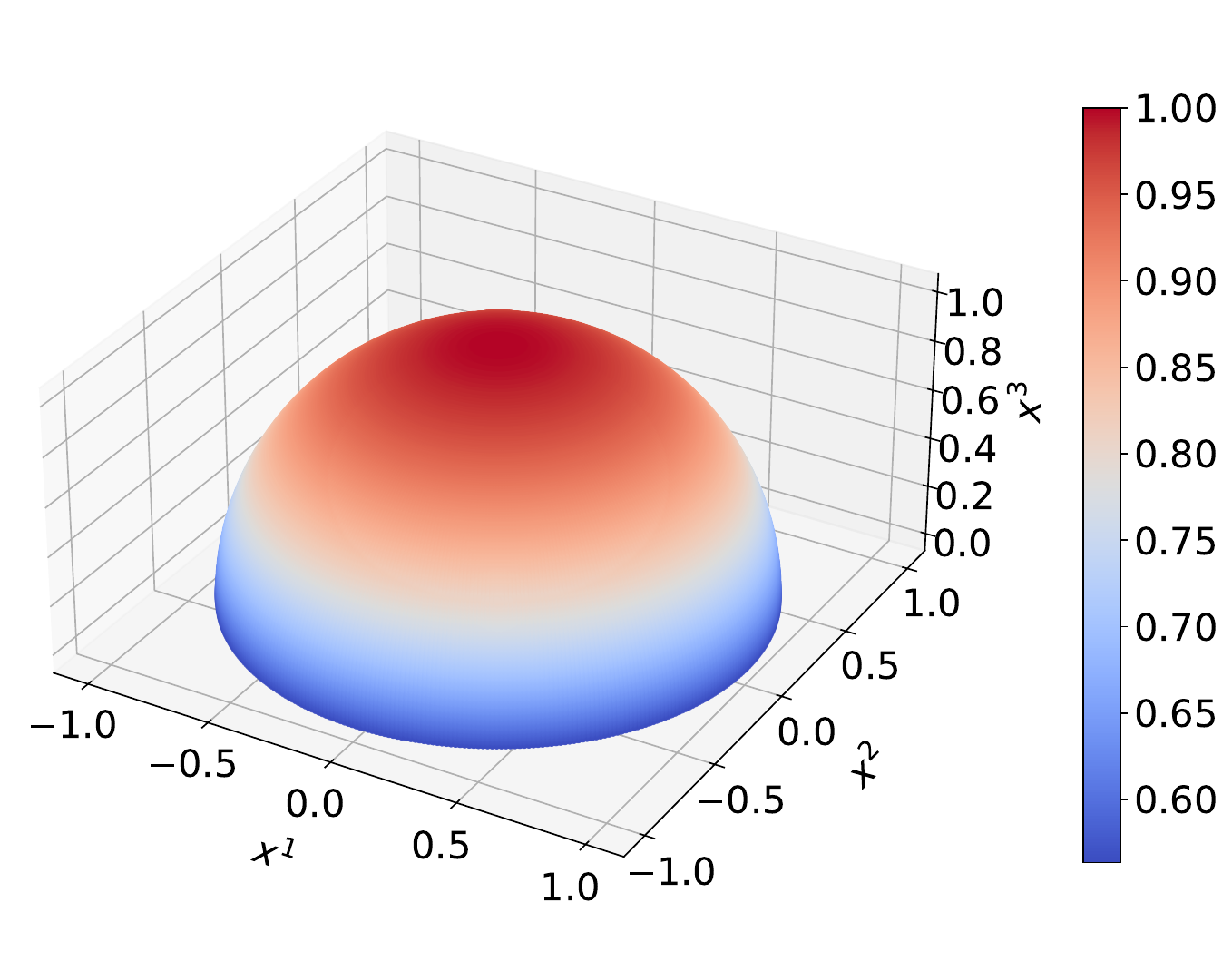}
		\end{subfigure}\hfill
		\begin{subfigure}{0.24\textwidth}
			\centering
			\caption{Mode 1 (Sphere)\\RBF-FD-QP}
			\includegraphics[width=\linewidth]{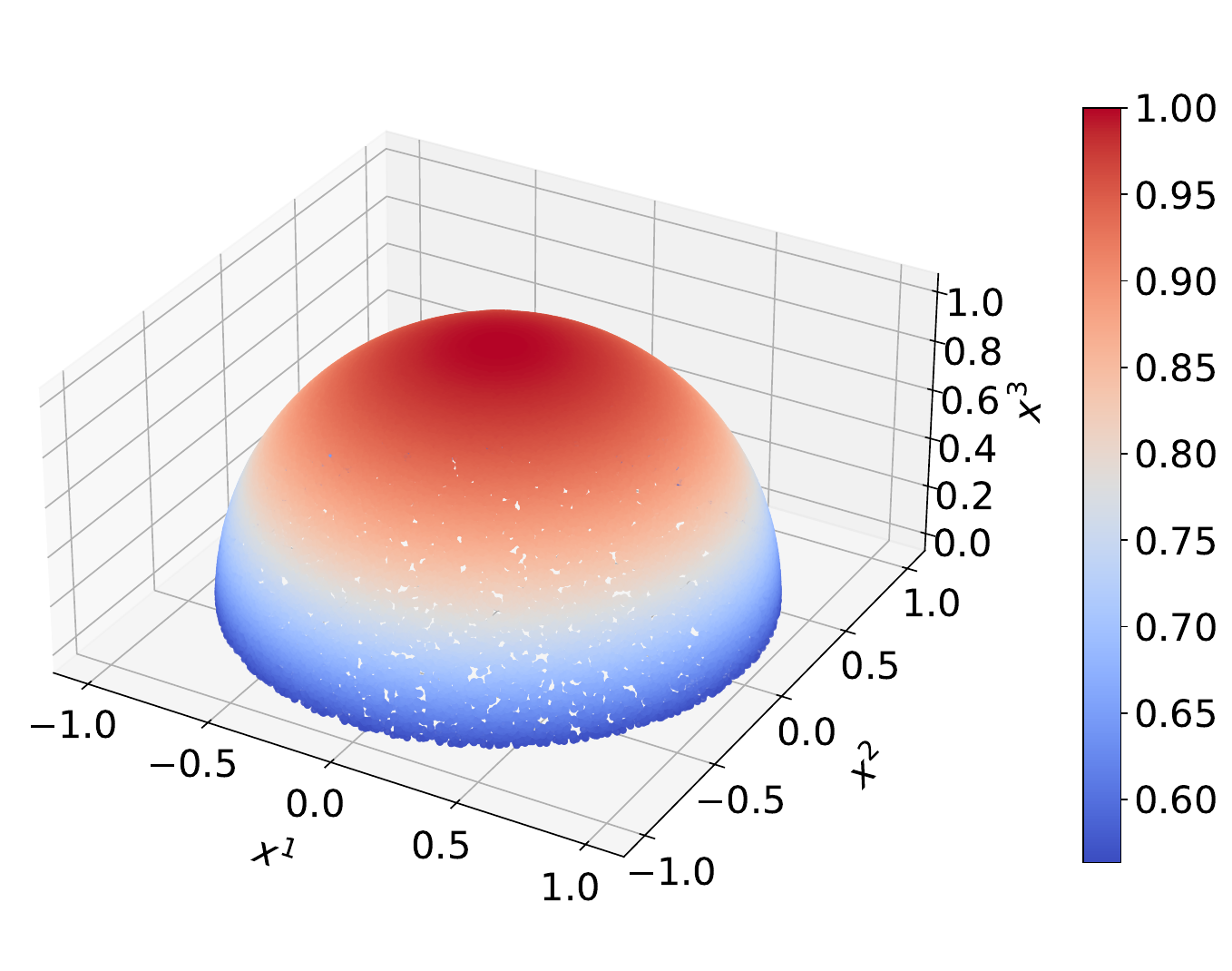}
		\end{subfigure}\hfill
		\begin{subfigure}{0.24\textwidth}
			\centering
			\caption{Mode 4 (Sphere)\\FEM}
			\includegraphics[width=\linewidth]{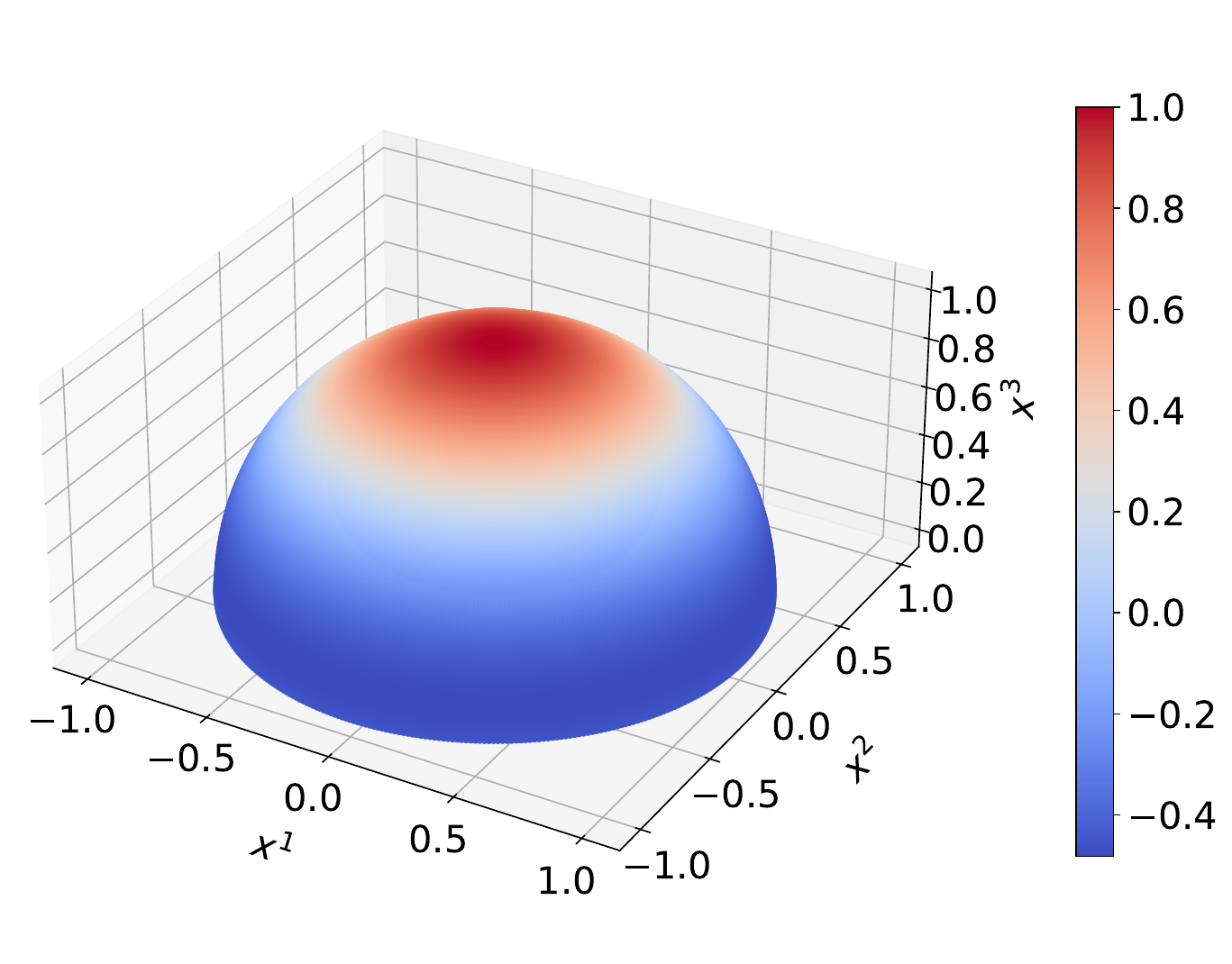}
		\end{subfigure}\hfill
		\begin{subfigure}{0.24\textwidth}
			\centering
			\caption{Mode 4 (Sphere)\\RBF-FD-QP}
			\includegraphics[width=\linewidth]{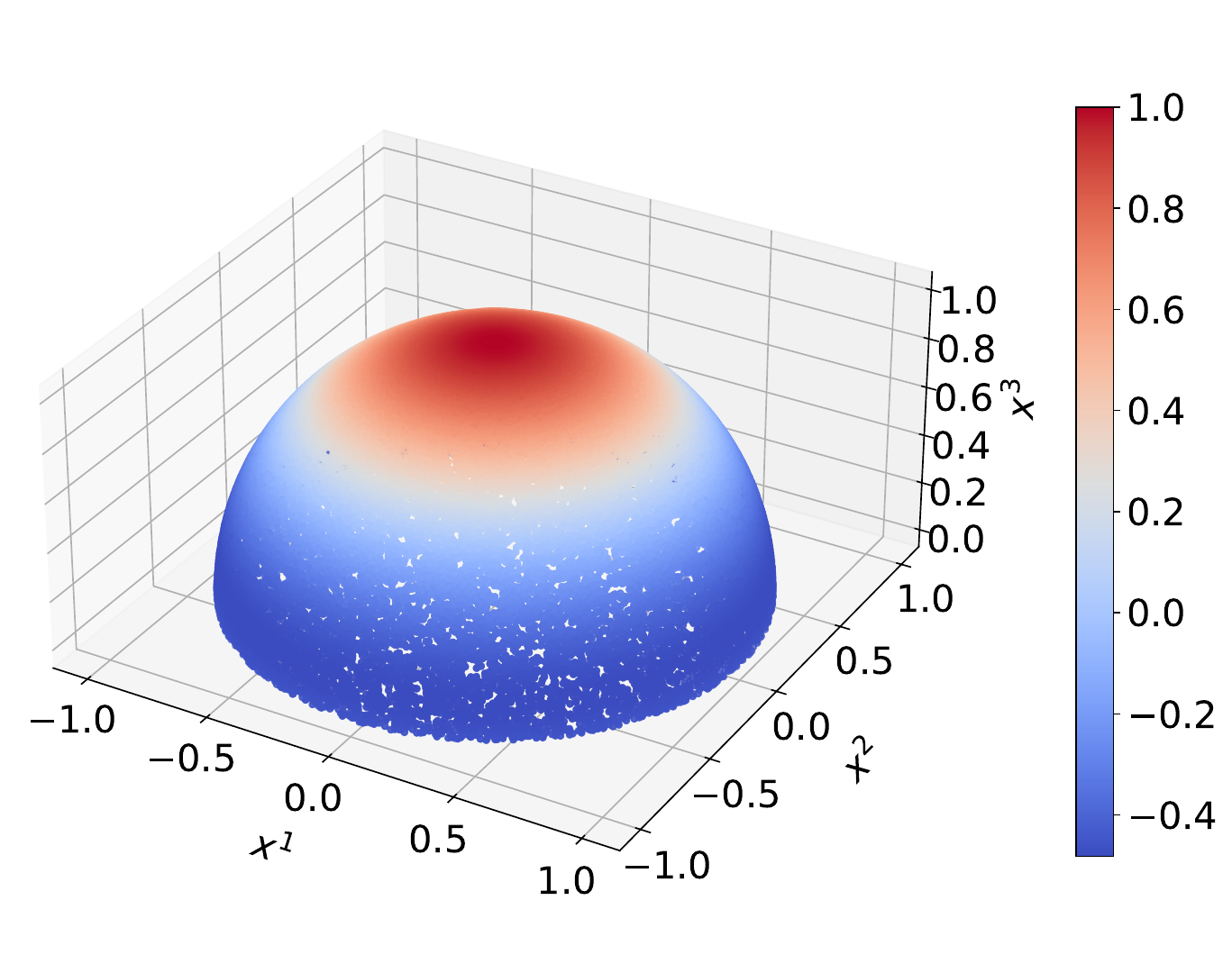}
		\end{subfigure}
		
		\vspace{3ex}
		\begin{subfigure}{0.24\textwidth}
			\centering
			\caption{Mode 1 (Torus)\\FEM}
			\includegraphics[width=\linewidth]{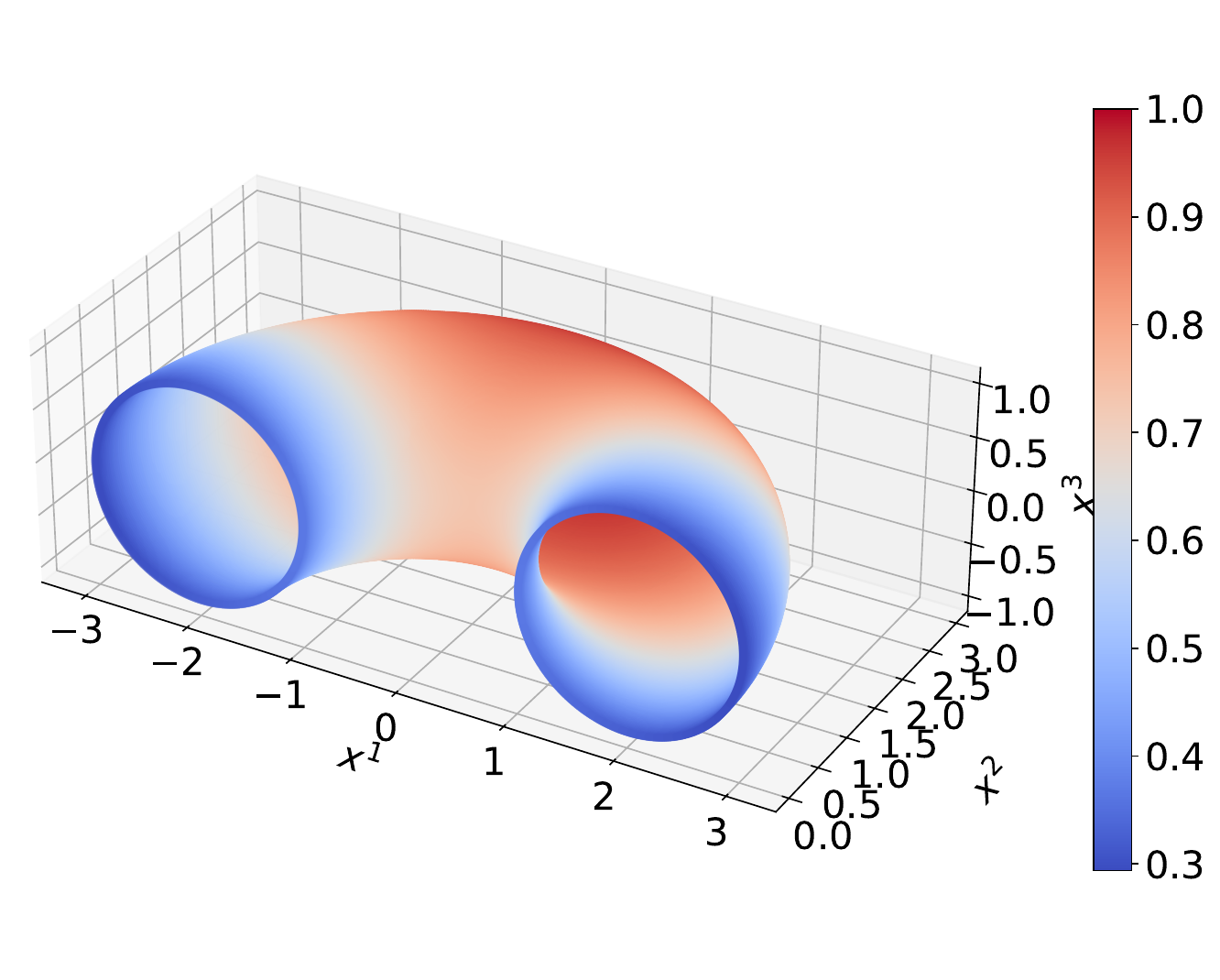}
		\end{subfigure}\hfill
		\begin{subfigure}{0.24\textwidth}
			\centering
			\caption{Mode 1 (Torus)\\RBF-FD-QP}
			\includegraphics[width=\linewidth]{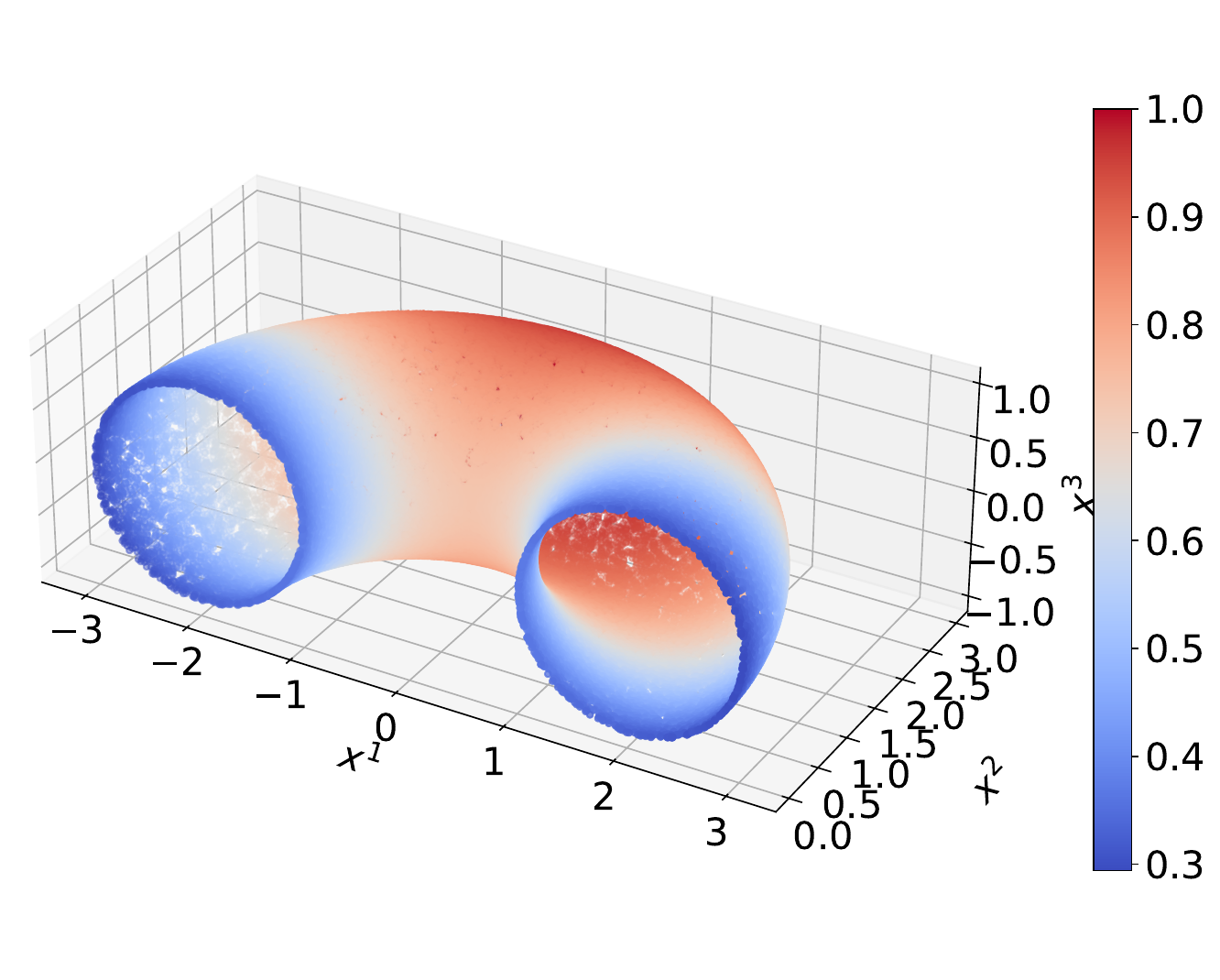}
		\end{subfigure}\hfill
		\begin{subfigure}{0.24\textwidth}
			\centering
			\caption{Mode 8 (Torus)\\FEM}
			\includegraphics[width=\linewidth]{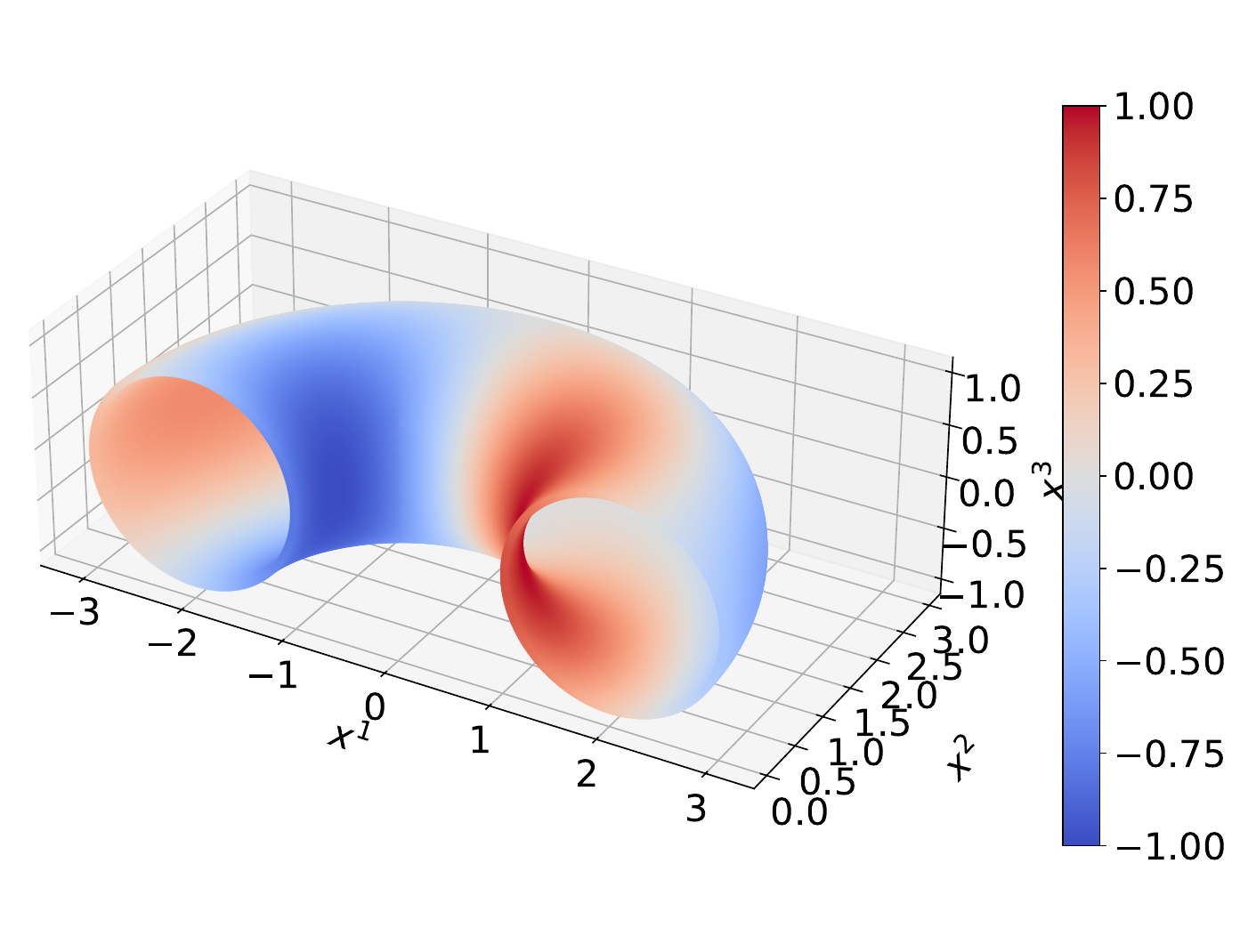}
		\end{subfigure}\hfill
		\begin{subfigure}{0.24\textwidth}
			\centering
			\caption{Mode 8 (Torus)\\RBF-FD-QP}
			\includegraphics[width=\linewidth]{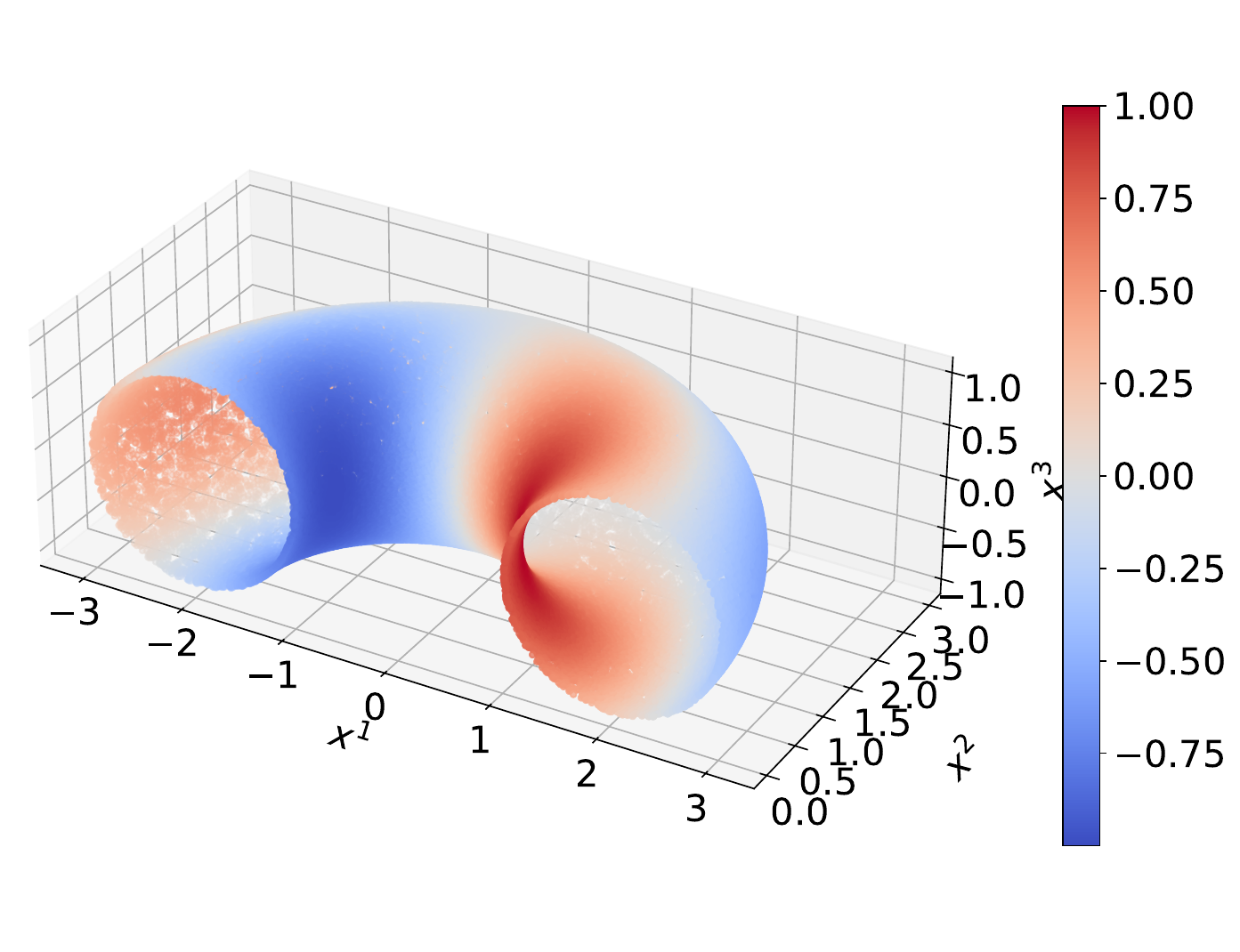}
		\end{subfigure}
		\caption{\textbf{Comparison of normalized eigenfunctions.} Shown are the normalized eigenfunctions calculated by the FEM reference and the proposed RBF-FD-QP approach. Top row: Modes 1 and 4 on the semi-sphere. Bottom row: Modes 1 and 8 on the semi-torus. For RBF-FD-QP, we use $l=4$, $l_{\mathrm{bd}} = 4$ and $N = 51200$.}
		\label{fig:evec_comparison}
	\end{figure*}
	
	In Table~\ref{tab:eigenv_convergence},  we show the mean, standard deviation (std), and convergence rate of the absolute errors for several leading eigenvalues on the two surfaces using 12 independent trials. It can be seen that as the number of nodes increases, the numerical eigenvalues of RBF-FD-QP converge toward the FEM reference values, with most modes exhibiting approximately second-order convergence. The small standard deviations further confirm the stability of the approximation across different random samplings.
		
		Figure~\ref{fig:evec_comparison} compares the eigenfunctions obtained from the proposed RBF-FD-QP method with those from the FEM reference. The numerical eigenfunctions produced by RBF-FD-QP exhibit the same spatial patterns as the FEM results, including the symmetry features, on both the semi-sphere and semi-torus. These results indicate that the proposed method can accurately capture both the eigenvalues and the corresponding eigenfunctions of the Laplace--Beltrami operator on surfaces with boundary.
	
	
	\subsection{Linear time-dependent equations}\label{sec:diffusion}
	
	\begin{figure*}[tbp]
		\centering
		\small
		
		\begin{subfigure}{0.32\textwidth}
			\centering
			\caption{\textbf{IE}}
			\label{fig:time:sphere:ie}
			\includegraphics[height=4.4cm]{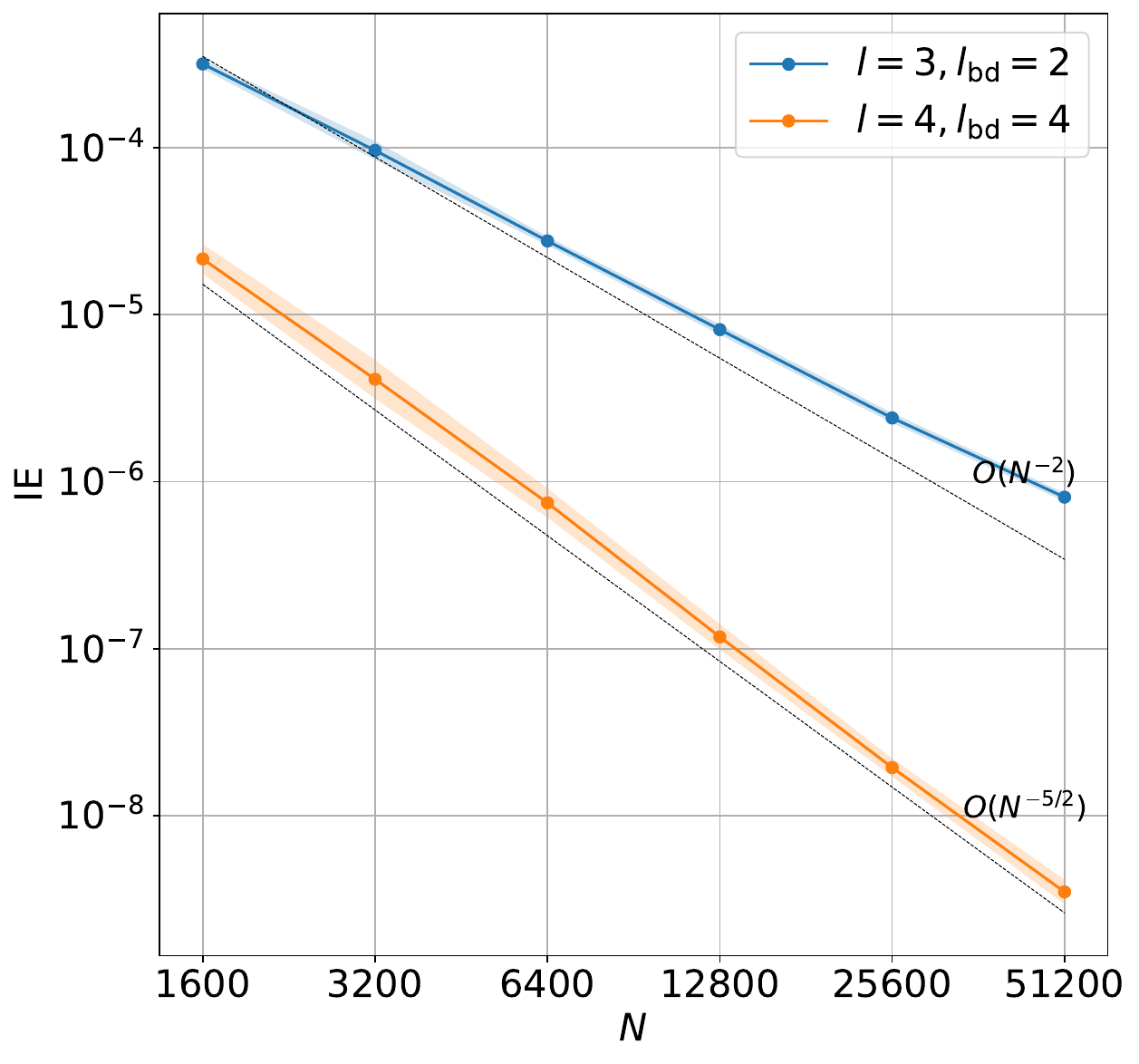}
		\end{subfigure}\hfill
		\begin{subfigure}{0.32\textwidth}
			\centering
			\caption{Numerical solution}
			\label{fig:time:sphere:sol}
			\includegraphics[height=4.4cm]{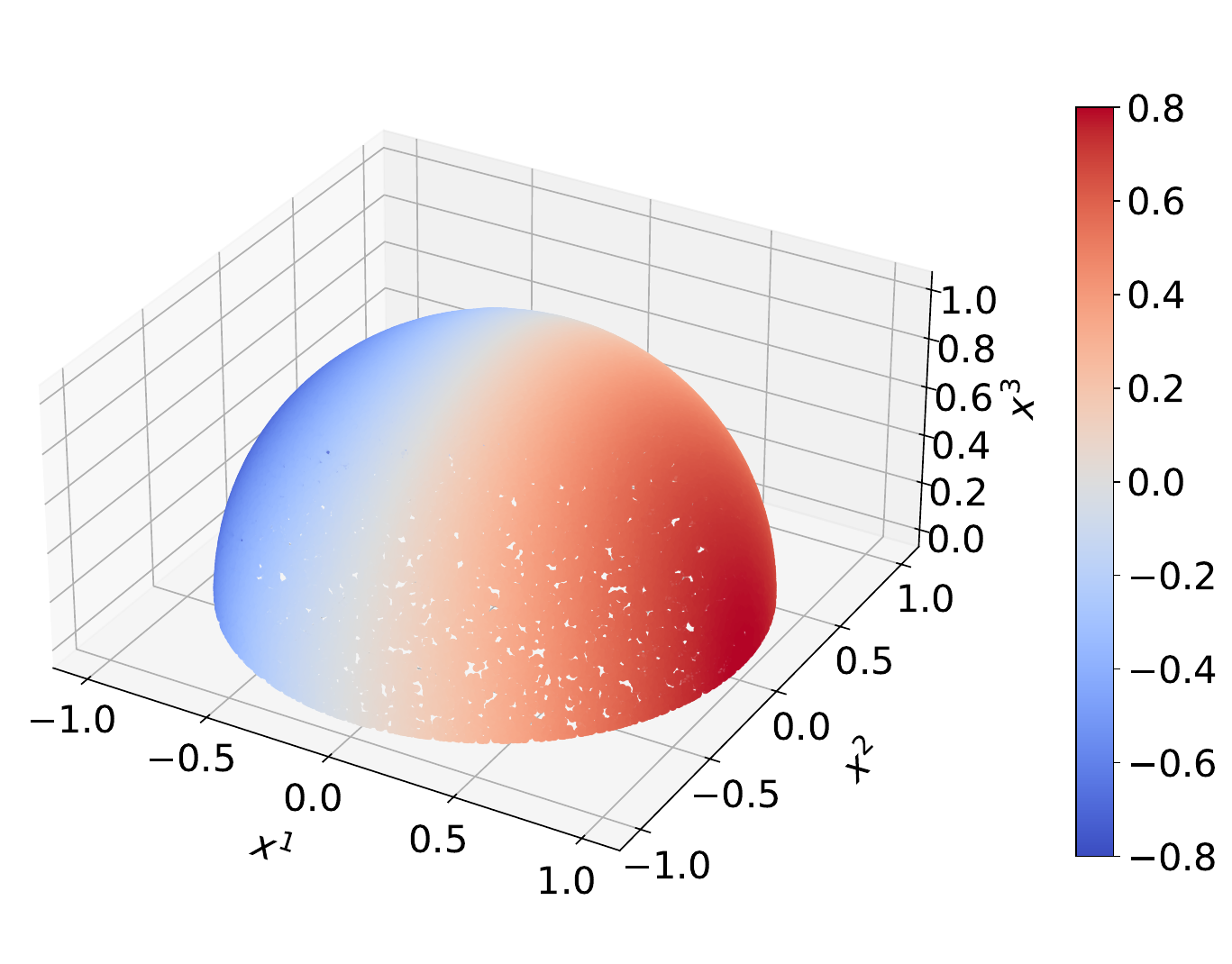}
		\end{subfigure}\hfill
		\begin{subfigure}{0.32\textwidth}
			\centering
			\caption{Pointwise absolute error}
			\label{fig:time:sphere:err}
			\includegraphics[height=4.4cm]{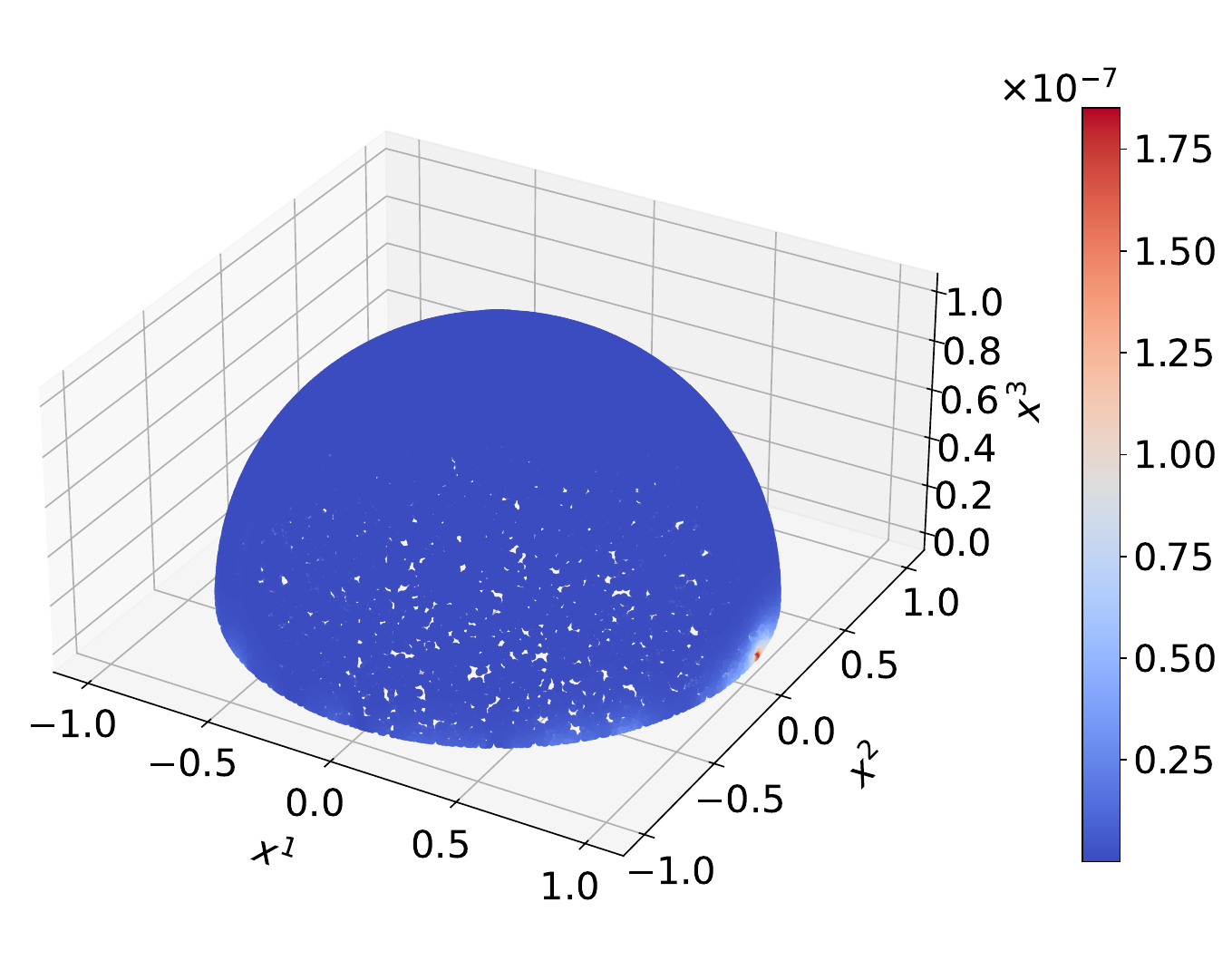}
		\end{subfigure}
		
		\vspace{4ex}
		
		\begin{subfigure}{0.32\textwidth}
			\centering
			\caption{\textbf{IE}}
			\label{fig:time:pipe:ie}
			\includegraphics[height=4.4cm]{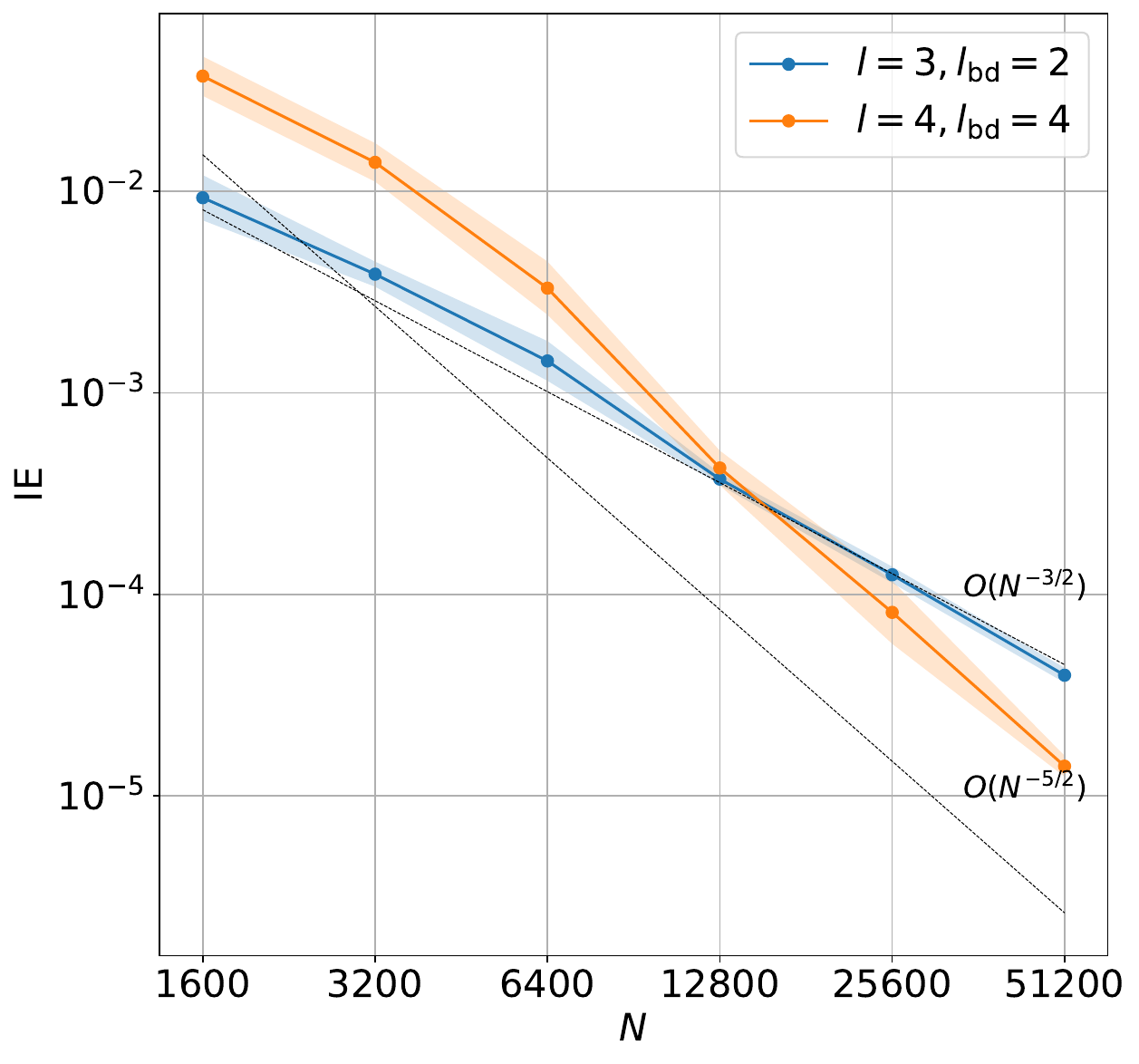}
		\end{subfigure}\hfill
		\begin{subfigure}{0.32\textwidth}
			\centering
			\caption{Numerical solution}
			\label{fig:time:pipe:sol}
			\includegraphics[height=4.4cm]{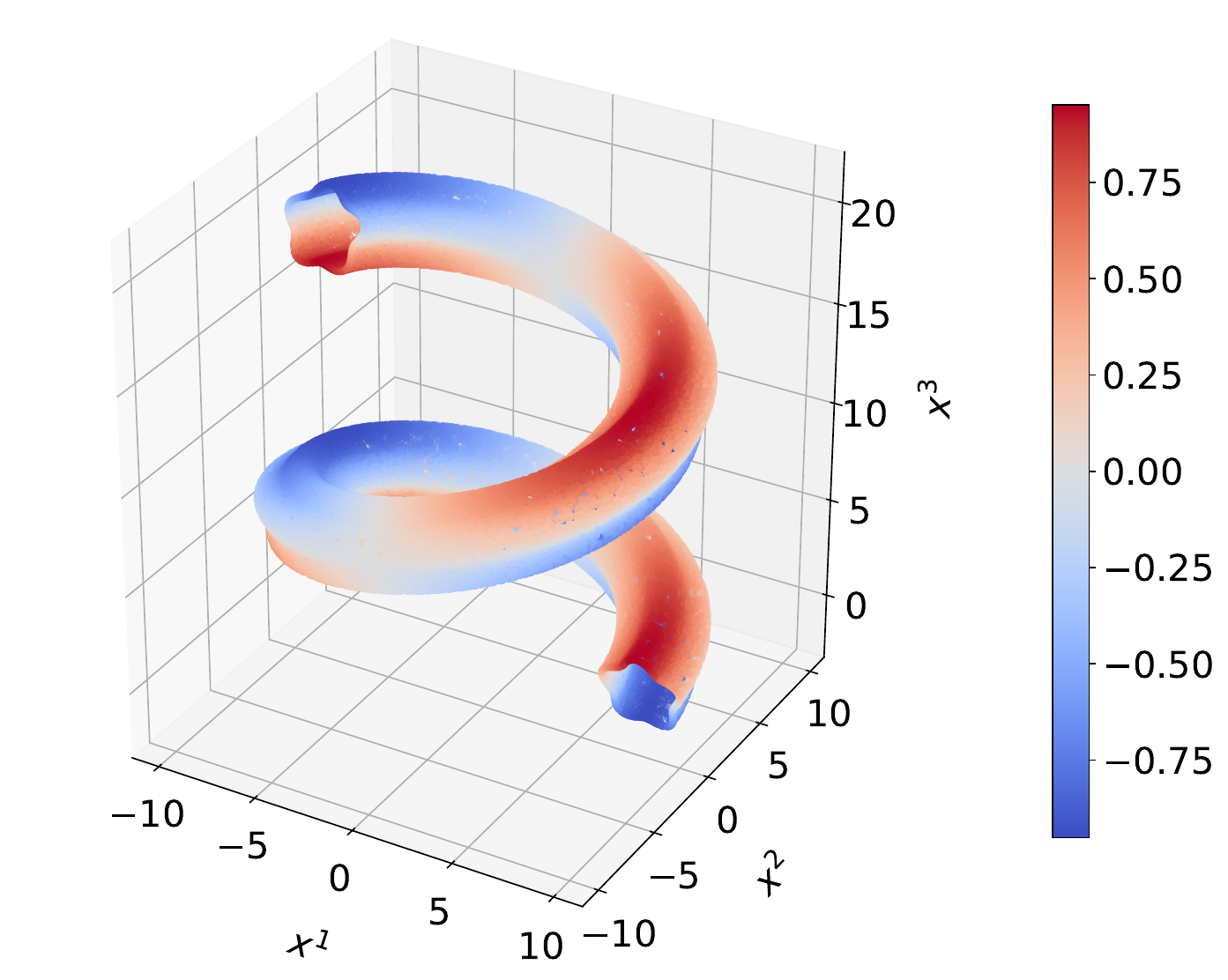}
		\end{subfigure}\hfill
		\begin{subfigure}{0.32\textwidth}
			\centering
			\caption{Pointwise absolute error}
			\label{fig:time:pipe:err}
			\includegraphics[height=4.4cm]{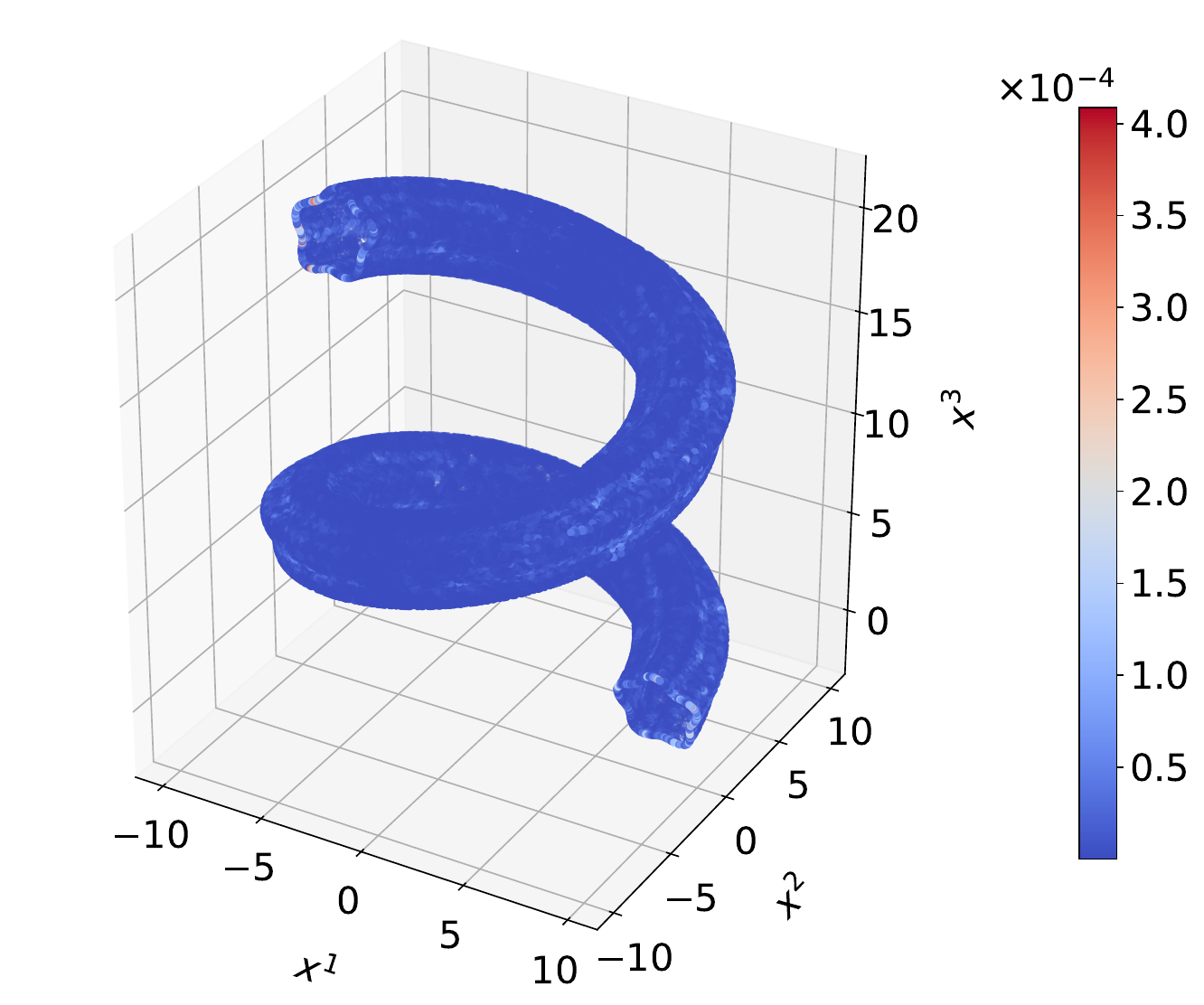}
		\end{subfigure}
		
		\vspace{2ex}
		\caption{\textbf{Numerical results for the time-dependent heat equation on surfaces.} (a)--(c) Simulation on a semi-sphere using the time step $\Delta t = 10^{-3}$ at the final time $T_{\mathrm{end}} = 0.05$. (d)--(f) Simulation on a helical pipe using the time step $\Delta t = 10^{-3}$ at the final time $T_{\mathrm{end}} = 0.05$. The left column shows the convergence of inverse errors (\textbf{IE}s) over various independent trials, examining the stability of the spatial RBF-FD-QP discretization coupled with the implicit BDF4 time integration. The middle and right columns display the numerical solutions and the corresponding pointwise absolute errors, respectively, using  $l=4$, $l_\mathrm{bd} = 4$ and $N=51200$. All simulations are run with 12 independent trials, each with a set of randomly sampled data points.}
		\label{fig:time_evolution_results}
	\end{figure*}
	
	
	Next, we evaluate the performance of the proposed framework on a linear time-dependent diffusion equation defined on a surface,
	\begin{equation}
		u_t = \nu \Delta_M u + f,
	\end{equation}
	where the diffusion coefficient $\nu=0.1$.
	For temporal discretization, we employ the implicit fourth-order backward differentiation formula (BDF4), which allows larger time steps than explicit schemes. In all simulations, we fix the time step with $\Delta t = 10^{-3}$ and run the integration to a final time of $T_{\text{end}} = 0.05$ to examine the convergence.
	
	
	
	For the semi-sphere  defined in above Sec. \ref{sec:eigen}, we set the true solution to be $u(\mathbf{x},t) = e^{-t} \sin(x^1) \cos(x^2)$. For the second example, we consider a geometrically more complex helical pipe featuring a star-shaped cross-section \cite{hu2026pipe}, parameterized by
	\begin{equation}
		\begin{aligned}
			x^1 &= a \cos\psi^2 + R(\psi^1,\psi^2) (\beta \sin\psi^1 \sin\psi^2 - \cos\psi^1 \cos\psi^2), \\
			x^2 &= a \sin\psi^2 - R(\psi^1,\psi^2) (\beta \sin\psi^1 \cos\psi^2 + \cos\psi^1 \sin\psi^2), \\
			x^3 &= b\psi^2 + R(\psi^1,\psi^2) \alpha \sin\psi^1,
		\end{aligned}
	\end{equation}
	where the moving frame constants are defined as $\alpha = \frac{a}{\sqrt{a^2 + b^2}}$ and $\beta = \frac{b}{\sqrt{a^2 + b^2}}$. The helical centerline is defined by the spiral radius $a = 8$ and the vertical growth rate $b = 1$, defined over the domain $\psi^2 \in [0, 3\pi]$. The surface radius for the star-shaped cross-section is explicitly given by:
	\begin{equation}
		R(\psi^1,\psi^2) = \frac{3}{5} + \frac{3}{40} \sin(5\psi^1),
	\end{equation}
	where the first parameter $\psi^1 \in [0,2\pi)$.
	On this surface, we prescribe the exact true solution to be $u(\psi^1, \psi^2, t) = e^{-t} \sin(\psi^1) \cos(\psi^2)$. Then, for both examples, the manufactured $f$ are calculated by $f:= u_t - \nu \Delta_M u $ and the manufactured $h$ are calculated by $h:=\frac{\partial u}{\partial \boldsymbol{n}}+u$ in the Robin boundary conditions  so that the prescribed exact solutions satisfy the governing equations and the boundary conditions.
	

	In Figs.~\ref{fig:time_evolution_results}(a)(d), we  plot the \textbf{IE}s as functions of $N$ at the final time $T_{\text{end}}=0.05$. For both the semi-sphere and the helical pipe, the \textbf{IE}s decrease as the number of nodes increases, indicating convergence of the solutions obtained from the proposed method. We also show the true solutions in Figs. \ref{fig:time_evolution_results}(b)(e), and the pointwise absolute errors in Figs.~\ref{fig:time_evolution_results}(c)(f). The convergence results and pointwise error distributions in Fig.~\ref{fig:time_evolution_results} demonstrate that  the RBF-FD-QP spatial discretization, combined with the BDF4 time-marching scheme, yields  stable and accurate numerical solutions for linear time-dependent diffusion equations on surfaces with boundary, even when the surface geometry is complicated and the surfaces are identified by randomly sampled data points.
	
	
	

	

	\subsection{Interface problem} \label{sec:interface}
	In this subsection,  we consider solving the elliptic interface problems on surfaces:
	\begin{equation}
		-\nabla_M \cdot (c \nabla_M u) + \sigma u = f.
	\end{equation}
	The surface is partitioned into two subdomains, $M^+$ and $M^-$, separated by a smooth interface curve $\tilde{\Lambda}$. Across this interface, the solution satisfies the jump conditions,
	\begin{equation*}
		[u] = u^+ - u^- = q_0,\quad [c \nabla u \cdot \boldsymbol{n}] = c^+ (\nabla u^+ \cdot \boldsymbol{n}) - c^- (\nabla u^- \cdot \boldsymbol{n}) = q_1,
	\end{equation*}
	where $\boldsymbol{n}$ denotes the unit normal to the interface within the tangent plane of the surface.
	While resolving such jumps may involve computationally expensive interface-fitted meshes, our mesh-free RBF-FD-QP framework offers an alternative approach by handling these discontinuities directly from point clouds.
	
	
	
\begin{figure*}[tbp]
		\centering
		\small
		
		\begin{subfigure}{0.32\textwidth}
			\centering
			\caption{\textbf{IE} convergence}
			\label{fig:if:sphere:ie}
			\includegraphics[height=4.4cm]{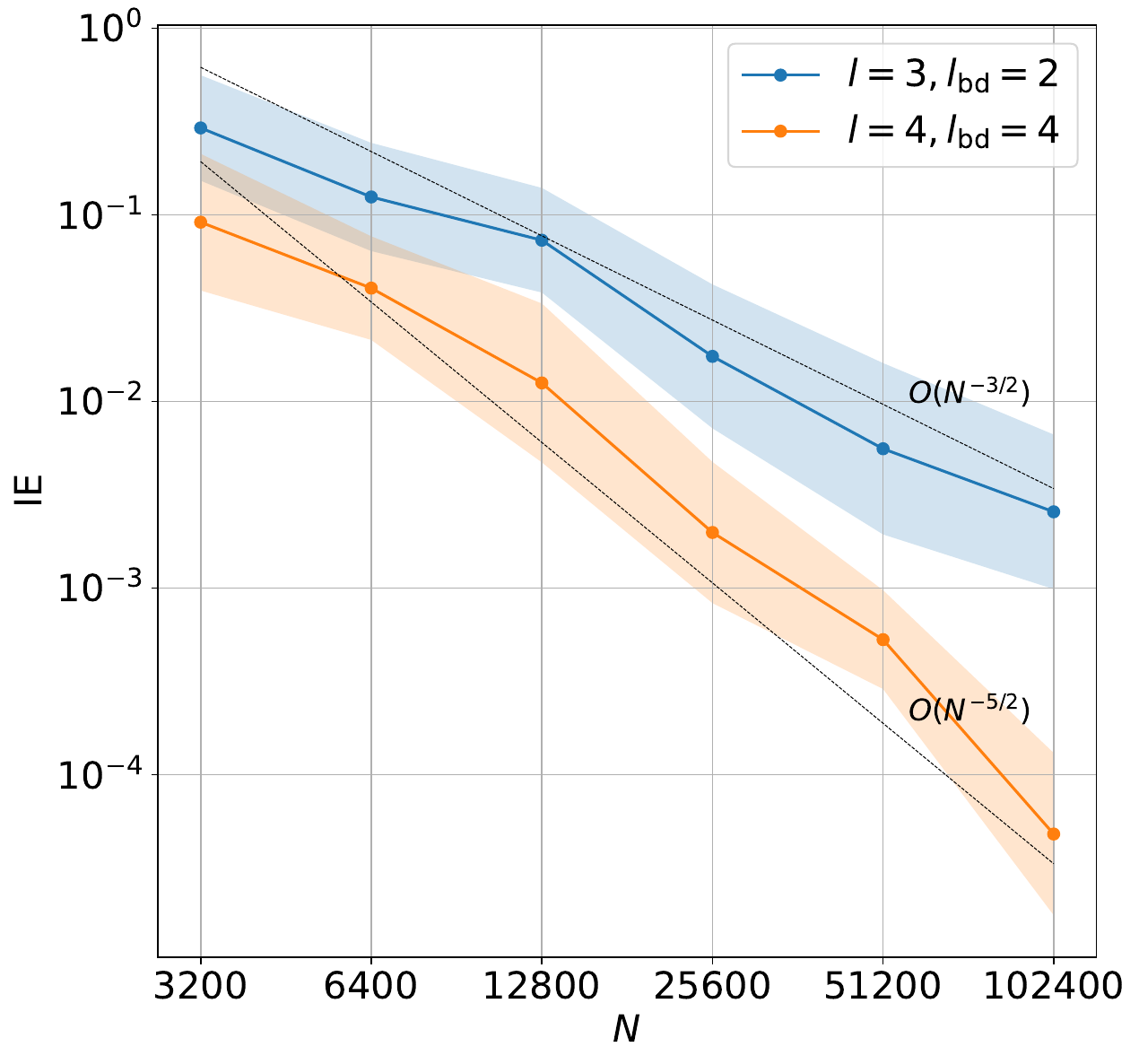}
		\end{subfigure}\hfill
		\begin{subfigure}{0.32\textwidth}
			\centering
			\caption{Numerical solution}
			\label{fig:if:sphere:sol}
			\includegraphics[height=4.4cm]{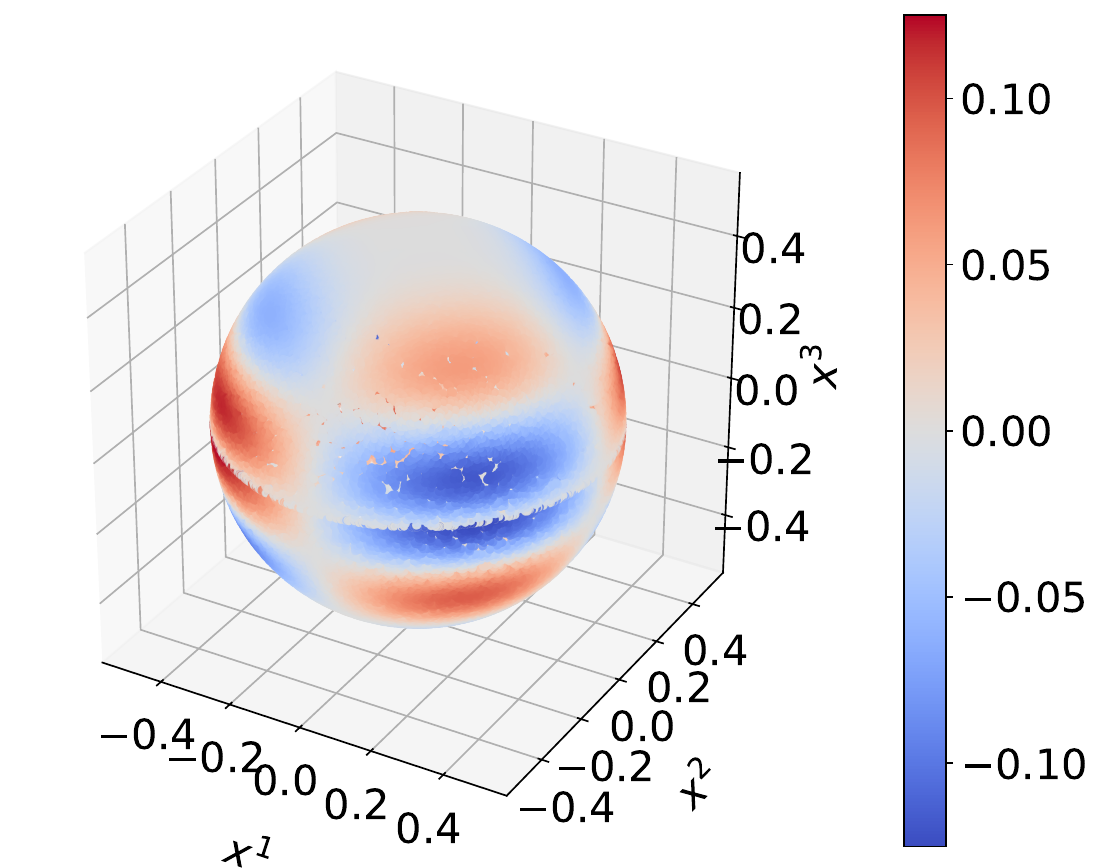}
		\end{subfigure}\hfill
		\begin{subfigure}{0.32\textwidth}
			\centering
			\caption{Pointwise absolute error}
			\label{fig:if:sphere:err}
			\includegraphics[height=4.4cm]{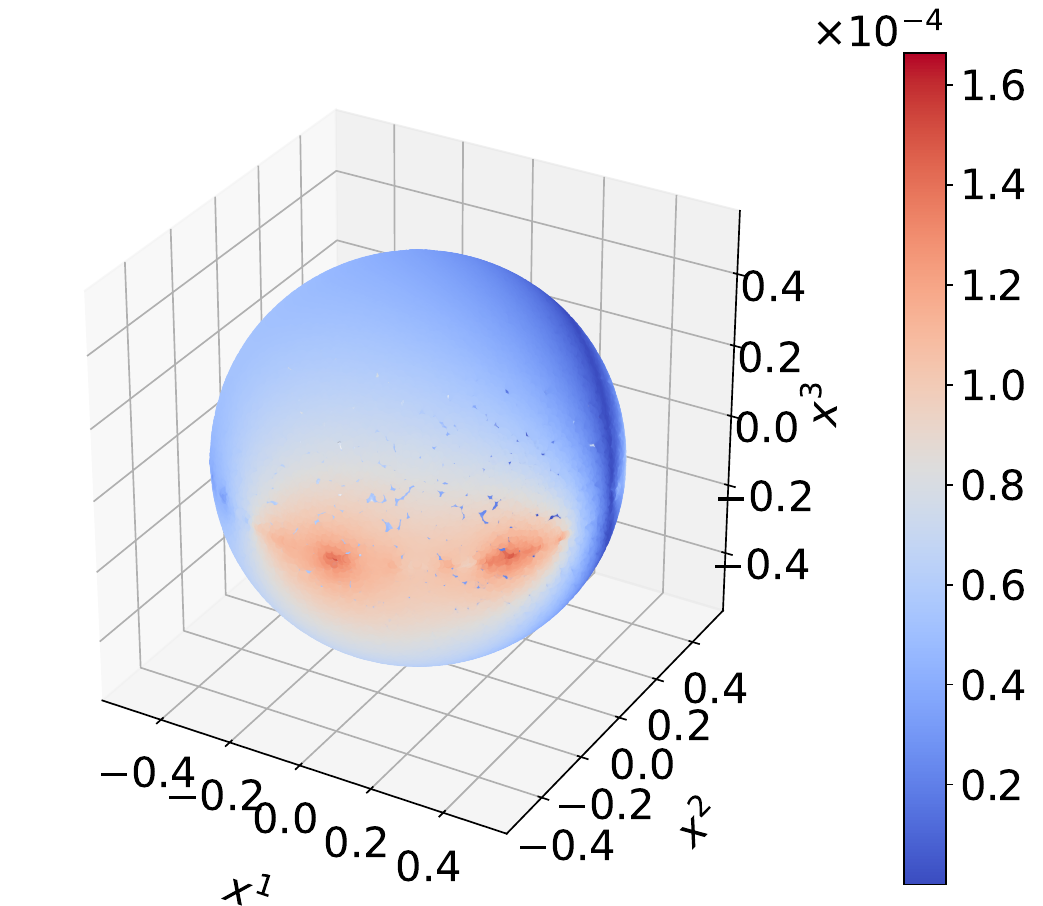}
		\end{subfigure}
		
		\vspace{4ex}
		
		\begin{subfigure}{0.32\textwidth}
			\centering
			\caption{\textbf{IE} convergence}
			\label{fig:if:para:ie}
			\includegraphics[height=4.4cm]{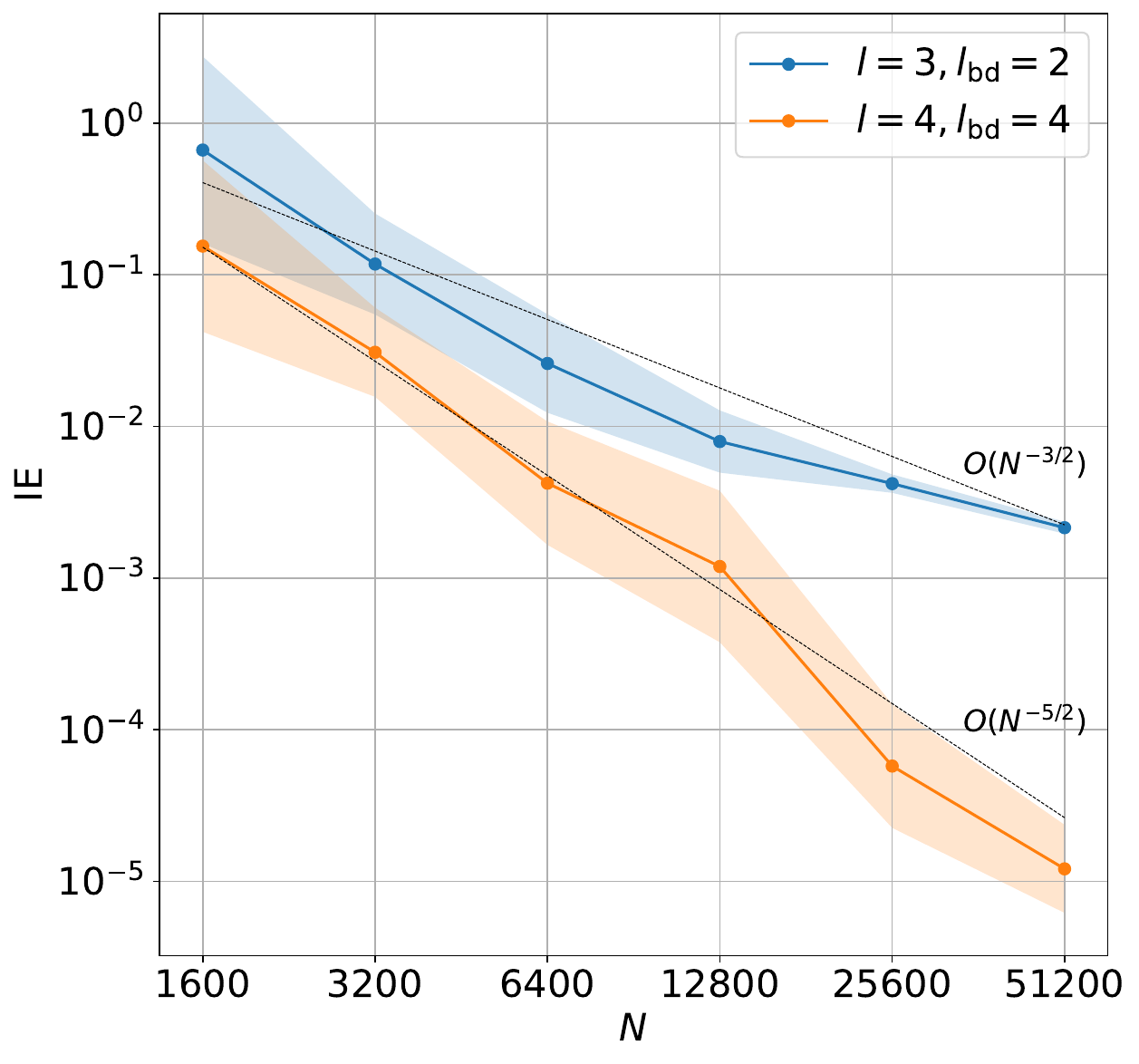}
		\end{subfigure}\hfill
		\begin{subfigure}{0.32\textwidth}
			\centering
			\caption{Numerical solution}
			\label{fig:if:para:sol}
			\includegraphics[height=4.4cm]{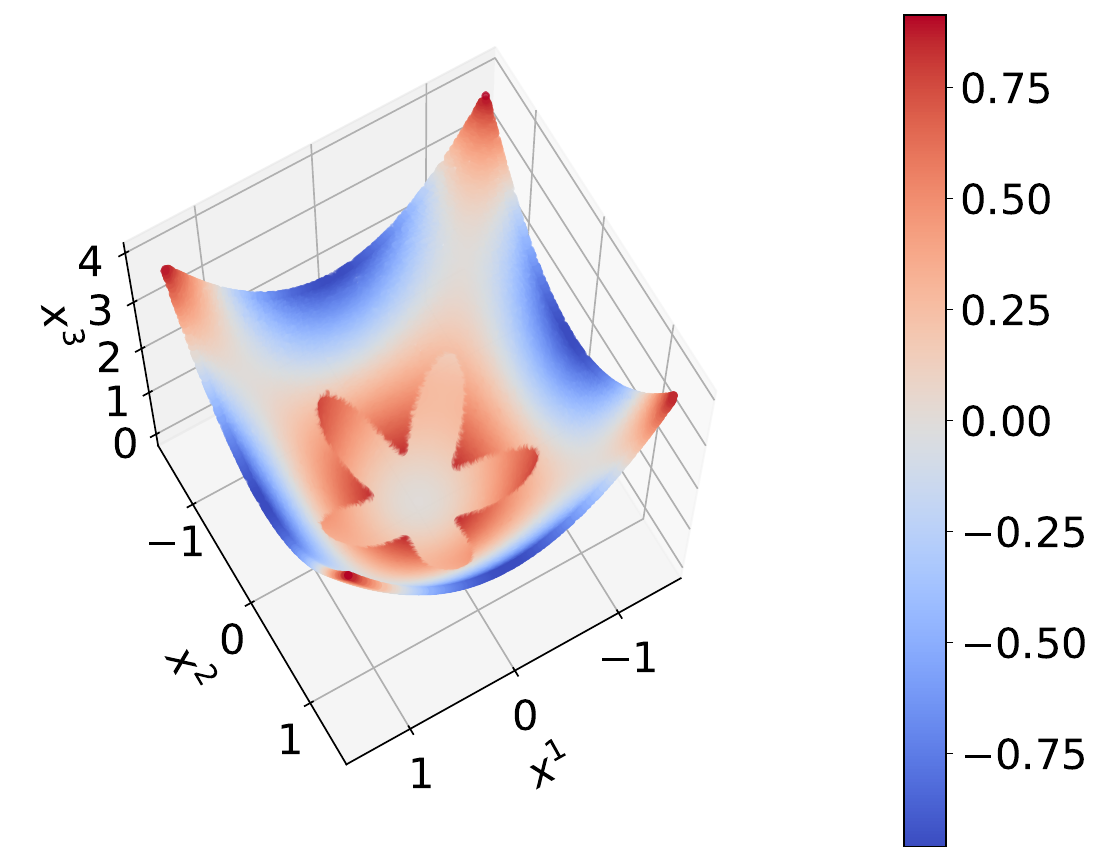}
		\end{subfigure}\hfill
		\begin{subfigure}{0.32\textwidth}
			\centering
			\caption{Pointwise absolute error}
			\label{fig:if:para:err}
			\includegraphics[height=4.4cm]{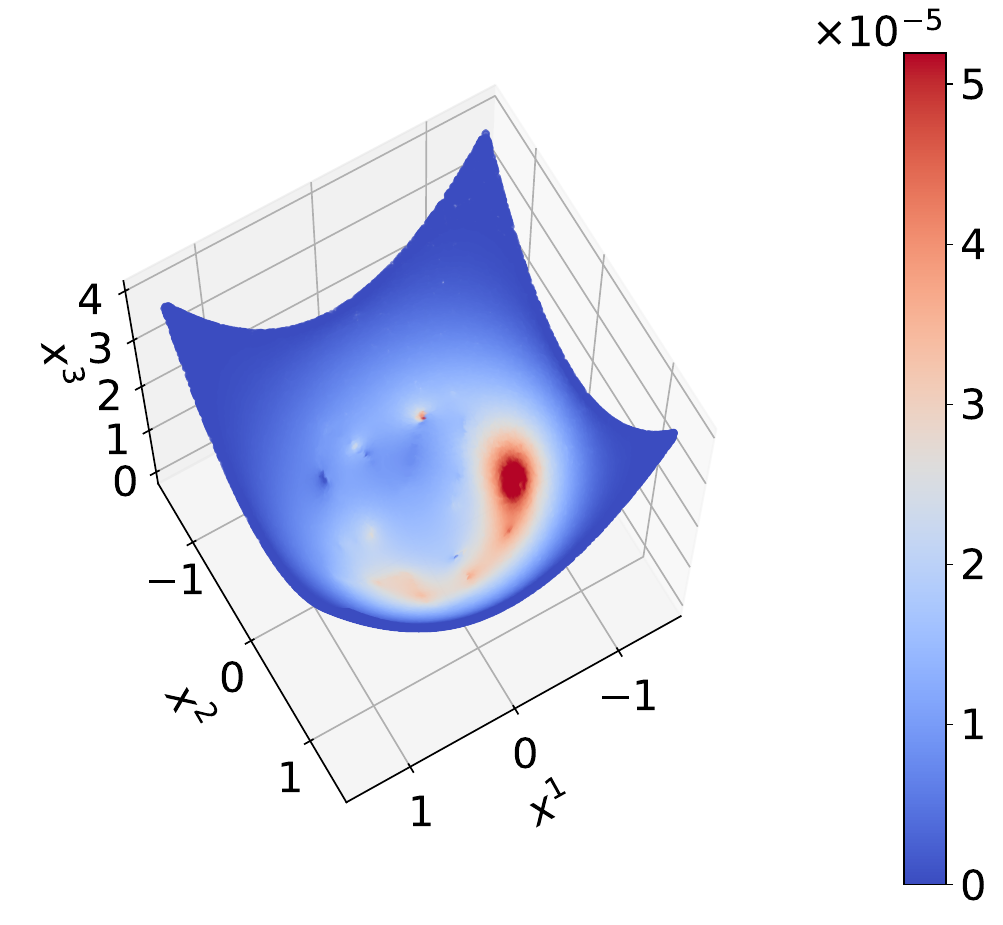}
		\end{subfigure}
		
		\vspace{2ex}
		\caption{\textbf{Numerical results for elliptic interface problems.} (a) A sphere featuring an equatorial jump interface with a 1:10 diffusion contrast. (b) A paraboloid featuring a highly non-convex, five-fold star interface. Left column: Stable high-order convergence of the inverse error (\textbf{IE}). Middle column: The numerical solution, capturing the sharp discontinuities across the interfaces. Right column: Pointwise absolute error distributions. For panels (b), (c), (e) and (f), we use $l=4$, $l_\mathrm{bd} = 4$ and $N=51200$. All simulations are run with 12 independent trials, each with a set of randomly sampled data points.}
		\label{fig:interface_results}
	\end{figure*}

	We test our method on two interface problems, with the numerical results detailed in Fig.~\ref{fig:interface_results}:
	\begin{enumerate}
		\item \textbf{Equatorial Interface on a Sphere \cite{guo2023sphere}:} The spherical domain is parameterized by standard spherical coordinates with a radius of $R = 0.5$. The parameterization mapping $\mathbf{X}: [0, \pi] \times [0, 2\pi) \rightarrow \mathbb{R}^3$ is defined as:
		\begin{equation}
			\begin{cases}
				x^1(\psi^1, \psi^2) = 0.5 \sin\psi^1 \cos\psi^2, \\
				x^2(\psi^1, \psi^2) = 0.5 \sin\psi^1 \sin\psi^2, \\
				x^3(\psi^1, \psi^2) = 0.5 \cos\psi^1,
			\end{cases}
		\end{equation}
		where $\psi^1$ is the polar angle and $\psi^2$ is the azimuthal angle. The domain is partitioned into two subdomains by an equatorial interface, which corresponds to the line $\psi^1 = \pi/2$ in the parameter space (effectively the plane $x^3 = 0$). We impose a high-contrast diffusion ratio ($c^+ = 1.0$, $c^- = 10.0$) alongside a uniform coefficient $\sigma = 1.0$. The exact solution exhibits oscillatory behavior: $u^+ = x^1 x^2 \sin(4\pi x^3)$ and $u^- = x^1 x^2 \cos(4\pi x^3)$.
		
		\item \textbf{Star-shaped Interface on a Paraboloid \cite{yin2025para}:} We embed a highly non-convex, five-fold star interface onto an elliptic paraboloid surface. The underlying surface is parameterized by the mapping $\mathbf{X}: [-1.4, 1.4]^2 \rightarrow \mathbb{R}^3$, where $(\psi^1, \psi^2) \mapsto (\psi^1, \psi^2, (\psi^1)^2 + (\psi^2)^2)$. On this parametric $(\psi^1,\psi^2)$-plane, the interface is defined by a star-shaped domain parameterized by
		\begin{equation}
			\begin{cases}
				\psi^1(t) = [0.7 + 0.3 \cos(5(t - 11\pi/13))]\cos(t), \\
				\psi^2(t) = [0.7 + 0.3 \cos(5(t - 11\pi/13))]\sin(t),
			\end{cases}
		\end{equation}
		with $t \in [0, 2\pi)$.
		The coefficients are set to be $c^+ = 2.0$, $c^- = 1.0$ and $\sigma^+ = 1.0$, $\sigma^- = 3.0$. The exact solution takes different forms across the interface: a rapidly expanding polynomial field $u^+(\mathbf{x}) = ((x^1)^2-1)((x^2)^2-1)$ on one side, and a coupled trigonometric field $u^-(\mathbf{x}) = \cos(x^1+x^2)\sin(x^3)$ on the other.
	\end{enumerate}
	
	Fig. \ref{fig:interface_results} shows convergence of the \textbf{IE}s for both interface configurations. In addition, the numerical solutions sharply resolve the prescribed discontinuities across the interface. The pointwise error plots indicate that the error does not exhibit pronounced deterioration near either the equatorial interface on the sphere or the non-convex star-shaped interface on the paraboloid. These results show that the proposed RBF-FD-QP framework can accurately handle elliptic interface problems on surfaces using scattered nodes only.
	
	Figs.~\ref{fig:interface_results}(a) and \ref{fig:interface_results}(d) show the convergence of the \textbf{IE}s for both interface problems. In addition, our mesh-free method resolves the prescribed discontinuity across the interface.
	As shown in Figs.~\ref{fig:interface_results}(c) and \ref{fig:interface_results}(f),  the pointwise errors remain well-behaved near both the equatorial interface on the sphere and the non-convex star-shaped interface on the paraboloid.
	These results demonstrate that the proposed RBF-FD-QP framework can effectively handle elliptic interface problems on surfaces using only scattered nodes.

	\section{Conclusions}
	\label{sec:conclusion}
	
	In this paper, we considered the two-step RBF-FD method coupled with quadratic programming, referred to as RBF-FD-QP, for solving PDEs on surfaces with derivative boundary conditions, where the surfaces are represented by randomly sampled point clouds. We observed that using the RBF-FD coefficients alone, even when combined with an auto-tuned
	$K$ and a $1/K$ weight function, could still produce unstable approximations of the Laplacian operator and, consequently, non-convergent solutions.
	The instability was found to arise from interior points near the boundary, where stencils are one-sided. To overcome this instability, we introduced a quadratic-programming-based approach to compute weights that yield stable approximations while preserving the consistency constraints of the Laplacian. At boundary points, we applied a similar quadratic optimization formulation to approximate the outward co-normal derivative to further enhance stability. This was accompanied by a specific stencil construction based on restricted $K$-nearest neighbors, selected using a weighted quadratic norm biased in the co-normal direction.
	In addition, we numerically investigated compatible choices for the interior polynomial degree $l$ and the boundary degree $l_{\mathrm{bd}}$. By combining all the techniques described above, we validated the proposed RBF-FD-QP approach and demonstrated numerical convergence of the solutions for various types of surface PDEs with derivative boundary conditions.
	
	There are two questions remaining for future investigation. First, the quadratic-programming-based approach used here is related to the generalized moving least squares (GMLS) method, which only involves polynomial constraints. A natural extension is to incorporate radial basis functions into the quadratic optimization formulation, thereby establishing a connection to the RBF-FD method.
	Second, a promising direction is to apply the proposed RBF-FD-QP approach to variable-coefficient diffusion-type equations, such as weighted diffusion equations and Fokker–Planck equations, with derivative boundary conditions. This is particularly relevant for modeling highly heterogeneous media with discontinuities, as well as the time evolution of distributions for stochastic differential equations with jumps, which mathematically correspond to highly oscillatory or discontinuous coefficients. A key question is whether the proposed RBF-FD-QP approach remains stable and robust for such problems under numerical investigation.
	

	
	\section*{Acknowledgment}
	
	S. J. was supported by the NSFC Grant No. 12471412, the ShanghaiTech University Grant No.
	2024X0303-902-01, and the HPC Platform of ShanghaiTech University. Q. Y. was
	supported by the Simons Foundation grant 601937, DNA.
	
	\bibliographystyle{abbrv}
	\bibliography{arxiv_kme_cleaned}

@article{MR3043640,
  AUTHOR = {Hecht, F.},
  TITLE = {New development in FreeFem++},
  JOURNAL = {J. Numer. Math.},
  FJOURNAL = {Journal of Numerical Mathematics},
  VOLUME = {20}, YEAR = {2012},
  NUMBER = {3-4}, PAGES = {251--265},
  ISSN = {1570-2820},
  MRCLASS = {65Y15},
  MRNUMBER = {3043640},
  URL = {https://freefem.org/}
}

@article{guo2023sphere,
title = {A generalized finite difference method for solving elliptic interface problems with non-homogeneous jump conditions on surfaces},
journal = {Engineering Analysis with Boundary Elements},
volume = {157},
pages = {259-271},
year = {2023},
issn = {0955-7997},
doi = {https://doi.org/10.1016/j.enganabound.2023.09.006},
url = {https://www.sciencedirect.com/science/article/pii/S0955799723004642},
author = {Changyin Guo and Xufeng Xiao and Lina Song and Zhijun Tan and Xinlong Feng},
}

@article{yin2025para,
  title={A kernel-free boundary integral method for elliptic interface problems on surfaces},
  author={Yin, Pengsong and Ying, Wenjun and Zhang, Yulin and Zhou, Han},
  journal={arXiv preprint arXiv:2508.16061},
  year={2025}
}

@article{hu2026pipe,
  author  = {Hu, Shuaifei and Jiao, Yujian and Kong, Desong and Wang, Li-Lian},
  title   = {Solving {PDEs} on Surfaces of Pipe Geometries Using New Coordinate Transformations and High-Order Compact Finite Differences},
  journal = {Journal of Scientific Computing},
  year    = {2026},
  volume  = {107},
  number  = {1},
  pages   = {9},
  doi     = {10.1007/s10915-026-03214-x},
  url     = {https://doi.org/10.1007/s10915-026-03214-x},
}

@article{liu2006stabilized,
  title={A stabilized least-squares radial point collocation method (LS-RPCM) for adaptive analysis},
  author={Liu, GR and Kee, Bernard BT and Chun, Lu},
  journal={Computer methods in applied mechanics and engineering},
  volume={195},
  number={37-40},
  pages={4843--4861},
  year={2006},
  publisher={Elsevier}
}

@article{hangelbroek2024generalized,
  title={Generalized local polynomial reproductions},
  author={Hangelbroek, Thomas and Rieger, Christian and Wright, Grady B},
  journal={arXiv preprint arXiv:2410.12973},
  year={2024}
}

@article{asai2014self,
  title={On self-similar solutions to the surface diffusion flow equations with contact angle boundary conditions},
  author={Asai, Tomoro and Giga, Yoshikazu},
  journal={Interfaces and Free Boundaries},
  volume={16},
  number={4},
  pages={539--573},
  year={2014}
}

@article{lai2010numerical,
  title={Numerical simulation of moving contact lines with surfactant by immersed boundary method},
  author={Lai, Ming-Chih and Tseng, Yu-Hau and Huang, Huaxiong},
  journal={Communications in Computational Physics},
  volume={8},
  number={4},
  pages={735--757},
  year={2010}
}

@article{alphonse2018coupled,
  title={A coupled ligand-receptor bulk-surface system on a moving domain: well posedness, regularity, and convergence to equilibrium},
  author={Alphonse, Amal and Elliott, Charles M and Terra, Joana},
  journal={SIAM Journal on Mathematical Analysis},
  volume={50},
  number={2},
  pages={1544--1592},
  year={2018},
  publisher={SIAM}
}

@article{nielsen2025high,
  title={High-order numerical method for solving elliptic partial differential equations on unfitted node sets},
  author={Nielsen, Morten E and Fornberg, Bengt},
  journal={Computers \& Mathematics with Applications},
  volume={193},
  pages={103--116},
  year={2025},
  publisher={Elsevier}
}

@article{TominecLarssonHeryudono2021,
  author  = {Tominec, Igor and Larsson, Elisabeth and Heryudono, Alfa},
  title   = {A Least Squares Radial Basis Function Finite Difference Method with Improved Stability Properties},
  journal = {SIAM Journal on Scientific Computing},
  volume  = {43},
  number  = {2},
  pages   = {A1441--A1471},
  year    = {2021},
}

@article{li2025two,
  title={Two-step Generalized RBF-Generated Finite Difference Method on Manifolds},
  author={Li, Rongji and Di, Haichuan and Jiang, Shixiao Willing},
  journal={arXiv preprint arXiv:2511.18049},
  year={2025}
}

@article{mohammadi2021divergence,
  title={A divergence-free generalized moving least squares approximation with its application},
  author={Mohammadi, Vahid and Dehghan, Mehdi},
  journal={Applied Numerical Mathematics},
  volume={162},
  pages={374--404},
  year={2021},
  publisher={Elsevier}
}

@article{halada2025overview,
  title={An Overview of Meshfree Collocation Methods},
  author={Halada, Tomas and Yaskovets, Serhii and Singh, Abhinav and Benes, Ludek and Suchde, Pratik and Sbalzarini, Ivo F},
  journal={arXiv preprint arXiv:2509.20056},
  year={2025}
}

@article{barreira2011surface,
  title={The surface finite element method for pattern formation on evolving biological surfaces},
  author={Barreira, Raquel and Elliott, Charles M and Madzvamuse, Anotida},
  journal={J. Math. Biol.},
  volume={63},
  number={6},
  pages={1095--1119},
  year={2011},
  publisher={Springer}
}

@article{dziuk2013finite,
  title={Finite element methods for surface PDEs},
  author={Dziuk, Gerhard and Elliott, Charles M},
  journal={Acta Numerica},
  volume={22},
  pages={289--396},
  year={2013},
  publisher={Cambridge University Press}
}

@article{lin2020radial,
  title={The radial basis function differential quadrature method with ghost points},
  author={Lin, Ji and Zhao, Yuxiang and Watson, Daniel and Chen, CS},
  journal={Mathematics and Computers in Simulation},
  volume={173},
  pages={105--114},
  year={2020},
  publisher={Elsevier}
}

@article{chen2020novel,
  title={A novel RBF collocation method using fictitious centres},
  author={Chen, CS and Karageorghis, Andreas and Dou, Fangfang},
  journal={Appl. Math. Lett.},
  volume={101},
  pages={106069},
  year={2020},
  publisher={Elsevier}
}

@inproceedings{iske2003approximation,
  title={On the approximation order and numerical stability of local Lagrange interpolation by polyharmonic splines},
  author={Iske, Armin},
  booktitle={Modern Developments in Multivariate Approximation: 5th International Conference, Witten-Bommerholz (Germany), September 2002},
  pages={153--165},
  year={2003},
  organization={Springer}
}

@article{lipman2009stable,
  title={Stable moving least-squares},
  author={Lipman, Yaron},
  journal={J. Approx. Theory},
  volume={161},
  number={1},
  pages={371--384},
  year={2009},
  publisher={Elsevier}
}

@inproceedings{liang2012geometric,
  title={Geometric understanding of point clouds using Laplace-Beltrami operator},
  author={Liang, Jian and Lai, Rongjie and Wong, Tsz Wai and Zhao, Hongkai},
  booktitle={2012 IEEE Conference on Computer Vision and Pattern Recognition},
  pages={214--221},
  year={2012},
  organization={IEEE}
}

@article{levin1998approximation,
  title={The approximation power of moving least-squares},
  author={Levin, David},
  journal={Math. Comput.},
  volume={67},
  number={224},
  pages={1517--1531},
  year={1998}
}

@article{petras2018rbf,
  title={An RBF-FD closest point method for solving PDEs on surfaces},
  author={Petras, Argyrios and Ling, Leevan and Ruuth, Steven J},
  journal={J. Comput. Phys.},
  volume={370},
  pages={43--57},
  year={2018},
  publisher={Elsevier}
}

@article{piret2012orthogonal,
  title={The orthogonal gradients method: A radial basis functions method for solving partial differential equations on arbitrary surfaces},
  author={Piret, C{\'e}cile},
  journal={J. Comput. Phys.},
  volume={231},
  number={14},
  pages={4662--4675},
  year={2012},
  publisher={Elsevier}
}

@article{alvarez2021local,
  title={A local radial basis function method for the Laplace--Beltrami operator},
  author={{\'A}lvarez, Diego and Gonz{\'a}lez-Rodr{\'\i}guez, Pedro and Kindelan, Manuel},
  journal={J. Sci. Comput.},
  volume={86},
  year={2021},
  publisher={Springer}
}

@article{suchde2021meshfree,
  title={A meshfree Lagrangian method for flow on manifolds},
  author={Suchde, Pratik},
  journal={Internat. J. Numer. Methods Fluids},
  volume={93},
  number={6},
  pages={1871--1894},
  year={2021},
  publisher={Wiley Online Library}
}

@article{harlim2023radial,
  title={Radial basis approximation of tensor fields on manifolds: from operator estimation to manifold learning},
  author={Harlim, John and Jiang, Shixiao Willing and Peoples, John Wilson},
  journal={J. Mach. Learn. Res.},
  volume={24},
  number={345},
  pages={1--85},
  year={2023}
}

@article{jiang2024generalized,
  title={Generalized finite difference method on unknown manifolds},
  author={Jiang, Shixiao Willing and Li, Rongji and Yan, Qile and Harlim, John},
  journal={J. Comput. Phys.},
  volume={502},
  pages={112812},
  year={2024},
  publisher={Elsevier}
}

@book{pressley2010elementary,
  title={Elementary differential geometry},
  author={Pressley, Andrew N},
  year={2010},
  publisher={Springer Science \& Business Media}
}

@book{monge1809application,
  title={Application de l'analyse {\`a} la g{\'e}om{\'e}trie {\`a} l'usage de l'Ecole imp{\'e}riale polytechnique},
  author={Monge, Gaspard},
  year={1809},
  publisher={Veuve Bernard}
}

@article{mirzaei2012generalized,
  title={On generalized moving least squares and diffuse derivatives},
  author={Mirzaei, Davoud and Schaback, Robert and Dehghan, Mehdi},
  journal={IMA J. Numer. Anal.},
  volume={32},
  number={3},
  pages={983--1000},
  year={2012},
  publisher={OUP}
}

@article{gross2020meshfree,
  title={Meshfree methods on manifolds for hydrodynamic flows on curved surfaces: A Generalized Moving Least-Squares (GMLS) approach},
  author={Gross, Ben J and Trask, Nathaniel and Kuberry, Paul and Atzberger, Paul J},
  journal={J. Comput. Phys.},
  volume={409},
  pages={109340},
  year={2020},
  publisher={Elsevier}
}

@article{jones2023generalized,
  title={Generalized moving least squares vs. radial basis function finite difference methods for approximating surface derivatives},
  author={Jones, Andrew M and Bosler, Peter A and Kuberry, Paul A and Wright, Grady B},
  journal={Comput. Math. Appl.},
  volume={147},
  pages={1--13},
  year={2023},
  publisher={Elsevier}
}

@article{lehto2017radial,
  title={A radial basis function (RBF) compact finite difference (FD) scheme for reaction-diffusion equations on surfaces},
  author={Lehto, Erik and Shankar, Varun and Wright, Grady B},
  journal={SIAM J. Sci. Comput.},
  volume={39},
  number={5},
  pages={A2129--A2151},
  year={2017},
  publisher={SIAM}
}

@article{shankar2015radial,
  title={A radial basis function (RBF)-finite difference (FD) method for diffusion and reaction--diffusion equations on surfaces},
  author={Shankar, Varun and Wright, Grady B and Kirby, Robert M and Fogelson, Aaron L},
  journal={J. Sci. Comput.},
  volume={63},
  number={3},
  pages={745--768},
  year={2015},
  publisher={Springer}
}

@article{flyer2016role,
  title={On the role of polynomials in RBF-FD approximations: I. Interpolation and accuracy},
  author={Flyer, Natasha and Fornberg, Bengt and Bayona, Victor and Barnett, Gregory A},
  journal={J. Comput. Phys.},
  volume={321},
  pages={21--38},
  year={2016},
  publisher={Elsevier}
}

@article{bayona2017role,
  title={On the role of polynomials in RBF-FD approximations: II. Numerical solution of elliptic PDEs},
  author={Bayona, Victor and Flyer, Natasha and Fornberg, Bengt and Barnett, Gregory A},
  journal={J. Comput. Phys.},
  volume={332},
  pages={257--273},
  year={2017},
  publisher={Elsevier}
}

@article{bayona2019role,
  title={On the role of polynomials in RBF-FD approximations: III. Behavior near domain boundaries},
  author={Bayona, V{\'\i}ctor and Flyer, Natasha and Fornberg, Bengt},
  journal={J. Comput. Phys.},
  volume={380},
  pages={378--399},
  year={2019},
  publisher={Elsevier}
}

@article{liang2013solving,
  title={Solving partial differential equations on point clouds},
  author={Liang, Jian and Zhao, Hong-Kai},
  journal={SIAM J. Sci. Comput.},
  volume={35},
  number={3},
  pages={A1461--A1486},
  year={2013},
  publisher={SIAM}
}

@book{Wendland2005Scat,
  title={Scattered Data Approximation},
  author={H. Wendland},
  publisher={Cambridge University Press},
  year={2005},
}

@article{Fuselier2009Stability,
  title={Stability and Error Estimates for Vector Field Interpolation and Decomposition on the Sphere with RBFs},
  author={Fuselier, Edward J  and Wright, Grady  B},
  journal={SIAM J. Numer. Anal.},
  volume={47},
  pages={3213--3239},
  year={2009},
}

@article{fasshauer2012stable,
  title={Stable evaluation of Gaussian radial basis function interpolants},
  author={Fasshauer, Gregory E and McCourt, Michael J},
  journal={SIAM J. Sci. Comput.},
  volume={34},
  number={2},
  pages={A737--A762},
  year={2012},
  publisher={SIAM}
}

@article{jiang2023ghost,
	Author = {Jiang, Shixiao W and Harlim, John},
	journal = {Commun. Pure Appl. Math.,},
  volume = {76},
  number = {2},
  pages = {337-405},
	Title = {Ghost point diffusion maps for solving elliptic PDEs on manifolds with classical boundary conditions},
	Year = {2023}
}

@article{yan2023kernel,
  title={Kernel-based methods for solving time-dependent advection-diffusion equations on manifolds},
  author={Yan, Qile and Jiang, Shixiao W and Harlim, John},
  journal={J. Sci. Comput.},
  volume={94},
  number={1},
  pages={5},
  year={2023},
  publisher={Springer}
}

@article{fuselier2013high,
  title={A high-order kernel method for diffusion and reaction-diffusion equations on surfaces},
  author={Fuselier, Edward J and Wright, Grady B},
  journal={J. Sci. Comput.},
  volume={56},
  number={3},
  pages={535--565},
  year={2013},
  publisher={Springer}
}

@article{flyer2009radial,
  title={A radial basis function method for the shallow water equations on a sphere},
  author={Flyer, Natasha and Wright, Grady B},
  journal={Proc. Roy. Soc. A.},
  volume={465},
  number={2106},
  pages={1949--1976},
  year={2009},
  publisher={The Royal Society London}
}

@book{shaw2019radial,
  title={Radial basis function finite difference approximations of the Laplace-Beltrami operator,  thesis},
  author={Shaw, Sage Byron},
  year={2019},
  publisher={Boise State University Graduate College}
}

@article{wright2023mgm,
  title={MGM: a meshfree geometric multilevel method for systems arising from elliptic equations on point cloud surfaces},
  author={Wright, Grady B and Jones, Andrew and Shankar, Varun},
  journal={SIAM J. Sci. Comput.},
  volume={45},
  number={2},
  pages={A312--A337},
  year={2023},
  publisher={SIAM}
}

@article{gh2019,
title = {Approximating solutions of linear elliptic PDE's on a smooth manifold using local kernel},
journal = "J. Comput. Phys.",
volume = "395",
pages = "563 - 582",
year = "2019",
author = "Faheem Gilani and John Harlim"
}

@article{li2017point,
	title={Point integral method for solving Poisson-type equations on manifolds from point clouds with convergence guarantees},
	author={Li, Zhen and Shi, Zuoqiang and Sun, Jian},
	journal={Commun. Comput. Phys.},
	volume={22},
	number={1},
	pages={228--258},
	year={2017},
	publisher={Cambridge University Press}
}

@article{li2024generalized,
  title={Generalized Moving Least-Squares for Solving Vector-valued PDEs on Unknown Manifolds},
  author={Li, Rongji and Yan, Qile and Jiang, Shixiao W},
  journal={arXiv e-prints},
  pages={arXiv--2406},
  year={2024}
}

@inproceedings{tian2009segmentation,
  title={Segmentation on surfaces with the closest point method},
  author={Tian, Li and Macdonald, Colin B and Ruuth, Steven J},
  booktitle={2009 16th IEEE International Conference on Image Processing (ICIP)},
  pages={3009--3012},
  year={2009},
  organization={IEEE}
}

@article{rauter2018finite,
  title={A finite area scheme for shallow granular flows on three-dimensional surfaces},
  author={Rauter, Matthias and Tukovi{\'c}, {\v{Z}}eljko},
  journal={Comput. \& Fluids},
  volume={166},
  pages={184--199},
  year={2018},
  publisher={Elsevier}
}

@article{elliott2010modeling,
  title={Modeling and computation of two phase geometric biomembranes using surface finite elements},
  author={Elliott, Charles M and Stinner, Bj{\"o}rn},
  journal={J. Comput. Phys.},
  volume={229},
  number={18},
  pages={6585--6612},
  year={2010},
  publisher={Elsevier}
}
	
\end{document}